\documentclass{article}

\usepackage{arxiv}

\usepackage[utf8]{inputenc} 
\usepackage[T1]{fontenc}    
\usepackage{hyperref}       
\usepackage{url}            
\usepackage{booktabs}       
\usepackage{amsfonts}       
\usepackage{nicefrac}       
\usepackage{microtype}      
\usepackage{lipsum}
\usepackage{graphicx}
\graphicspath{ {./images/} }

\usepackage{ wrapfig}

\usepackage{subcaption}
\usepackage{float}

\usepackage{cite}

\usepackage{amsfonts,amsthm}
\usepackage{amssymb,amsmath,latexsym,tensor}
\usepackage{mathrsfs}
\usepackage{enumerate}
\usepackage{savesym}
\usepackage{graphicx}
\usepackage{adjustbox}
\usepackage[pdftex]{color}
\usepackage[font=footnotesize,labelfont=bf]{caption}
\savesymbol{AND}
\usepackage{epstopdf}
\usepackage{algorithm}
\usepackage{algorithmic}

\newcommand{\R}{\mathbb{R}}

\DeclareMathOperator{\E}{\mathbb{E}}

\DeclareMathOperator*{\argmin}{arg\,min}
\DeclareMathOperator*{\argmax}{arg\,max}

\usepackage{mathtools}
\DeclarePairedDelimiterX{\norm}[1]{\lVert}{\rVert}{#1}
\DeclarePairedDelimiterX{\inp}[2]{\langle}{\rangle}{#1, #2}

\newtheorem{theorem}{Theorem}[section]
\newtheorem{corollary}{Corollary}[theorem]
\newtheorem{lemma}[theorem]{Lemma}
\newtheorem{remark}[theorem]{Remark}

\usepackage{xcolor}
\usepackage{mdframed}

\usepackage{array,multirow}

\title{Sampling Kaczmarz Motzkin Method for Linear Feasibility Problems: Generalization \& Acceleration}

\author{
 Md Sarowar Morshed \\
Department of Mechanical $\&$ Industrial Engineering\\
Northeastern University\\
Boston, MA \\
  \texttt{morshed.m@northeastern.edu} \\
   \And
 Md Saiful Islam \\
Department of Mechanical $\&$ Industrial Engineering\\
Northeastern University\\
Boston, MA \\
  \texttt{islam.m@northeastern.edu} \\
  \And
 Md. Noor-E-Alam \\
 Department of Mechanical $\&$ Industrial Engineering\\
Northeastern University\\
Boston, MA \\
  \texttt{mnalam@neu.edu} \\
}

\begin{document}
\maketitle
\begin{abstract}
\textit{Randomized Kaczmarz} (RK), \textit{Motzkin Method} (MM) and \textit{Sampling Kaczmarz Motzkin} (SKM) algorithms are commonly used iterative techniques for solving a system of linear inequalities (i.e., $Ax \leq b$). As linear systems of equations represent a modeling paradigm for solving many optimization problems, these randomized and iterative techniques are gaining popularity among researchers in different domains. In this work, we propose a \textit{Generalized Sampling Kaczmarz Motzkin} (GSKM) method that unifies the iterative methods into a single framework. In addition to the general framework, we propose a Nesterov type acceleration scheme in the SKM method called as \textit{Probably Accelerated Sampling Kaczmarz Motzkin} (PASKM). We prove the convergence theorems for both GSKM and PASKM algorithms in the $L_2$ norm perspective with respect to the proposed sampling distribution. Furthermore, we prove sub-linear convergence for the Cesaro average of iterates for the proposed GSKM and PASKM algorithms.From the convergence theorem of the GSKM algorithm, we find the convergence results of several well-known algorithms like the Kaczmarz method, Motzkin method and SKM algorithm. We perform thorough numerical experiments using both randomly generated and real-world (classification with support vector machine and Netlib LP) test instances to demonstrate the efficiency of the proposed methods. We compare the proposed algorithms with SKM, \textit{Interior Point Method} (IPM) and \textit{Active Set Method} (ASM) in terms of computation time and solution quality. In the majority of the problem instances, the proposed generalized and accelerated algorithms significantly outperform the state-of-the-art methods. 
\end{abstract}


\section{Introduction}
\label{sec:intro}
We consider the following \textit{Linear Feasibility} (LF) problem:
\begin{align}
\label{eq:1}
Ax \leq b, \ \ b \in \R^m, \ A \in \R^{m\times n}.
\end{align}
We confine the scope of our work in the regime of thin/tall coefficient matrix $A$ ($m \gg n$), as iterative methods are more competitive for such problems. Note that, while almost all of the classical methods are deterministic in nature, recent advances \cite{strohmer:2008,lewis:2010,needell:2010,drineas:2011,zouzias:2013,lee:2013,ma:2015,gower:2015,qu:2016,needell:2016,haddock:2017,razaviyayn:2019} suggest that randomized iterative methods can outperform existing deterministic methods for solving many computational problems including linear feasibility, linear systems and convex optimization problems. From an algorithmic point of view, our work generalizes the SKM method and furthermore explores the possibility of faster variants of these methods. Before we delve into the contributions of this work, we give brief descriptions of some of the classical and modern techniques related to solving LF problems with iterative methods.

\paragraph{\textbf{Randomized Kaczmarz (RK)}}
Kaczmarz method is one of the popular methods for solving linear systems due to its algorithmic simplicity \cite{kaczmarz:1937}. Originally proposed in 1937 by Kaczmarz \cite{kaczmarz:1937}, the Kaczmarz method remained hidden to the research community until the early 1980s, when Gordon \textit{et. al} proposed \textit{Algebraic Reconstruction Techniques} (ART) in the area of image reconstruction  \cite{gordon:1970}. Later, it has found applications in several areas like computer tomography \cite{Censor:1988,herman:2009}, digital signal processing \cite{lorenz:2014}, distributed computing \cite{elble:2010,fabio:2012} and many other engineering and physics problems. It has been rediscovered several times as a  family of methods including component solution, successive projection, row-action and cyclic projection methods (see \cite{censor:1981}). Given a current point $x_k$, the Kaczmarz method generates new update $x_{k+1}$ based on the orthogonal projection of $x_k$ onto the hyper-plane $a_{i^{*}}^Tx_k \leq b_{i^{*}}$,
\begin{align}
\label{eq:KM}
  x_{k+1} = x_k- \delta \frac{\left(a_{i^{*}}^Tx_k -b_{i^{*}}\right)^+}{\|a_{i^{*}}\|_2^2} a_{i^{*}}.
\end{align}
The differences between the old and modern Kaczmarz schemes are the choice of projection hyper-planes in the update formula of equation \eqref{eq:KM} at each iteration and the choice of projection parameter $\delta$. The original Kaczmarz method chooses hyper-planes by $i^{*} \equiv k \mod m, k = 1,2,3,...,m$ with parameter $\delta =1$. Strohmer \textit{et. al} \cite{strohmer:2008} showed that instead of using cyclic rules, convergence can be improved by choosing $i^{*}$ from $\{1,2,...,m\}$ at random with probability proportional to $\|a_i^{*}\|_2^2$. This randomization scheme is very efficient for the linear system as well \cite{lewis:2010}. The projection parameter $\delta$ can be chosen any value in the range of $(0,2]$ \cite{haddock:2017}.


\paragraph{\textbf{Motzkin Method (MM)}}
Another classical method for solving LF problems is the Motzkin method (MM) discovered by Motzkin \textit{et. al} in the early 1950s \cite{agamon,motzkin}. The work of Motzkin was rediscovered several times by other researchers in the field of \textit{Machine Learning} (ML). For instance, the so-called perceptron algorithm in ML \cite{rosenblatt,ramdas:2014,ramdas:2016} can be classified as a member of Motzkin type methods. Furthermore, MM can be sought as the Kaczmarz method with ``maximal-residual control" or with ``most violated constraint control"   \cite{censor:1981,nutini:2016,petra:2016}. The MM starts with an initial point $x_k$ and finds the next update $x_{k+1}$ as the projection of $x_k$ onto the most violated hyper-plane defined in the equation \eqref{eq:1}. Given the current point $x_k$, find the next projection hyper-plane $a_{i^*}$ as the maximum violated constraint (i.e., select $i^* = \argmax_{i \in \{1,2,...,m\}} \{a_i^Tx_k-b_i\}$) and then update $x_{k+1}$ as follows
\begin{align}
\label{eq:MM}
  x_{k+1} = (1-\delta) x_k + \delta \ \mathcal{P}_{H_{i^*}}(x_k) ,
  \end{align}
with the choice $ 0 \leq \delta < 2$, where $\mathcal{P}_{H_{i^*}}(x_k)$ denotes the orthogonal projection of $x_k$ onto the hyper-plane $H_{i^*} = a_{i^*}^Tx_k \leq b_{i^*}$. The analysis of the MM depends on the so-called Hoffman constant (see Lemma \ref{lem0} and Table \ref{tab:1}). The main drawback of the standard MM is that it fails to terminate when the LF problem of \eqref{eq:1} is infeasible. In the late 1980s, MM resurfaced for its connection to the ellipsoid method \cite{telgen:1982}. For rational data, it's proven that the system can detect infeasibility and for totally unimodular data, the scheme gives strong polynomial-time algorithms \cite{Maurras1981}. Recently, Chubanov \cite{Chubanov:2012,Chubanov:2015} developed a modified method compared to the traditional relaxation type methods \cite{motzkin}, where instead of projecting on the original hyper-plane, one projects the new point to an induced hyper-plane.

In recent time, Kaczmarz type methods gained immense popularity in the research community. The work of Strohmer \textit{et. al} \cite{strohmer:2008} encouraged numerous extensions and variants of the RK method (see \cite{lewis:2010,needell:2010,zouzias:2013,lee:2013,ma:2015,gower:2015,wright:2016}). For instance, in \cite{zouzias:2013,NEEDELL:2015}, authors analyzed variants of the Kaczmarz method for a least square setup. A significant breakthrough came from the work of Gower \textit{et. al} when they developed a generalized framework namely the \textit{Gower-Richtarik} (GR) sketch. The authors showed that several well-known algorithms like \textit{Randomized Kaczmarz} (RK), \textit{Randomized Newton} (RN) and \textit{Randomized Coordinate Descent} methods can be sought as special cases of the GR algorithm. For different choices of sampling distribution and a positive definite matrix, one can recover all of the above algorithms as special cases (see \cite{gower:2015,anna:2015,needell:2016,hefny:2017} for a detailed discussion). 


Another area of research spurred when Gower \textit{et. al} provided the extension of the GR sketching method to combine several Quasi-Newton methods into one framework \cite{gower2016linearly}. They showed that almost all of the available Quasi-Newton algorithms like \textit{Bad Broyden} (BB), \textit{Powell-Symmetric-Broyden} (PSB), \textit{Good Broyden} (GB), \textit{Broyden–Fletcher–Goldfarb–Shanno} (BFGS) and \textit{Davidon–Fletcher–Powell} (DFP) can be derived as special cases of the GR sketch. In another work, they extended the GR method for finding the pseudo-inverse of a matrix \cite{gower:2017}. Several variants of acceleration have been explored recently for the GR sketch \cite{richtrik2017stochastic,NIPS:2018}. Special block variants of RK methods have been analyzed by Needell \textit{et. al} \cite{NEEDELL:2014,blockneddel:2015,needell2019block}. From a linear programming perspective, Chubanov developed a polynomial-time algorithm for solving the $0-1$ linear system \cite{Chubanov:2012,basu:2014,VEGH:2014} and $0-1$ LF problem \cite{Chubanov:2015}. In recent time, other variants of both RK and SKM algorithms have been developed that deal with various types of sampling strategies \cite{eldar:2011,agaskar:2014,needell:2016,bai:2018,greedbai:2018}.

Moreover, a large number of scientific computing and machine learning tasks aim to solve the unconstrained minimization problem $x^{*} = \argmin \Phi(x)$ with a differentiable function $\Phi: \mathbb{R}^n \rightarrow \mathbb{R}$ \cite{kovachki2019analysis}. \textit{Gradient Descent} (GD) and its variants have been the de facto choice in the artificial intelligence and machine learning community to solve such problems \cite{ruder2016overview}. However, GD suffers from slow convergence as soon as the current solution approaches $x^{*}$. To achieve faster convergence, one of the major algorithmic development is the idea of momentum. The momentum method was first studied by Polyak \cite{polyak1964some} in the sense of rolling a heavy ball along with a well-defined cost function. However, despite its intuitiveness, Polyak's heavy ball momentum was difficult to analyze mathematically. Nesterov's acceleration method, proposed by Nesterov in his seminal work \cite{nesterov:1983} for the GD provides the mathematical rigor that Polyak's method lacks and exhibits the worst-case convergence rate of $O(\frac{1}{k^2})$ for minimizing smooth convex functions compared to the original convergence rate of $O(\frac{1}{k})$.
Since the introduction of Nesterov's work, numerous work has been done on algorithmic development of the first-order accelerated methods (for a detailed discussion see \cite{nesterov:2005,nesterov:2013,nesterov:2014,nesterov:2012}). From then on, Nesterov and Polyak's work has been integrated into several well-known projection-based algorithms like \textit{Coordinate Descent} \cite{nesterov:2012}, \textit{Randomized Kaczmarz} \cite{wright:2016}, \textit{Momentum Induced GR Sketching} \cite{loizou:2017}, \textit{Affine Scaling} \cite{morshed:2018}, \textit{Accelerated Quasi-Newton} \cite{NIPS:2018}, \textit{Randomized Gossip} \cite{peter:2019}, \textit{Sampling Kaczmarz Motzkin} \cite{Morshed2019} and the references therein. Particularly, Morshed \textit{et. al} \cite{Morshed2019} investigated the acceleration scheme of Nesterov in the SKM algorithm for $\delta =1$. 

In this work, we develop a generalized framework namely the GSKM method that extends the SKM algorithm and proves the existence of a family of SKM type methods for solving LF problems. This general framework will provide an ideal platform for the researchers to experiment with a wide range of iterative projection methods and to design efficient algorithms for solving optimization problems in areas like artificial intelligence, machine learning, data mining, and engineering. In addition to the general framework, we propose a Nesterov type acceleration scheme in the SKM method ($0 < \delta < 2$) that outperforms state-of-the-art methods in terms of computation time and solution quality. With the convergence analysis of the GSKM algorithm, we synthesize the convergence analysis of SKM type methods into one convergence theorem from which one can derive convergence results of RK, MM and SKM methods. We also prove convergence of the average iterate (i.e., Cesaro average) generated by both GSKM and PASKM method. We prove sub-linear convergence rate for the Cesaro average under somewhat weaker conditions. We carry out thorough numerical experiments to show the effectiveness of the proposed methods in comparison with state-of-the-art methods for solving a wide range of linear feasibility test instances. Although the proposed methods deal with the case of linear feasibility problem with systems of inequalities, it can be noted that with some modification, like the one stated in the work of Lewis \textit{et. al} \cite{lewis:2010}, one can apply this method to linear systems with both equality and inequality constraints.


The remainder of the paper is organized as follows. The proposed algorithms are discussed in section \ref{sec:contr}, and the convergence analysis of the proposed algorithms is given in section \ref{sec:conv}. In section \ref{sec:num}, we perform extensive numerical experiments on artificial and real test instances for a better understanding of the behavior of the proposed generalized and accelerated schemes. Besides, we compared the effectiveness of the proposed acceleration schemes with state-of-the-art techniques (i.e., SKM, IPM and ASM). And finally, the paper is concluded in section \ref{sec:colc} with concluding remarks and future research directions.

\section{Preliminaries \& Contributions}
\label{sec:contr}
In this section, we discuss the SKM algorithm and some preliminary technical tools to analyze the SKM type methods. We first discuss the notations and assumptions that will be used throughout the paper. We then briefly discuss the SKM method along with the expectation induced by the sampling distribution of the SKM method. To make the analysis easier and more formal, we introduce the function $f(x)$. Finally, we conclude the section with the proposed GSKM method and the PASKM method and their geometric interpretations.

\subsection{Notation}

We follow the standard linear algebra notation in this work. $\mathbb{R}^n$ denotes the $n$ dimensional real space, $\mathbb{R}^{m\times n}$ denotes the set of $m \times n$ real-valued matrices. For any matrix $A\in \mathbb{R}^{m\times n}$, $A^T$ denotes the transpose matrix $A$ and $a_i^T$ for $i = 1,2,..,m$ denotes the rows of matrix $A$. Furthermore, $P = \{x \in \R^n | \ Ax \leq b\}$ denotes the feasible region of the feasibility problem and $\mathcal{P}(x)$ denotes the projection of $x \in \R^n $ onto the feasible region $P$. The notation $d(x,P)$ denotes the distance between $x \in \R^n $ and the feasible region $P$, i.e., $d(x,P) \ = \ \inf_{z \in P} \|x-z\| \ = \ \|x-\mathcal{P}(x)\|$. For any matrix $A$, the spectral norm and Frobenius norm are denoted by $\|A\|$  and $\|A\|_F$, respectively. For any function $f:X\mapsto Y$, we use $\nabla f$ to represent the gradient of $f$. Finally, $\langle x,y\rangle = x^Ty$ denotes the standard inner product and $\|x\| = \sqrt{\langle x, x \rangle}$ as the euclidean ($L_2$) norm. The notation $x^+$ denotes the positive part of any real number (ie., $x^+ = \max \{x,0\}$). For any two arbitrary matrices $M,\ N$, the notation $M \succ N$ implies the positive definiteness of the matrix $M-N$. The notation $\E_\mathbb{S}[\cdot]$ is used to denote the expectation with respect to the sampling distribution $\mathbb{S}$. 
\subsection{Assumptions}
Throughout the paper, we assume that the system $Ax \leq b$ is consistent and the matrix $A$ has no zero rows. We also assumed that the rows of matrix $A$ are normalized (i.e., $\|a_i\|^2=1$ for all $i$). Note that, normalization simplifies the convergence analysis considerably. The normalization doesn't impact the computational time significantly (we could simply normalize each row for the first time it occurs during the computation). Moreover, normalization simplifies the convergence analysis considerably. In the description of algorithms, we do not enforce the assumption. Furthermore, it can be noted that the proposed algorithms generate the same iterates irrespective of normalization.

\subsection{Sampling Kaczmarz Motzkin}

The SKM method (Algorithm \ref{alg:skm}) for solving LF problems, proposed by De Loera \textit{et. al} \cite{haddock:2017}, combines the ideas of both Kaczmarz and Motzkin method. The authors provided a generalized convergence Theorem and a certificate of feasibility which synthesizes the convergence analysis of the Kaczmarz method and Motzkin method for solving LF problems. The proposed method requires only $O(n)$ memory storage and is much more efficient than the state-of-the-art techniques such as Kaczmarz type methods, IPMs and ASMs. The main advantage of SKM can be ascribed to its innovative way of projection plane selection. The hyper-plane selection goes as follows: at iteration $k$ the SKM algorithm selects a collection of $\beta$ rows namely $\tau_k$ uniformly at random out of $m$ rows of the constraint matrix $A$, then out of these $\beta$ rows the row with maximum positive residual is selected (i.e., choose row $i^*$ as $i^* = \argmax_{i \in \tau_k} \{a_i^Tx_k-b_i, 0\}$) and finally the next point $x_{k+1}$ is updated as follows
\begin{align}
\label{eq:SKM}
  x_{k+1} = x_k- \delta \frac{\left(a_{i^*}^Tx_k -b_{i^*}\right)^+}{\|a_{i^*}\|_2^2} a_{i^*}.
\end{align}
For ease of analysis, we denote the above sampling distribution as $\mathbb{S}_k$ at iteration $k$, i.e., at each iteration $k$ choose $\tau_k \sim \mathbb{S}_k$ and denote $i^{*}$ as $i^* = \argmax_{i \in \tau_k \sim \mathbb{S}_k} \left(a_i^Tx_k-b_i\right)^{+}$.

\begin{algorithm}
\caption{SKM Algorithm: $x_{k+1} = \textbf{SKM}(A,b,x_0,K,\delta, \beta )$}
\label{alg:skm}
\begin{algorithmic}
\STATE{Initialize $k \leftarrow 0$;}
\WHILE{$k \leq K$}
\STATE{Choose a sample of $\beta$ constraints, $\tau_k$}, uniformly at random from the rows of matrix $A$.
\STATE{From these $\beta$ constraints, choose $i^* = \argmax_{i \in \tau_k} \{a_i^Tx_k-b_i, 0\}$;}
\STATE{Update $x_{k+1} = x_k- \delta \frac{\left(a_{i^*}^Tx_k -b_{i^*}\right)^+}{\|a_{i^*}\|^2} a_{i^*}$;}
\STATE{$k \leftarrow k+1$;}
\ENDWHILE
\RETURN $x$
\end{algorithmic}
\end{algorithm}

The SKM method generalizes RK and MM, and it also combines their strength in choosing a constraint at each iteration. It has a cheaper per iteration cost compared to Motzkin's method and converges faster compared to the Kaczmarz method. Several extensions of the SKM method in terms of acceleration \cite{Morshed2019}, improved rate \cite{haddock:2019} have been proposed recently.

\subsection{Expectation}
For the convergence analysis of Algorithm \ref{alg:skm} and its variations (any algorithm that uses that specific type of sampling distribution), we need to discuss a specific expectation calculation. First of all, let us sort the residual vector $(Ax-b)^+$ from smallest to largest for any iterate $x$ and denote $(Ax-b)^+_{\underline{\mathbf{i_j}}}$ as the $(\beta+j)^{th}$ entry on the sorted list \footnote{We use the notation $(Ax-b)^+_{\underline{\mathbf{i_j}}}$ throughout the paper to express the underlying expectation, where the indices $\underline{\mathbf{i_j}}$ represent the sampling process of \ref{eq:sampling}.}, i.e.,
\begin{align}
\label{eq:sampling}
  \underbrace{(Ax-b)^+_{\underline{\mathbf{i_0}}}}_{\beta^{th}} \ \leq ... \leq \ \underbrace{(Ax-b)^+_{\underline{\mathbf{i_j}}}}_{(\beta+j)^{th}} \ \leq ... \leq \ \underbrace{(Ax-b)^+_{\underline{\mathbf{i_{m-\beta}}}}}_{m^{th}}.
\end{align}
Now, consider the list with all of the entries of the residual vector $(Ax-b)^+$, then we need to calculate the probability that particular entry of the residual vector is selected at any given iteration. Note that, the probability that any sample is selected is $\frac{1}{\binom{m}{\beta}}$ and each sample has an equal probability of selection. Another intersecting fact can be noted that the size of the residual list controls the order and frequency that each entry of the residual vector will be expected to be selected. From now on, we will denote this specific choice of sampling distribution as $\mathbb{S}$ for any point $x \in \R^n$ \footnote{For ease of notation, throughout the paper, we will use $\mathbb{S}_k$ to denote the sampling distribution corresponding to any random iterate $x_k \in \R^n$ }. To calculate the resulting expectation with respect to the above-mentioned sampling distribution, let us first denote, $\tau \sim \mathbb{S}$ as the set of sampled $\beta$ constraints and $i^{*}$ as \footnote{ Similarly, we will use $\tau_k \sim \mathbb{S}_k$ to denote the sampled set and $i^* = \argmax_{i \in \tau_k \sim \mathbb{S}_k} \{a_i^Tx_k-b_i, 0\} \ = \ \argmax_{i \in \tau_k \sim \mathbb{S}_k} (A_{\tau_k}x_k-b_{\tau_k})^+_i$ for any iterate $x_k \in \R^n$.}
\begin{align}
\label{def:i1}
i^* = \argmax_{i \in \tau \sim \mathbb{S}} \{a_i^Tx-b_i, 0\} \ = \ \argmax_{i} (A_{\tau}x-b_{\tau})^+_i,
\end{align}
where, $A_{\tau}$ denotes the collection of rows of $A$ restricted to the index set $\tau$ and $(A_{\tau}x-b_{\tau})_i$ denotes the $i^{th}$ entry of $A_{\tau}x-b_{\tau}$. Using the above discussion with the list provided in equation \eqref{eq:sampling}, we have the following:
{\allowdisplaybreaks
\begin{align}
\label{def:exp}
\E_{\mathbb{S}} \left[\big |(a_{i^*}^Tx-b_{i^*})^{+}\big |^2\right] = \frac{1}{\binom{m}{\beta}} \sum\limits_{j = 0}^{m-\beta} \binom{\beta-1+j}{\beta-1} \big | (Ax-b)^{+}_{\underline{\mathbf{i_j}}} \big |^2, \end{align}}
where, $\E_{\mathbb{S}}$ denotes the required expectation corresponding to the sampling distribution $\mathbb{S}$. The above expectation calculation was first used by De Loera \textit{et.al} in their work \cite{haddock:2017} where they first introduced the SKM method.

\subsection{Function $f(x)$}
In this section, we formalize the definition of function $f:\R^n \rightarrow \R$. Throughout section \ref{sec:conv}, we will use the properties of function $f(x)$ \footnote{Similar type of functions with uniform sampling have been studied in \cite{needell:2016} \cite{razaviyayn:2019} in the context of stochastic gradient descent and alternating projection algorithms respectively.}. First, for any index $i$, let us define the following function
\begin{align}
\label{def:function_i}
    f_i(x) = \frac{1}{2} |(a_{i}^Tx-b_{i})^+|^2, \quad \nabla f_i(x) = (a_{i}^Tx-b_{i})^+ a_{i}.
\end{align}
Then to simplify the expectation expression of \eqref{def:exp} further, we define the function $f$ and the gradient of $f$ as follows:
\begin{align}
\label{def:function}
    f(x) =  \E_{\mathbb{S}}\left[f_{i^*}(x)\right], \quad \nabla f(x) = \E_{\mathbb{S}} \left[\nabla f_{i^*}(x)\right],
\end{align}
where, the index $i^*$ is selected by the rule provided in \eqref{def:i1}.

\subsection{Contributions}
\noindent \textbf{\textit{Generalized Sampling Kaczmarz Method (GSKM)}}. For obtaining a generalized version of the SKM method, we suggest using history information in updating the current update. In particular, we take two random iterates $x_{k-1}$ and $x_k$ generated by successive SKM iteration and then update the next iterate $x_{k+1}$ as an affine combination of the previous two updates. Starting with $x_0 = x_1 \in \R^n$, for $k\geq 1$, we update
\begin{align*}
    x_{k+1}  = (1-\xi) z_k+ \xi z_{k-1},
\end{align*}
where $z_k = x_k - \delta  \frac{\left(a_{i^*}^Tx_k -b_{i^*}\right)^+}{\|a_{i^*}\|^2} a_{i^*}$ is the $k^{th}$ update of the SKM algorithm. Note that, by taking $\xi = 0$, one can recover the original SKM algorithm. For simpler representation, we denote this method as a generalized SKM method or GSKM method. GSKM method is formally provided in Algorithm \ref{alg:gskm} and the convergence analysis is provided in subsection \ref{subsec:gskm conv}. Our convergence analysis suggests that for any $0 < \delta < 2$, one could choose any $\xi$ such that $\xi \in Q$ \footnote{see \eqref{eq:s}.}.

\begin{algorithm}
\caption{GSKM Algorithm: $x_{k+1} = \textbf{GSKM}(A,b,x_0,K, \delta, \beta, \xi)$}
\label{alg:gskm}
\begin{algorithmic}
\STATE{Choose $0 < \delta < 2, \  \xi \in Q$}\ 
\STATE{Initialize $x_1 =  x_0, \ z_1 =  z_0, \  k = 0$;}
\WHILE{$1 \leq k \leq K$}
\STATE{Choose a sample of $\beta$ constraints, $\tau_k$, uniformly at random from the rows of matrix $A$. From these $\beta$ constraints, choose $i^* = \argmax_{i \in \tau_k} \{a_i^Tx_k-b_i, 0\}$ and update,}
\begin{align}
& z_k  = x_k - \delta  \frac{\left(a_{i^*}^Tx_k -b_{i^*}\right)^+}{\|a_{i^*}\|^2} a_{i^*}; \label{eq:3a} \\
& x_{k+1}  = (1-\xi) z_k+ \xi z_{k-1}; \label{eq:3b}
\end{align}
\STATE{$k \leftarrow k+1$;}
\ENDWHILE
\RETURN $x$
\end{algorithmic}
\end{algorithm}

\begin{table}[h!]
\centering
\caption{Algorithms \& convergence results obtained from GSKM.}
\begin{tabular}{|c|c|c|c|}
\hline
\begin{tabular}[c]{@{}c@{}}Parameters, $\beta, \ \delta, \ \xi $ \end{tabular} &
  \begin{tabular}[c]{@{}c@{}}Row selection  Rule, ($i^{*}$) \end{tabular} &
  Convergence Rate &
  Algorithm \\ \hline
\begin{tabular}[c]{@{}c@{}}$\beta =1, \ \delta = 1, \ \xi = 0$ \end{tabular} &
  \begin{tabular}[c]{@{}c@{}} $\mathbb{P}(i^{*}) = \frac{\|a_i\|^2}{\|A\|^2_F}$ \end{tabular} &
  $\E \left[ r_{k}^2 \right] \leq \left(1-\frac{\lambda_{\min}}{\|A\|^2_F}\right)^k r_{0}^2$ &
  RK \cite{strohmer:2008} \\ \hline
\begin{tabular}[c]{@{}c@{}}$\beta =m, \ \delta = 1, \ \xi = 0$ \end{tabular} &
  \begin{tabular}[c]{@{}c@{}} $i^{*} = \argmax_{j} e_j(x_{k-1}) $ \end{tabular} &
  $ r_{k}^2 \leq \left(1-\frac{\lambda_{\min}}{m}\right)^k r_{0}^2$ &
  MM \cite{motzkin} \\ \hline
\begin{tabular}[c]{@{}c@{}}$0 < \delta < 2, \ \xi = 0$ \end{tabular} &
  \begin{tabular}[c]{@{}c@{}} $\tau_k \sim \mathbb{S}_k$ \\ $i^{*} = \argmax_{j \in \tau_k} e_j(x_{k-1}) $ \end{tabular} &
  $\E \left[ r_{k}^2 \right] \leq \left(1-\frac{\eta}{mL^2}\right)^k r_{0}^2$ &
  SKM \cite{haddock:2017} \\ \hline
\end{tabular}
\label{tab:1}
\end{table}

In Table \ref{tab:1}, we list the algorithms and their respective convergence Theorems recovered from the GSKM algorithm with different parameter choices. To simplify the notation, we denote, $r_{k} =  d(x_k,P), \ \eta = 2\delta - \delta^2, \ \lambda_{\min} = \lambda_{\min}^{+}(A^TA), \ e_j(x) = a_j^Tx-b_j$.

\noindent \textbf{\textit{Probably Accelerated Sampling Kaczmarz Method (PASKM)}}. We propose an accelerated randomized projection method based on the SKM method and Nesterov accelerated gradient (NAG). Note that, NAG generates sequences $\{y_k\}$ and $\{v_k\}$ using the following update formulas:
\begin{align}
\label{eq:nag}
 & y_k  = \alpha_k v_k + (1-\alpha_k) x_k, \quad x_{k+1}  = y_k - \theta_k \nabla f(y_k), \nonumber \\
 & v_{k+1} = \omega_k v_k + (1- \omega_k) y_k - \gamma_k \nabla f (y_k).
\end{align}
In equation \eqref{eq:nag}, $\nabla f$ is the gradient of the given function and $\alpha_k, \omega_k, \theta_k$ are the step sequences. Nesterov used updated values for the sequences $\alpha_k, \omega_k, \theta_k$ and obtained a better convergence rate for the acceleration of standard gradient descent. There are two available works directly involve applying Nesterov's acceleration in Kaczmarz type methods \footnote{Recently, heavy ball momentum method has been proposed in the context of SKM method \cite{morshed:momentum}}, first one is by Wright \textit{et. al} \cite{wright:2016} where the accelerated RK method is proposed for linear systems, the second one deals with applying acceleration in SKM for $\delta =1$ \cite{Morshed2019}. 

\begin{algorithm}
\caption{PASKM Algorithm: $x_{k+1} = \textbf{PASKM}(A,b,x_0,K, \delta, \beta)$}
\label{alg:paskm}
\begin{algorithmic}
\STATE{Initialize $v_0 \leftarrow x_0, \  k \leftarrow 0$;}
\WHILE{$k \leq K$}
\STATE{Choose $\gamma, \omega,  \alpha$ considering either \eqref{cond} or \eqref{eq:askm12} and update
\begin{align}
\label{eq:yk}
  y_k = \alpha v_k + (1-\alpha)x_k;
\end{align}}
\STATE{Choose a sample of $\beta$ constraints, $\tau_k$, uniformly at random from the rows of matrix $A$. From these $\beta$ constraints, choose $i^* = \argmax_{i \in \tau_k} \{a_i^Ty_k-b_i, 0\}$; \ Update 
\begin{align}
x_{k+1} & = y_k- \delta \frac{\left(a_{i^*}^Ty_k -b_{i^*}\right)^+}{\|a_{i^*}\|^2} a_{i^*}; \label{eq:askm1} \\
v_{k+1} & = \omega v_k+ (1-\omega)y_k-\gamma \frac{\left(a_{i^*}^Ty_k -b_{i^*}\right)^+}{\|a_{i^*}\|^2} a_{i^*}; \label{eq:askm2}
\end{align}
}
\STATE{$k \leftarrow k+1$;}
\ENDWHILE
\RETURN $x$
\end{algorithmic}
\end{algorithm}
In this work, we consider the general case $0 < \delta < 2$ and develop a probably accelerated scheme for the SKM algorithm. The main difference between the proposed PASKM algorithm and the above-mentioned method is the choice of step sequences. We propose to use precomputed values for the parameters $\omega, \gamma, \alpha$ for every iterate compared to the iterative parameter selection process in \cite{nesterov:2012,wright:2016,Morshed2019}. Now, using the definition of function $f_i$ (see \eqref{def:function_i}) in \eqref{eq:nag}, we derive the following scheme:
\begin{align}
\label{eq:paskm}
 & y_k = \alpha v_k + (1-\alpha)x_k, \quad  x_{k+1}  = y_k- \delta \nabla f_{i^*}(x),   \nonumber \\
 & v_{k+1}  = \omega v_k+ (1-\omega)y_k-\gamma \nabla f_{i^*}(x),
\end{align}
with $i^{*}$ chosen as $i^{*} = \argmax_{j \in \tau_k} e_j(x_{k-1}) $, where $\tau_k \sim \mathbb{S}_k$. The PASKM method is formalized as Algorithm \ref{alg:paskm} and the detailed convergence analysis of the method is provided in Section \ref{sec:conv}. This method generally outperforms both the SKM and GSKM algorithms for almost all of the test instances considered in this work (see Section \ref{sec:num}).

\subsection{Geometric Interpretation}

The goal of this section is to provide a geometric interpretation of the proposed GSKM and PASKM methods. We shed more lights on how the proposed algorithms work in practice and the difference among SKM, GSKM and PASKM methods. 
\vspace{-10 pt}
\begin{figure}[h!]
\begin{subfigure}{0.48\textwidth}
\includegraphics[width=\linewidth]{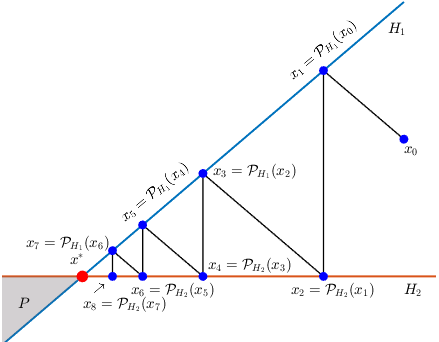}
\caption{SKM: $\delta = 1$}
\end{subfigure}
\hspace*{\fill}
\begin{subfigure}{0.48\textwidth}
\includegraphics[width=\linewidth]{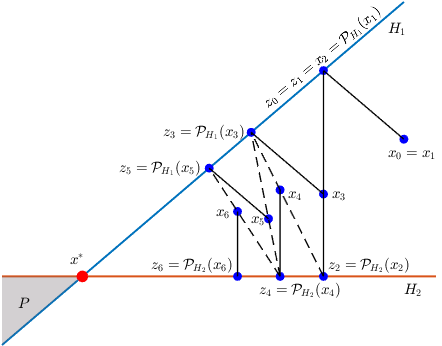}
\caption{GSKM: $\delta = 1, \ \xi = 0.4$}
\end{subfigure}
\caption{Graphical interpretation of the SKM method and the GSKM method with only two hyper-planes $H_j = \{x | a_j^Tx \leq b_j\}$} \label{fig:g1}
\end{figure}

\vspace{-10 pt}
\begin{figure}[h!]
\begin{subfigure}{0.48\textwidth}
\includegraphics[width=\linewidth]{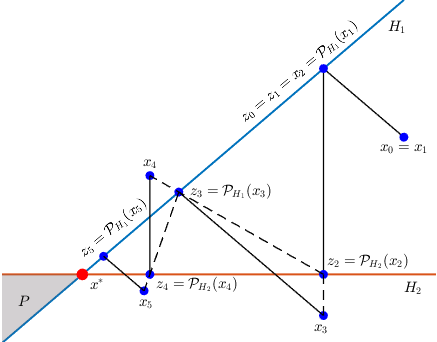}
\caption{GSKM: $\delta = 1, \ \xi = -0.2$}
\end{subfigure}
\hspace*{\fill}
\begin{subfigure}{0.48\textwidth}
\includegraphics[width=\linewidth]{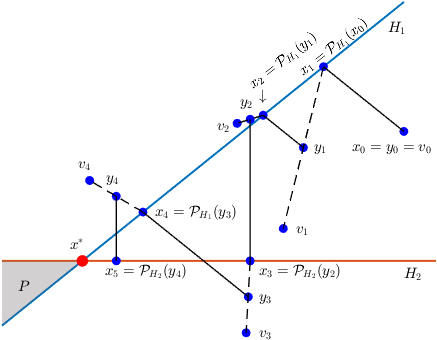}
\caption{PASKM: $\delta = 1, \ \omega = 0.3, \ \alpha = 0.5, \ \gamma = 1.5 $}
\end{subfigure}
\caption{Graphical interpretation of the GSKM method and the PASKM method with only two hyper-planes $H_j = \{x | a_j^Tx \leq b_j\}$} \label{fig:g2}
\end{figure}

In Figures \ref{fig:g1} and \ref{fig:g2}, we illustrate the differences among SKM, GSKM and PASKM methods in an $\R^2$ plane. Our goal is to show how each of the proposed algorithms progress at each iteration. For illustration purposes, We performed the experiment with only two hyper-planes and the selection of hyper-planes is done in an alternative fashion. The notation $\mathcal{P}_{H_1}(x)$ denotes the orthogonal projection of point $x$ onto the hyper-plane $H_1$. 
For comparison purposes, we started with the same starting point $x_0$ and drew the figures with the same scaling. For any given starting point $x_0$, each algorithm projects the point onto the most violated constraint from the sampled constraint set.

The projection step corresponds to the computation of the term $x_k - \delta  \frac{(a_{i^*}^Tx_k -b_{i^*})^+}{\|a_{i^*}\|^2} a_{i^*}$, which means that the current update $x_k$ is projected onto the violated hyper-plane. The projection parameter $\delta \in (0,2]$ defines the type of projection. When $\delta =1$, the projection is exact, that is the point $\mathcal{P}_H(x_k)$ belongs to the hyper-plane $H$. GSKM $(0 \leq \xi \leq 1)$ can be seen as a kind of convex projection update which is slower compared to SKM. From Figure \ref{fig:g2}, it can be seen that the GSKM method with $-1 < \xi <0$ proceeds faster compared to SKM and it requires an affine combination of the previous two successive projections (i.e., $z_{k-1}$ and $z_k$). Compared to SKM and GSKM, the PASKM method updates three different sequences $x_k, v_k, y_k$. From Figure \ref{fig:g2}, it can be noted that GSKM with negative $\xi$ and PASKM moves faster to the feasible region $P$ compared to the SKM method (later in the numerical section this comparison will become much more apparent for larger test instances).

\subsection{Connection between GSKM and PASKM}
\label{conn}
Assume, $-1 < \xi \leq 0$. Then, we can simplify the update formula of the GSKM method as
\begin{align}
\label{connection1}
    x_{k+1} = (1-\xi) x_k + \xi x_{k-1} - \delta (1-\xi)  \nabla f_{i^*}(x_k)- \delta \xi  \nabla f_{j^*}(x_{k-1}),
\end{align}
where the indices $i^*$ and $j^*$ are selected following the rule of \eqref{def:function_i} for the iterate $x_k$ and $x_{k-1}$, respectively. Furthermore, take $\omega (1-\alpha) = -\xi$ and $\gamma$ such that the condition $\alpha \gamma = \delta(1-\xi)$ holds, then from the update formula of the PASKM method we get,
\begin{align*}
    v_{k+1} & \overset{\eqref{eq:askm2}}{=} \omega v_k +(1-\omega) y_k - \gamma \nabla f_{i^*}(y_k) \overset{\eqref{eq:yk}}{=} \ \left(1- \frac{\xi}{\alpha}\right) y_k + \frac{\xi}{\alpha} x_k - \gamma \nabla f_{i^*}(y_k).
\end{align*}
Similarly, from the definition of $y_{k+1}$, we have
\begin{align}
    \label{connection3}
     y_{k+1} &  = \alpha v_{k+1} + (1-\alpha) x_{k+1} = (1-\xi) y_k + \xi x_k - [\alpha \gamma + \delta (1-\alpha)] \nabla f_{i^*}(y_k) \nonumber \\
    & = (1-\xi) y_k  + \xi \left[y_{k-1}-\delta \nabla f_{j^*}(y_{k-1})\right] - \left[\delta(1-\xi)+\delta (1-\alpha)\right] \nabla f_{i^*}(y_k)  \nonumber \\
    & = (1-\xi) y_k + \xi y_{k-1} - \delta (1-\xi)  \nabla f_{i^*}(y_k)- \delta \xi  \nabla f_{j^*}(y_{k-1}),
\end{align}
where the indices $i^*$ and $j^*$ are selected following the rule of \eqref{def:function_i} for the iterate $y_k$ and $y_{k-1}$, respectively. Considering update formulas \eqref{connection1} and \eqref{connection3}, we can conclude that if the conditions $0 \leq \omega (1-\alpha) = -\xi < 1$ and $\alpha \gamma = \delta(1-\xi)$ hold, then the sequence $x_k$ generated by the GSKM algorithm and  the sequence $y_k$ generated by the PASKM algorithm is the same sequence.

\section{Main Results}
\label{sec:conv}

In this section, we present the convergence analysis of the proposed algorithms. In the first subsection, we provided the necessary technical Lemmas $\&$ Theorems that will be used later for our convergence analysis. In the second subsection, we provided the convergence Theorems of the GASKM algorithm. Finally, the last subsection deals with the convergence analysis of the PASKM method.

\subsection{Technical Tools}

In this subsection, we will discuss two types of results. Most of the results derived are related to the properties of the function $f(x)$. Lemma \ref{lem0} is the famous result of Hoffman regarding the linear system of inequalities. Lemmas \ref{lem3}-\ref{lem:grad1} discuss the strong convexity and existence of Lipschitz constant along some restricted segment. Finally, Theorems \ref{th:0} and \ref{th:seq2} deal with developing decay bounds for some non-negative sequences. We will use Lemmas \ref{lem3}-\ref{lem:grad1} frequently in our convergence analysis. Theorems \ref{th:0} and \ref{th:seq2} will be used to derive the proposed convergence bounds of the quantities $\E[d(x_k,P)]$ and $\E[d(x_k,P)^2]$.

\begin{mdframed}[backgroundcolor=gray!15,   topline=false,   bottomline =false,   rightline=false,   leftline=false] \begin{lemma}
\label{lem0}
(Hoffman \cite{hoffman}, Theorem 4.4 in \cite{lewis:2010}) Let $x \in \R^n$ and $P$ be the feasible region, then there exists a constant $L > 0$ such that the following identity holds:
\begin{align*}
  d(x,P)^2 \leq L^2 \ \|(Ax-b)^+\|^2.
\end{align*}
\end{lemma} \end{mdframed}
The constant $L$ is the so-called Hoffman constant. Note that, for a consistent system of equations (i.e., there exists a unique $x^*$ such that $Ax = b$), $L$ can be expressed in terms of the smallest singular value of matrix $A$, i.e.,
\[L^2 = \frac{1}{\|A^{-1}\|^2} = \frac{1}{\lambda_{min}^{+}(A^TA)}.\]

\begin{mdframed}[backgroundcolor=gray!15,   topline=false,   bottomline =false,   rightline=false,   leftline=false] \begin{lemma}
\label{lem:skmseq}
(Lemma 2.1 in \cite{haddock:2017})  Let $\{x_k\},\ \{y_k\}$ be real non-negative sequences such that $x_{k+1} > x_k > 0$ and $y_{k+1} \geq y_k \geq 0$, then
\begin{align*}
  \sum\limits_{k=1}^{n} x_k y_k \ \geq \   \sum\limits_{k=1}^{n} \overline{x} y_k, \quad \text{where} \ \ \overline{x} = \frac{1}{n}\sum\limits_{k=1}^{n} x_k.
\end{align*}
\end{lemma} \end{mdframed}

\begin{mdframed}[backgroundcolor=gray!15,   topline=false,   bottomline =false,   rightline=false,   leftline=false] \begin{lemma}
\label{lem:distance}
For any $x \in \R^n$ and $\bar{x} \in P$, the following identity holds,
\begin{align*}
    d(x,P)^2 \ = \ \|x-\mathcal{P}(x) \|^2 \ \leq \ \|x-\bar{x} \|^2.
\end{align*}
\end{lemma} \end{mdframed}

\begin{mdframed}[backgroundcolor=gray!15,   topline=false,   bottomline =false,   rightline=false,   leftline=false] \begin{lemma}
\label{lem1}
Let $\lambda_j$ be the $j^{th}$ eigenvalue of the matrix $W = \E_{\mathbb{S}}\left[a_{i^*}a_{i^*}^T\right]$, then for all $j$, the bound $0 \leq \lambda_j \leq 1$ holds.
\end{lemma} \end{mdframed}

\begin{proof}
Since $W$ is positive semi-definite, we can write $\lambda_j \geq 0$ for all $j$. Also as the mapping $\lambda_{\max}(X)$ is convex, using Jensen's inequality we have,
\begin{align*}
    \lambda_{\max} (W) = \lambda_{\max} \left[ \E_{\mathbb{S}} \left[a_{i^*}a_{i^*}^T\right] \right] \leq \E_{\mathbb{S}} \left[\lambda_{\max}\left(a_{i^*}a_{i^*}^T\right)\right] \leq 1.
\end{align*}
\end{proof}

\begin{mdframed}[backgroundcolor=gray!15,   topline=false,   bottomline =false,   rightline=false,   leftline=false] \begin{lemma}
\label{lem2}
For any $1 \leq \beta \leq m$, we have the following:
\begin{align*}
 \E_{\mathbb{S}} \left[a_{i^*}a_{i^*}^T\right] \preceq \frac{\beta}{m} A^TA.
\end{align*}
\end{lemma} \end{mdframed}

\begin{proof}
See Appendix 1.
\end{proof}

\begin{mdframed}[backgroundcolor=gray!15,   topline=false,   bottomline =false,   rightline=false,   leftline=false] \begin{lemma}
\label{lem3}
For any $x \in \R^n$ with $\lambda_{\max} = \lambda_{\max}(A^TA)$, we have the following:
\begin{align*}
    \frac{\mu_1}{ 2} \ d(x,P)^2 \ \leq \ f(x) \ \leq \  \frac{\mu_2}{2}\ d(x,P)^2,
\end{align*}
with $ 0 < \mu_1 = \frac{1}{m L^2} \leq \ \mu_2 = \min \left\{1, \frac{\beta}{m} \lambda_{\max}\right\} \leq 1$. 
\end{lemma} \end{mdframed}

\begin{proof}
See Appendix 1.
\end{proof}

Lemma \ref{lem3} states that the function $f$ is strongly convex with constant $\mu_1 $ and has Lipschitz continuous gradient with constant $\mu_2$ when restricted along the segment $[x,\mathcal{P}(x)]$. Let, $f^* = \min_{x} f(x)$, then it can be easily checked that $f^{*} = f(x^*)= 0$. Here, $x^*$ is the optimal solution and it satisfies $Ax^* \leq b$. Moreover, the point $\mathcal{P}(x)$ satisfies the condition $\nabla f(\mathcal{P}(x)) = 0 $. Then we rewrite the inequalities of Lemma \ref{lem3} as follows
\begin{align}
    & \frac{\mu_1}{ 2} \|x-\mathcal{P}(x)\|^2 + \langle \nabla f(\mathcal{P}(x)),x-\mathcal{P}(x) \rangle  \ \leq \ f(x) - f^*, \label{eq:strongl}\\
    & f(x) - f^* \ \leq \  \langle \nabla f(\mathcal{P}(x)),x-\mathcal{P}(x) \rangle + \frac{\mu_2}{2}\ \|x-\mathcal{P}(x)\|^2. \label{eq:strongm}
\end{align}
Here, equation \eqref{eq:strongl} and \eqref{eq:strongm} represent the Lipschitz continuity condition and the strong convexity condition respectively along the line segment $[x,\mathcal{P}(x)]$. For our convergence analysis of Algorithm \ref{alg:gskm} and \ref{alg:paskm}, we will need inequalities like \eqref{eq:strongl} and \eqref{eq:strongm} along the segment  $[x,y]$ for any $x,y \in \R^n$. Following two Lemmas deal with the problem of finding such bounds.

\begin{mdframed}[backgroundcolor=gray!15,   topline=false,   bottomline =false,   rightline=false,   leftline=false] \begin{lemma}
\label{lem:grad}
For any $x, y \in \R^n$, we have the following:
\begin{align*}
 \langle  x-y, \E_{\mathbb{S}}  \left[(a_{i^*}^Ty-b_{i^*})^{+} a_{i^*}\right]  \rangle  & =  \langle  x-y, \nabla f(y)   \rangle   \\
 & \leq    f(x) - f(y) \leq     \frac{\mu_2}{2} \ d(x,P)^2 -  \frac{\mu_1}{2} \ d(y,P)^2.
\end{align*}
\end{lemma} \end{mdframed}

\begin{proof}
See Appendix 1.
\end{proof}

\begin{remark}
We note that the condition of Lemma \ref{lem:grad} is weaker than the traditional strong convexity, and it is also weaker than the essentially strong convexity condition defined in \cite{karimi:2016}. For instance, the essentially strong convexity requires the following identity:
\begin{align*}
    f(x)-f(y)  \leq \langle \nabla f(x), x-y \rangle - \frac{\epsilon}{2} \|x-y\|^2, \quad \forall x,y, \ \text{s.t.} \ \mathcal{P}(x) = \mathcal{P}(y),  
\end{align*}
for some $\epsilon>0$. The above condition clearly implies \eqref{eq:grad}. Moreover, the restricted secant inequality condition defined in \cite{karimi:2016} can be written as
\begin{align}
\label{eq:rsi}
    \langle \nabla f(x), x-\mathcal{P}(x) \rangle \geq \epsilon \|x-\mathcal{P}(x)\|^2.
\end{align}
Note that, with the choice $x = \mathcal{P}(y)$ in Lemma \ref{lem:grad}, we have the following:
\begin{align*}
    \langle \nabla f(y), y-\mathcal{P}(y) \rangle \geq \frac{\mu_2}{2} \|y-\mathcal{P}(y)\|^2.
\end{align*}
Here, we used the fact $d(\mathcal{P}(y),P)^2 = \|\mathcal{P}(y)-\mathcal{P}(y)\|^2 = 0$. This implies that the function $f(x)$ satisfies the restricted secant inequality condition of \eqref{eq:rsi} with constant $\epsilon = \frac{\mu_1}{2}$. Indeed it can be shown that the constant $\epsilon = \frac{\mu_1}{2}$ can be improved further (see the following Lemma). 
\end{remark}

\begin{mdframed}[backgroundcolor=gray!15,   topline=false,   bottomline =false,   rightline=false,   leftline=false] \begin{lemma}
\label{lem:grad1}
For any $y \in \R^n$ and $\bar{y}$ such that $A \bar{y} \leq b$, we have the following:
\begin{align*}
\langle  \bar{y}-y, \E_{\mathbb{S}} \left[a_{i^*}(a_{i^*}^Ty-b_{i^*})^{+}\right]  \rangle = \ \langle  \bar{y}-y, \nabla f(y)  \rangle    \leq  -2f(y)  \leq \ - \mu_1  d(y,P)^2.
\end{align*}
\end{lemma} \end{mdframed}

\begin{proof}
See Appendix 1.
\end{proof}

\begin{remark}
Substituting $ \bar{y} = \mathcal{P}(y)$, in Lemma \ref{lem:grad1} we have, 
\begin{align*}
\langle  \mathcal{P}(y)-y, \E_{\mathbb{S}} & \left[a_{i^*}(a_{i^*}^Ty-b_{i^*})^{+}\right]  \rangle \ \leq  -2f(y)  \leq  - \mu_1 \ d(y,P)^2.
\end{align*}
Note that, similar types of results can be found in the literature. For instance, in \cite{razaviyayn:2019}, authors obtained similar result with respect to a different expectation, they used $\E[x] = \frac{1}{n}\sum \nolimits_i x_i$ for any $x \in \R^n$, which is commonly used to analyze randomized Kaczmarz type methods (see \cite{strohmer:2008,lewis:2010}). Furthermore, we believe a better upper bound than the one obtained in Lemma \ref{lem:grad} can be obtained considering some restrictions on the data matrix $A$. To that end, one needs to obtain a better version of equation \eqref{eq:grad}, i.e., one needs to show that the function $f(x)-\frac{\epsilon}{2} \|x\|^2 $ is convex along the line segment $[x,y]$.
\end{remark}

\begin{mdframed}[backgroundcolor=gray!15,   topline=false,   bottomline =false,   rightline=false,   leftline=false] \begin{lemma}
\label{lem4}
For any $x \in \R^n$ and $0 < \delta < 2$, we have the following:
\begin{align*}
  \E_{\mathbb{S}} \left[ d(z,P)^2\right] = \E_{\mathbb{S}} \left[ \Big \| x- \mathcal{P} (x) - \delta \left(a_{i^*}^Tx-b_{i^*}\right)^+ a_{i^*}\Big \|^2 \right]  \leq   h(\delta) \ d(x,P)^2,
\end{align*}
where, $z = x - \delta \left(a_{i^*}^Tx-b_{i^*}\right)^+ a_{i^*}$, $ \eta = 2 \delta -\delta^2$ and $h(\delta) = 1- \eta \mu_1 < 1$.
\end{lemma} \end{mdframed}

\begin{proof}
See Appendix 1.
\end{proof}

Before we delved into the main Theorems, for any $ \phi_1, \phi_2 \geq 0$, let us define the following parameters:
\begin{align}
\label{def:0}
   &  \phi = \frac{-\phi_1 +\sqrt{\phi_1^2 +4\phi_2 }}{2},  \ \ \rho = \phi + \phi_1, \nonumber \\
   & R_1 = \frac{1+\phi}{\phi+\rho}, \ R_2 = \frac{1-\rho}{\phi+\rho}, \ R_3 = \frac{\rho+\phi_2}{\phi+\rho}, \ R_4 = \frac{\phi-\phi_2}{\phi+\rho}.
\end{align}
The following two Theorems deal with the growth of non-negative real sequences. We will use these results in our main analysis of GSKM and PASKM method.

\begin{mdframed}[backgroundcolor=gray!15,   topline=false,   bottomline =false,   rightline=false,   leftline=false] \begin{theorem}
\label{th:0}
Let $\{G_k\}$ be a non-negative real sequence satisfying the following relation:
\allowdisplaybreaks{\begin{align*}
    G_{k+1} \leq \phi_1 G_k + \phi_2 G_{k-1}, \ \forall k\geq 1 \quad G_0 = G_1 \geq 0,
\end{align*}
if  $\phi_1, \phi_2  \geq 0$ and $\phi_1+\phi_2 < 1$ then the following bounds hold:
\begin{enumerate}
    \item  (Lemma 9 in \cite{loizou:2017}) Let, $\phi$ be the largest root of $\phi^2 + \phi_1 \phi - \phi_2 = 0$, then
    \begin{align*}
    G_{k+1} \leq (1+\phi) (\phi+ \phi_1)^k \ G_0, \ \forall k\geq 1.
\end{align*}
    \item Define $\rho = \phi + \phi_1$, then we have the following:
\allowdisplaybreaks{\begin{align*}
   \begin{bmatrix}
G_{k+1}  \\[6pt]
G_k 
\end{bmatrix} & \leq  \begin{dcases}
    \begin{bmatrix}
  R_1 \rho^{k+1}+ R_2 \phi^{k+1} \\[6pt]
R_1 \rho^{k}- R_2 \phi^{k} 
\end{bmatrix} \ G_0     \qquad k \ \text{even}; \\
  \begin{bmatrix}
R_3 \rho^{k}- R_4 \phi^{k} \\[6pt]
R_3 \rho^{k-1}+ R_4 \phi^{k-1} 
\end{bmatrix} \ G_0     \qquad   k \ \text{odd},
  \end{dcases}
\end{align*}}
where, $ 0 \leq \phi < 1$ and $ 0 < \rho = \phi + \phi_1 < 1$.
\end{enumerate}}

\end{theorem} \end{mdframed}

\begin{proof}
See Appendix 1.
\end{proof}


\allowdisplaybreaks{\begin{mdframed}[backgroundcolor=gray!15,   topline=false,   bottomline =false,   rightline=false,   leftline=false] \begin{theorem}
\label{th:seq2}
Let the real sequences $H_k \geq 0$ and $F_k \geq 0$ satisfy the following recurrence relation:
\begin{align}
\label{t-1}
\begin{bmatrix}
H_{k+1} \\
F_{k+1} 
\end{bmatrix} & \leq   \begin{bmatrix}
\Pi_1 & \Pi_2 \\
\Pi_3  & \ \Pi_4
\end{bmatrix} \begin{bmatrix}
H_{k} \\
F_{k}
\end{bmatrix}, 
 \end{align}
where, $\Pi_1, \Pi_2, \Pi_3, \Pi_4 \geq 0$ such that the following relations
\begin{align}
    \label{t0}
   \Pi_1\Pi_4 - \Pi_2\Pi_3 \geq 0,  \qquad \Pi_1+ \Pi_4 < 1+ \min\{1, \Pi_1\Pi_4 - \Pi_2\Pi_3 \},
\end{align}
hold. Then the sequence $\{H_k\}$ and $ \{F_k\}$ converges and the following result holds:
\allowdisplaybreaks{\begin{align*}
  \begin{bmatrix}
H_{k+1}  \\[6pt]
F_{k+1} 
\end{bmatrix} & \leq   \begin{bmatrix}
\Pi_1 & \Pi_2 \\
\Pi_3  & \ \Pi_4
\end{bmatrix}^k \begin{bmatrix}
H_{1} \\
F_{1}
\end{bmatrix} =  \begin{bmatrix}
\Gamma_2 \Gamma_3 (\Gamma_1-1) \ \rho_1^{k}+ \Gamma_1 \Gamma_3 (\Gamma_2+1)\ \rho_2^{k} \\[6pt]
\Gamma_3 (\Gamma_1-1) \ \rho_1^{k}+ \Gamma_3 (\Gamma_2+1)\ \rho_2^{k}
\end{bmatrix} \ \begin{bmatrix}
H_{1} \\
F_1 
\end{bmatrix}.
\end{align*}}
where,
\allowdisplaybreaks{\begin{align}
    \label{t1}
    & \Gamma_1 = \frac{\Pi_1-\Pi_4+\sqrt{(\Pi_1-\Pi_4)^2+4\Pi_2\Pi_3}}{2\Pi_3},  \nonumber \\
    & \Gamma_2 = \frac{\Pi_1-\Pi_4-\sqrt{(\Pi_1-\Pi_4)^2+4\Pi_2\Pi_3}}{2\Pi_3}, \ \Gamma_3 = \frac{\Pi_3}{\sqrt{(\Pi_1-\Pi_4)^2+4\Pi_2\Pi_3}}, \nonumber \\
    & \rho_1 = \frac{1}{2} \left[\Pi_1+\Pi_4 - \sqrt{(\Pi_1-\Pi_4)^2+4\Pi_2\Pi_3}\right], \nonumber \\
    & \rho_2 = \frac{1}{2} \left[\Pi_1+\Pi_4 + \sqrt{(\Pi_1-\Pi_4)^2+4\Pi_2\Pi_3}\right],
\end{align}}
and $  \Gamma_1, \Gamma_3 \geq 0$ and $ 0 \leq \rho_1 \leq \rho_2 < 1$.

\end{theorem} \end{mdframed}}

\begin{proof}
See Appendix 1.
\end{proof}


\subsection{Convergence Analysis of the GSKM Method}
\label{subsec:gskm conv}
In this subsection, we study convergence properties of the proposed GSKM method, i.e., we study the convergence behavior of the quantities of  $\E[\|x_k-\mathcal{P}(x_k)\|]$ and $\E[f(x_k)]$. For any $\xi \in \R$, let us define the sets $Q, Q_1, Q_2$ as
\begin{align}
\label{eq:s}
  & Q_1 = \{\xi \ | \ 0 \leq \xi \leq 1\}, \ \ Q = Q_1 \cup Q_2, \nonumber \\
  & Q_2 = \{ -1 < \xi \leq 0 \ | \  (1+\xi) \ \sqrt{h(\delta)} -\xi \left(1+ \delta \sqrt{\mu_2}\right) < 1 \}.
\end{align}
We proved that whenever $\xi \in Q$ and $0< \delta < 2$, the proposed GSKM method enjoys a global linear rate. We also provided convergence analysis of the function values (i.e., $f(x_k)$) with respect to the Cesaro average. Our results are global in nature and to the best of our knowledge, this is the first of its kind result for the SKM method.

\allowdisplaybreaks{\begin{mdframed}[backgroundcolor=gray!15,   topline=false,   bottomline =false,   rightline=false,   leftline=false] \begin{theorem}
\label{th:1}
Let $\{x_k\}$ be the sequence of random iterates generated by algorithm \ref{alg:gskm}. With the choice of parameters, $0 < \delta < 2$ and $0 \leq \xi  \leq 1 \ (\xi \in Q_1)$, the sequence of iterates $\{x_k\}$ converges and the following results hold:
\allowdisplaybreaks{\begin{enumerate}
\item Take $ \phi_1 = (1-\xi) h(\delta), \ \phi_2 = \xi h(\delta)$ and $\rho, \phi $ as in equation \eqref{def:0}, then
\begin{align*}
\E [d(x_{k+1},P)^2]  \leq  \rho^{k}  (1+\phi)  d(x_{0},P)^2 \ \text{and} \ \E [f(x_k)]  \leq  \frac{\mu_2(1+\phi)}{2}  \rho^{k}  d(x_{0},P)^2.
\end{align*}
\item Take $ \phi_1 = (1-\xi) h(\delta)$ and $ \phi_2 = \xi h(\delta)$, then
\allowdisplaybreaks{\begin{align*}
  \E \begin{bmatrix}
d(x_{k+1},P)^2  \\[6pt]
d(x_{k},P)^2 
\end{bmatrix} & \leq  \begin{dcases}
    \begin{bmatrix}
  R_1 \rho^{k+1}+ R_2 \phi^{k+1} \\[6pt]
R_1 \rho^{k}- R_2 \phi^{k} 
\end{bmatrix} \ d(x_{0},P)^2    \qquad  k \ \text{even}; \\
  \begin{bmatrix}
R_3 \rho^{k}- R_4 \phi^{k} \\[6pt]
R_3 \rho^{k-1}+ R_4 \phi^{k-1} 
\end{bmatrix} \ d(x_{0},P)^2     \qquad   k \ \text{odd},
  \end{dcases}
\end{align*}}
where, the constants $R_1, R_2, R_3, R_4$ are defined in equation \eqref{def:0} and $ 0 \leq \phi , \phi_1, \phi_2 < 1$ and $ 0 < \rho = \phi + \phi_1 < 1$.
\item Also the average iterate $\Tilde{x}_k = \sum \limits_{l=1}^{k} x_l$ for all $k \geq 0$ satisfies the following
\begin{align*}
    \E[d(\Tilde{x}_k,P)^2]  \leq \frac{(1+\phi) \ d(x_0,P)^2}{ k(1-\rho)} \quad \text{and} \quad \E[f(\Tilde{x}_k)] \leq \frac{(1+\xi) d(x_0,P)^2}{2 \delta k (2-\delta)}.
\end{align*}

\end{enumerate}}
\end{theorem} \end{mdframed}}

\begin{proof}
See Appendix 2.

\end{proof}

In the above Theorem, we obtain a global linear rate for the GSKM method with $0\leq \xi \leq 1$. Note that, when $0\leq \xi \leq 1$, we have,
\begin{align*}
   \rho =  \phi + \phi_1 & =  \frac{(1-\xi)h(\delta)+\sqrt{(1-\xi)^2 h^2(\delta)+4\xi h(\delta)}}{2} \\
    & \geq \frac{(1-\xi) h(\delta)+\sqrt{(1-\xi)^2 h^2(\delta)}}{2} = (1-\xi) h(\delta).
\end{align*}
Since the maximum value of $(1-\xi) h(\delta)$ can be derived as $h(\delta)$, the above inequality attains equality when $\xi = 0$ (see the next Corollary).  This gives us $1 > \rho = \phi + \phi_1 \geq h (\delta)$. Since the rate of the SKM algorithm is given by $h(\delta)$, we can say that the theoretical convergence rate of Algorithm \ref{alg:gskm} is always worse or equal compared to SKM whenever $0 \leq \xi \leq 1$.

\begin{mdframed}[backgroundcolor=gray!15,   topline=false,   bottomline =false,   rightline=false,   leftline=false] \begin{corollary} 
\label{cor:1}
(Theorem 1.3 in \cite{haddock:2017}) Let $\{x_k\}$ be the sequence of random iterates generated by the SKM method (algorithm \ref{alg:skm}) starting with $x_0 \in \R^n$. With $0 < \delta < 2$, the sequence of iterates $\{x_k\}$ converges and the following result holds:
\begin{align*}
\E \left[d(x_{k+1}, P)^2\right]  \leq \ \left[h(\delta)\right]^{k} \ d(x_{0}, P)^2.
\end{align*}
\end{corollary} \end{mdframed}
\begin{proof}
Note that, if we let $\xi = 0$ in the GSKM method, then we have $x_{k+1} = z_k$, which is precisely the SKM method. Now, take $\xi = 0$ in Theorem \ref{th:1}, then considering the first part of the Theorem, we have $\rho = h(\delta)$. Furthermore, from the second part, we have $R_3 \rho^k -R_4 \phi^k = R_1 \rho^{k+1} + R_2 \phi^{k+1} = \rho^k = \left(h(\delta)\right)^k$. This proves the result of Corollary \ref{cor:1} which is precisely the convergence rate obtained in \cite{haddock:2017} for the SKM method.
\end{proof}

Our next Theorem, states that, for a range of negative values of the parameter $\xi$, the GSKM method enjoys a global linear rate.

\allowdisplaybreaks{\begin{mdframed}[backgroundcolor=gray!15,   topline=false,   bottomline =false,   rightline=false,   leftline=false] \begin{theorem}
\label{th:2}
Let $\{x_k\}$ be the sequence of random iterates generated by algorithm \ref{alg:gskm} and let $0 < \delta < 2$ and $\xi \in Q_2$. Define
\begin{align}
\label{6}
    \Pi_1 = \sqrt{h(\delta)}, \ \Pi_2 = |\xi|, \ \Pi_3 =\delta \sqrt{\mu_2 h(\delta)}, \ \Pi_4 = |\xi| \left(1+ \delta \sqrt{\mu_2}\right),
\end{align}
and $\Gamma_1, \Gamma_2, \Gamma_3, \rho_1, \rho_2$ as in \eqref{t0} with the parameter choice of \eqref{6}. Then the sequence of iterates $\{x_k\}$ converges and the following result holds:
\allowdisplaybreaks{\begin{align*}
  \E \begin{bmatrix}
d(x_{k+1}, P)  \\[6pt]
\|z_{k+1}-z_k\| 
\end{bmatrix} & \leq   \begin{bmatrix}
-\Gamma_2 \Gamma_3 \ \rho_1^{k}+ \Gamma_1 \Gamma_3 \ \rho_2^{k} \\[6pt]
- \Gamma_3 \ \rho_1^{k}+ \Gamma_3 \ \rho_2^{k} 
\end{bmatrix} \ d(x_{0}, P),
\end{align*}}
where $  \Gamma_1, \Gamma_3 \geq 0$ and $ 0 \leq \rho_1 \leq \rho_2 < 1$.

\end{theorem} \end{mdframed}}

\begin{proof}
See Appendix 2.
\end{proof}

\paragraph{\textbf{Parameter Choice for GSKM}} Now we discuss allowable parameter selection for the GSKM algorithm based on Theorem \ref{th:1} and \ref{th:2}.

\begin{figure}[h!]
 \centering
    \includegraphics[scale = 0.8]{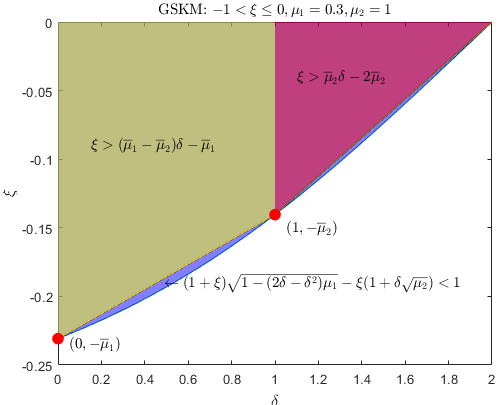}
    \caption{Allowable parameter range}
    \label{fig:param1}
\end{figure}
From Theorem \ref{th:1}, it can be noted that the GASKM method will converge for any $0 \leq \xi \leq 1$. 
Whenever $\xi$ is negative (i.e., $\xi \in Q_2$), the allowable range for $\xi$ can be shown in the following figure. In Figure \ref{fig:param1}, we plot the feasible region for allowable $\xi$ values for $\mu_1 = 0.3, \mu_2 =1$. Denote, $\Tilde{\mu}_1 = \frac{\mu_1}{\mu_1+\sqrt{\mu_2}}$ and $\Tilde{\mu}_2 = \frac{1-\sqrt{1-\mu_1}}{1-\sqrt{1-\mu_1}+\sqrt{\mu_2}}$. Then the feasible region of Figure \ref{fig:param1} can be approximated piece-wise as $\xi > (\Tilde{\mu_1} -\Tilde{\mu_2} )\delta - \Tilde{\mu_1}$ for $ 0 < \delta \leq 1 $ and $\xi > \Tilde{\mu_2} \delta - 2\Tilde{\mu_2}$ for $1 \leq \delta < 2$. Moreover, any $(\xi, \delta)$ pair that resides inside the region $\{0 < \delta < 2, \ -1< \xi <0, \  \xi \ \geq \ 0.5 \Tilde{\mu_1} (\delta-2) \}$ also resides inside the feasible region of Theorem \ref{th:2}.

\paragraph{Cesaro Average:} In the next Theorem, we propose the convergence analysis of the function values $f(x)$, with respect to the Cesaro average. Instead of bounding $\E[f(x_k)]$ in terms of initial function value $f(x_0)$, we bound the decay in terms of a larger quantity that results in a better convergence rate. To the best of our knowledge, this is the first result that shows $\mathcal{O}(\frac{1}{k})$ convergence of the Kaczmarz type methods for solving linear feasibility problems \footnote{Several works exits for the Kaczmarz type methods for solving linear systems \cite{richtrik2017stochastic, loizou:2017}.}. An interesting corollary of our method is the Cesaro average result for the SKM method. Furthermore, the result holds under weaker assumptions than the previous Theorems.

\begin{mdframed}[backgroundcolor=gray!15,   topline=false,   bottomline =false,   rightline=false,   leftline=false] \begin{theorem}
\label{th:cesaro}
Let $\{x_k\}$ be the random sequence generated by Algorithm \ref{alg:gskm}. Take, $-1 < \xi \leq 0 $ and $0 < \delta < \frac{2(1+\xi)}{1-2\xi}$. Define $\Tilde{x_k} = \frac{1}{k} \sum \limits_{l =1}^{k}x_l$ and $f(x)$ as in \eqref{def:function}, then
\begin{align*}
    \E \left[f(\bar{x}_k)\right] \leq \frac{ (1+\xi) (1+\xi-2\delta \xi \mu_2) \ d(x_0,P)^2+ 2 \xi \delta (\delta \xi -\delta  -1) f(x_0)}{2 \delta k \left(2+2 \xi + 2 \delta \xi -\delta\right)}.
\end{align*}
\end{theorem} \end{mdframed}

\begin{proof}
See Appendix.
\end{proof}

\begin{mdframed}[backgroundcolor=gray!15,   topline=false,   bottomline =false,   rightline=false,   leftline=false] \begin{corollary}
\label{cor:cesaro}
Let $\{x_k\}$ be the random sequence generated by SKM method (algorithm \ref{alg:skm}). Define $\Tilde{x_k} = \frac{1}{k} \sum \limits_{l =1}^{k}x_l$ and $f(x)$ as in \eqref{def:function}, then 
\begin{align*}
    \E \left[f(\bar{x}_k)\right] \leq \frac{d(x_0,P)^2}{2 \delta k \left(2- \delta\right)},
\end{align*}
holds for any  $ 0 < \delta  < 2$.
\end{corollary} \end{mdframed}

\begin{proof}
Take $\xi = 0$ in Theorem \ref{th:cesaro}, then the result follows.
\end{proof}

The next Theorem is an extension of the result obtained in \cite{haddock:2017} and to a certain extent, it can be taken as an extension of Telgen's result \cite{telgen:1982}. The Theorem gives one a certificate of feasibility after a finite number of GSKM iterations. Before delving into the Theorem, we will provide some known Lemmas for the SKM algorithm which holds for the GSKM algorithm too. We refer interested readers to the work of De-Loera \textit{et. al} \cite{haddock:2017} for detailed proof of these Lemmas (Lemma \ref{lem:skm1} to Lemma \ref{lem:skm4}).

\begin{mdframed}[backgroundcolor=gray!15,   topline=false,   bottomline =false,   rightline=false,   leftline=false] \begin{lemma}
\label{lem:skm1}
(Lemma 1 in \cite{haddock:2017}) Define, $\theta(x) = \left[\max_{i}\{a_i^Tx-b_i\}\right]^{+} $ as the maximum violation of point $x \in \R^n$ and the length of the binary encoding of a linear feasibility problem with rational data-points as
\begin{align*}
  \sigma = \sum \limits_{i}  \sum \limits_{j} \ln{\left(|a_{ij}|+1\right)} + \sum \limits_{i}  \ln{\left(|b_{i}|+1\right)} + \ln{(mn)} +2.
\end{align*}
Then if the rational system $Ax \leq b$ is infeasible, for any $x\in \R^n$, the maximum violation $\theta(x)$ satisfies the following lower bound:
\begin{align*}
    \theta(x) \ \geq \ \frac{2}{2^{\sigma}}.
\end{align*}
\end{lemma} \end{mdframed}

\begin{mdframed}[backgroundcolor=gray!15,   topline=false,   bottomline =false,   rightline=false,   leftline=false] \begin{lemma}
\label{lem:skm3}
(Lemma 4 in \cite{haddock:2017}) If $P$ is $n$-dimensional (full-dimensional) then the sequence of iterates $\{x_k\}$ generated by the GSKM method converges to a point $x \in P$.
\end{lemma} \end{mdframed}
\begin{proof}
Since, by assumption, $P$ is full dimensional, then the rest of the proof follows the same argument as Lemma 4 in \cite{haddock:2017}.
\end{proof}
\begin{mdframed}[backgroundcolor=gray!15,   topline=false,   bottomline =false,   rightline=false,   leftline=false] \begin{lemma}
\label{lem:skm4}
(\cite{KHACHIYAN:1980}) If the rational system $ Ax \leq b$ is feasible, then there is a feasible solution $x^{*}$ whose coordinates satisfy $|x^{*}_j| \leq \frac{2^{\sigma}}{2n}$ for $j = 1, ..., n$.
\end{lemma} \end{mdframed}

\paragraph{Certificate of feasibility:} To detect feasibility of the rational system $Ax \leq b$, one needs to find a point $x_k$ such that $\theta(x_k) < 2^{1-\sigma}$. Such a point if exists will be called a certificate of feasibility. When the system is feasible, one expects to find a certificate of feasibility after finitely many iterations, and that if one fails to find a certificate after finitely many iterations, one can obtain a lower bound on the probability that the system is infeasible. Moreover, as discussed in the next Theorem, if the system is feasible, one can bound the probability of finding a certificate of feasibility.

\begin{mdframed}[backgroundcolor=gray!15,   topline=false,   bottomline =false,   rightline=false,   leftline=false] \begin{theorem}
\label{th:3}
Suppose $A, b$ are rational matrices with binary encoding length, $\sigma$, and that we run the GSKM method ($0 < \delta < 2$, \  $ \xi \in Q$) on the system $ Ax \leq b\ (\|a_i\| = 1, i = 1,2,...,m)$ with $x_0 = 0$. Suppose the number of iterations $k$ satisfies the following lower bound:
\begin{align*}
   \frac{4 \sigma - 4 -\log n + \log (1+\phi)}{\log \left(\frac{1}{\bar{\rho}}\right)} < k.
\end{align*}
If the system $Ax \leq b$ is feasible, then,  
\begin{align*}
    p \ \leq H(\sigma, \phi, k, \bar{\rho}) = \sqrt{\frac{1+\phi}{n}} \ 2^{2\sigma -2} \ \bar{\rho}^{\frac{k}{2}},
\end{align*}
where $p$ is the probability that the current iterate is not a certificate of feasibility. And $\bar{\rho} = \max\{\rho, \rho_2^2\} < 1$, where $\rho$ and $\rho_2$ are defined in Theorem \ref{th:1} and Theorem \ref{th:2} for the choice $\xi \in Q_1$ and $\xi \in Q_2$, respectively. Also note that the function $H(\sigma, \phi, k, \bar{\rho})$ is a decreasing function with respect to $k$.
\end{theorem} \end{mdframed}

\begin{proof}
See Appendix.

\end{proof}

\begin{remark}
\label{rem:2}
Note that instead of a normalized system if we consider a non-normalized system $\overline{A}x \leq \overline{b}, \ \|\overline{a_i}\| \neq 1$ for some $i$, then suppose the number of iterations $k$ satisfies the following lower bound:
\begin{align*}
   \frac{4 \overline{\sigma} - 4 -\log n + \log (1+\phi) + 2 \log \psi}{\log \left(\frac{1}{\bar{\rho}}\right)} < k,
\end{align*}
where $\overline{\sigma}$ is the binary encoding length for $\overline{A}, \overline{b}$. If the system $\overline{A}x \leq \overline{b}$ is feasible, then,  
\begin{align*}
    p \ \leq \  \sqrt{\frac{1+\phi}{n}} \ 2^{2\overline{\sigma} -2} \ \psi \ \bar{\rho}^{\frac{k}{2}},
\end{align*}
where $p = $ probability that the current update $x_k$ is not a certificate of feasibility and $\psi = \max_{j} \|\overline{a_j}\|$.
\end{remark}

\begin{mdframed}[backgroundcolor=gray!15,   topline=false,   bottomline =false,   rightline=false,   leftline=false] \begin{corollary}
\label{cor:2}
(Theorem 1.5 in \cite{haddock:2017}) Suppose $\overline{A}, \overline{b}$ are rational matrices with binary encoding length, $\overline{\sigma}$, and that we run the SKM method on the system $ \overline{A}x \leq \overline{b}\ (\|\overline{a_i}\| \neq 1$ for some $i)$ and $x_0 = 0$. Suppose the number of iterations $k$ satisfies the following lower bound:
\begin{align*}
   \frac{4 \overline{\sigma} - 4 -\log n + 2 \log \psi}{\log \left(\frac{1}{h(\delta)}\right)} < k,
\end{align*}
where $\overline{\sigma}$ is the binary encoding length for $\overline{A}, \overline{b}$. If the system $\overline{A}x \leq \overline{b}$ is feasible, then,  
\begin{align*}
    p \ \leq \  \sqrt{\frac{1}{n}} \ 2^{2\overline{\sigma} -2} \ \psi \ \left[h(\delta)\right]^{\frac{k}{2}},
\end{align*}
where $p = $the probability that the current update $x_k$ is not a certificate of feasibility and $\psi = \max_{j} \|\overline{a_j}\|$.
\end{corollary} \end{mdframed}
\begin{proof}
Take $\xi = 0$ in Theorem 3. Then, we have, $\phi = 0, \ \rho = \phi + \phi_1 = h(\delta) = \rho_2^2$. It can be easily checked that the GSKM method with $\xi = 0$ is just the SKM method. Now, considering Theorem \ref{th:3} with the above parameter choice, we can get the bound of Corollary \ref{cor:2}.
\end{proof}


\subsection{Convergence Analysis of the PASKM Method}
\label{subsec:paskm conv}

In this subsection, we study convergence properties of the proposed PASKM algorithm, i.e., we study the convergence behavior of the quantities of  $\E[\|v_k-\mathcal{P}(v_k)\|^2]$, $ \ \E[\|x_k-\mathcal{P}(x_k)\|^2]$, $ \ \E[\|y_k-\mathcal{P}(y_k)\|^2]$ and $\E[f(x_k)]$ generated by the PASKM method. We proved that for a range of step parameters $\alpha, \gamma, \omega$, the proposed PASKM method enjoys a global linear rate. We also provided convergence analysis of the function values $f(x_k)$, with respect to the Cesaro average. The next Theorem deals with the convergence of the sequences $\{v_k\}$ and $\{y_k\}$ as well as the function values $f(x_k)$ generated by the PASKM algorithm.

\allowdisplaybreaks{\begin{mdframed}[backgroundcolor=gray!15,   topline=false,   bottomline =false,   rightline=false,   leftline=false] \begin{theorem}
\label{th:4}
Let $\{x_k\}$ be the sequence of random iterates generated by algorithm \ref{alg:paskm} and let $0 < \delta < 2$ and $ 0 \leq \alpha, \omega  \leq 1$ such that $\gamma +3 \omega -2 \leq 0$, $\omega  h(\delta) (1-\alpha) (1+\gamma) < 1$ and the following condition
\begin{align}
\label{cond}
    \omega (1+\gamma) + h(\delta) (1-\alpha) + \alpha (1-\omega) & + \alpha \gamma \mu_1 (\gamma + 3 \omega -2) \nonumber \\
    & -  \omega  h(\delta) (1-\alpha) (1+\gamma) < 1,
\end{align}
holds. Define, $\Pi_1 = \omega (1+\gamma), \ \Pi_2 = (1-\omega)+\gamma \mu_1 (\gamma + 3 \omega -2), \ \Pi_3 = \alpha \omega (1+\gamma)$, $ \ \Pi_4 = (1-\alpha) h(\delta) + \alpha (1-\omega)+ \alpha \gamma \mu_1 (\gamma + 3 \omega -2)$ and $\Gamma_1, \Gamma_2, \Gamma_3, \rho_1, \rho_2$ as in \eqref{t0}. Then the sequence of iterates $\{v_k\}$ and $\{y_k\}$ converges and the following results hold:
\allowdisplaybreaks{\begin{align*}
  \E \begin{bmatrix}
d(v_{k+1}, P)^2  \\[6pt]
d(y_{k+1}, P)^2
\end{bmatrix} & \leq  \begin{bmatrix}
\Gamma_2 \Gamma_3 (\Gamma_1-1) \ \rho_1^{k+1}+ \Gamma_1 \Gamma_3 (\Gamma_2+1)\ \rho_2^{k+1} \\[6pt]
\Gamma_3 (\Gamma_1-1) \ \rho_1^{k+1}+ \Gamma_3 (\Gamma_2+1)\ \rho_2^{k+1} 
\end{bmatrix} \ d(y_0,P)^2,
\end{align*}}
and
\begin{align*}
    \E \left(f(y_{k+1})\right) \leq \frac{\mu_2 }{2}  \left[\Gamma_3 (\Gamma_1-1) \ \rho_1^{k+1}+ \Gamma_3 (\Gamma_2+1) \ \rho_2^{k+1}  \right] \  d(y_0,P)^2.
\end{align*}
where $  \Gamma_1, \Gamma_3 \geq 0$ and $ 0 \leq \rho_1 \leq \rho_2 < 1$.

\end{theorem} \end{mdframed}}

\begin{proof}
See Appendix 3.
\end{proof}


The next Theorem deals with the convergence of the sequences $\{v_k\}$ and $\{x_k\}$ generated by the PASKM algorithm.

\begin{mdframed}[backgroundcolor=gray!15,   topline=false,   bottomline =false,   rightline=false,   leftline=false] \begin{theorem}
\label{th:10}
Let, $v_{k+1}$ and $x_{k+1}$ are generated by Algorithm \ref{alg:paskm}. If we select the parameters $\omega, \ \gamma, \ \alpha$ as
\begin{align*}
\omega = 1-\frac{\zeta \mu_1^2+2\gamma \mu_1-\zeta \mu_1}{1+\zeta \mu_1^2} , \ \ \gamma = \sqrt{ \zeta \eta \mu_1}, \ \ \alpha = \frac{ \eta}{ \eta + \gamma },
\end{align*}
where, $\zeta$ is chosen as $ 0 < \zeta  < \frac{4\eta \mu_1}{(1-\mu_1)^2} $ if $\mu_1 < 1$, otherwise choose any $\zeta > 0$. Then, for any $0 < \delta < 2$, the sequence of iterates $\{v_k\}, \ \{x_k\}$ converges and the following result holds:
\begin{align*}
\E \left[ d(v_{k+1},P)^2 + \zeta \mu_1 \  d(x_{k+1},P)^2 \right] \ & \leq \ \omega^{k+1} \ \E \left[ d(v_{0},P)^2+ \zeta \mu_1 \  d(x_{0},P)^2 \right] \\
& = (1+\zeta \mu_1) \ \omega^{k+1} \ d(x_{0},P)^2.
\end{align*}  
This theorem implies that the PASKM algorithm converges linearly with a rate of $\omega$, which accumulates to a total of $\mathcal{O}(\frac{1+\zeta \mu_1^2}{\zeta \mu_1^2+2\gamma \mu_1-\zeta \mu_1}\log{1/\epsilon})$ iterations to bring the given error below $\epsilon > 0$.
\end{theorem} \end{mdframed}

\begin{proof}
See Appendix 3.
\end{proof}

In the next Theorem, we present the convergence analysis of the function $f(x)$ with respect to the Cesaro average for the PASKM algorithm. We showed that the Cesaro average of the PASKM iterates converges to the optimum at a rate of $\mathcal{O}(1/k)$ where $k$ is the number of iterations.

\begin{mdframed}[backgroundcolor=gray!15,   topline=false,   bottomline =false,   rightline=false,   leftline=false] \begin{theorem}
\label{th:cesaro2}
Let $\{y_k\}$ be the random sequence generated by Algorithm \ref{alg:paskm}. Take, $\ 0 \leq 1-\alpha, \omega < 1 $, $ \ 0 < \delta < \frac{2(1-\omega + \alpha \omega)}{1+2\omega -2\alpha \omega}$ and $\alpha \gamma = \alpha \delta + \omega \delta (1-\alpha) $. Define $\Tilde{y_k} = \frac{1}{k} \sum \limits_{l =1}^{k}y_l$ and $f(y)$ as in \eqref{def:function}, then
\begin{align*}
    \E \left[f(\bar{y}_k)\right] \leq \frac{ (1-\omega + \alpha \omega )^2 \ d(y_0,P)^2+ 2\delta (\delta -2+ 3 \omega - 3 \alpha \omega + \delta \omega - \delta \alpha \omega) f(y_0)}{2 \delta k \left(2-2 \omega + 2 \alpha \omega - 2\delta \omega + 2 \delta \alpha \omega -\delta\right)}.
\end{align*}
\end{theorem} \end{mdframed}

\begin{proof}
See Appendix 3.
\end{proof}


\paragraph{\textbf{Parameter selection for PASKM algorithm}}

In this section, we discuss allowable parameter selection for the PASKM algorithm based on Theorem \ref{th:4}. If the parameters $0 \leq \alpha, \omega \leq 1$ and $\gamma \geq 0$ satisfies $\gamma + 3 \omega -2$ and the condition of \eqref{cond} hold then the PASKM method will converge for any $0 < \delta < 2$ \footnote{When, $\delta = 2$, we have $h(\delta) = 1-\eta \mu_1 = 1$. In that case, we can simplify the condition of \eqref{cond} as $\omega < \frac{2-\gamma}{3+\frac{1}{\mu_1}}$. In other words, for $\delta =2$ the PASKM algorithm will converge if we select the parameters as $0 \leq \gamma < 2$, $0 \leq \alpha \leq 1$, $\omega < \frac{2-\gamma}{3+\frac{1}{\mu_1}}$ and $ \mu_1 = \frac{\lambda_{\min}^+(A^TA)}{m}$.}. To simplify the conditions for ease of implementation, let's take $0 \leq \gamma < 2 $ and $ \omega = \frac{2-\gamma}{3+p}$ for some $0 \leq p \leq \frac{1}{\mu_1} $ \footnote{Note that for the choice $p > \frac{1}{\mu_1}$ the condition \eqref{eq:alpha} trivially holds as the right hand side of \eqref{eq:alpha} is always greater than $1$.}.
\begin{figure}[h!]
 \centering
    \includegraphics[scale = 0.9]{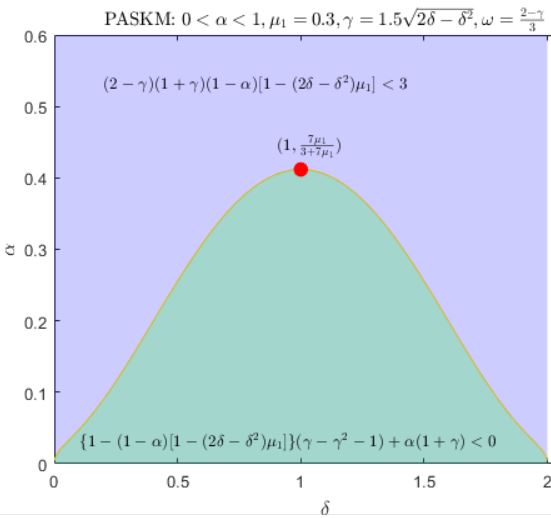}
    \caption{Allowable parameter range}
    \label{fig:param2}
\end{figure}

In Figure \ref{fig:param2}, we plot the feasible region considering the above parameter choice and the conditions of Theorem \ref{th:4}. Considering the choice of $\gamma$ and $\omega$, the condition $\omega  h(\delta) (1-\alpha) (1+\gamma) < 1$ simplifies to
\begin{align*}
    \alpha > 1-\frac{3+p}{(2-\gamma)(1+\gamma)h(\delta)} & = \frac{2h(\delta)+\gamma h(\delta) -\gamma^2 h(\delta)-3-p}{(2-\gamma)(1+\gamma)h(\delta)} \\
     = & \frac{(\gamma-\gamma^2-1-p)+\eta \mu_1(\gamma^2-2-\gamma)}{(2-\gamma)(1+\gamma)h(\delta)} \leq 0,
\end{align*}
where, we used the fact that the conditions $\gamma-\gamma^2-1 \leq 0$ and $\gamma^2-2-\gamma \leq 0$ hold for any $0 \leq \gamma \leq 2$. That implies for any $\alpha \geq 0$, the condition $\omega  h(\delta) (1-\alpha) (1+\gamma) < 1$ holds. Similarly, we can simplify the condition of \eqref{cond} as follows:
\begin{align}
\label{eq:alpha}
    \alpha & < \underbrace{\frac{(1+p-\gamma+\gamma^2)(1-h(\delta))}{1-h(\delta)+p+\gamma + (\gamma-p)h(\delta)- \gamma^2 h(\delta)+\mu_1p\gamma(\gamma-2)}}_{> \ 0 \ \text{for} \ 0 < \delta < 2} \nonumber \\
    & = \alpha(\gamma, \delta, p) \leq 1.
\end{align}
Therefore, if we choose $\gamma, \omega$ and $\alpha$ as
\begin{align}
   & \gamma = 1.5 \sqrt{2\delta - \delta^2}, \ p =0 , \  \omega = \frac{2-\gamma}{3+p}, \ \alpha = 0.99* \alpha(\gamma, \delta, p)  \label{eq:param1} \\
   & \gamma = 2 \sqrt{2\delta - \delta^2}, \ p =0, \  \omega = \frac{2-\gamma}{3+p}, \ \alpha = 0.99* \alpha(\gamma, \delta, p), \label{eq:param2} 
\end{align}
then the convergence result of Theorem \ref{th:4} holds for the PASKM algorithm.  We will use these two sets of parameter choices in our numerical experiments. Note that, our choice is empirical in nature. One can probably find a better combination of parameters than \eqref{eq:param1} and \eqref{eq:param2}. Similarly, if we choose $\gamma, \omega$ and $\alpha$ as
\begin{align}
\label{eq:askm12}
\zeta =\frac{3.99\eta \mu_1}{(1-\mu_1)^2} , \ \omega = 1-\frac{\zeta \mu_1^2+2\gamma \mu_1-\zeta \mu_1}{1+\zeta \mu_1^2} , \ \ \gamma = \sqrt{ \zeta \eta \mu_1}, \ \ \alpha = \frac{ \eta}{ \eta + \gamma },
\end{align}
then the convergence result of Theorem \ref{th:10} holds for the PASKM algorithm. The choice of \eqref{eq:askm12} is not of practical benefit as the value of $\frac{\lambda_{\min}^+(A^TA)}{m}$ is very small for most test cases. From \eqref{eq:askm12}, we have $\gamma \propto \frac{1}{m}$, which is very small for large test instances. Smaller $\gamma$ slows down the convergence of the PASKM algorithm as $\gamma$ can be seen as a projection parameter like $\delta$.

\section{Numerical Experiments}
\label{sec:num}

In this section, we discuss the numerical experiments performed to show the computational efficiency of the proposed algorithms (Algorithm \ref{alg:gskm} and \ref{alg:paskm}). As mentioned before, we limit our focus on the over-determined systems regime (i.e., $m \gg n$) where iterative methods are competitive in general. However, from our experiments, we see similar computational behavior for the under-determined systems as well. 

\subsection{Experiment Specifications}
We implemented the proposed GSKM and PASKM algorithms in \textit{MATLAB R2018b}  and performed the experiments in a Dell Precision 7510 workstation with 32GB RAM, Intel Core i7-6820HQ CPU, processor running at 2.70 GHz. To analyze computational performance, we perform the numerical experiments for a wide range of instances including both randomly generated and real-world test problems.
\begin{itemize}
    \item \textbf{Randomly generated problems:} Gaussian and highly correlated systems
    \item \textbf{Real-world test instances:} Standard ML data sets and Sparse Netlib LP instances
\end{itemize}
We compare SKM with two versions of the proposed GSKM and PASKM algorithms for a better understanding of the algorithmic behavior. In Table \ref{tab:2}, we provide the parameter choices for GSKM and PASKM algorithms. Throughout the numerical experiments section, we compared SKM with GSKM-1, GSKM-2 and PASKM-1, PASKM-2. 
\begin{table}[ht!]
\centering
\caption{Parameter choice of GSKM and PASKM algorithms for the numerical experiments.}
\begin{tabular}{|c|c|c|c|c|c|}
\hline
\multirow{2}{*}{Parameters} & \multicolumn{3}{c|}{\begin{tabular}[c]{@{}c@{}}GSKM (Algorithm \ref{alg:gskm}, $\xi \in Q$) \end{tabular}} & \multicolumn{2}{c|}{\begin{tabular}[c]{@{}c@{}}PASKM (Algorithm \ref{alg:paskm}, $\alpha, \omega, \gamma$) \end{tabular}} \\ \cline{2-6} 
 & \begin{tabular}[c]{@{}c@{}}  SKM\end{tabular} & \begin{tabular}[c]{@{}c@{}}  GSKM-1\end{tabular} & \begin{tabular}[c]{@{}c@{}}  GSKM-2\end{tabular} & \begin{tabular}[c]{@{}c@{}}  PASKM-1\end{tabular} & \begin{tabular}[c]{@{}c@{}}  PASKM-2\end{tabular} \\ \hline
\begin{tabular}[c]{@{}c@{}} $1 \leq \beta \leq m$\\ $0 < \delta < 2$\end{tabular} & $\xi = 0$ & \begin{tabular}[c]{@{}c@{}}$\xi = -0.1$\\ $\xi = -0.2$\end{tabular} & $\xi = 0.5$ & \begin{tabular}[c]{@{}c@{}} $\alpha, \omega, \gamma$ as in \eqref{eq:param1} \end{tabular} & \begin{tabular}[c]{@{}c@{}}$\alpha, \omega, \gamma$ as in \eqref{eq:param2}\end{tabular} \\ \hline
\end{tabular}
\label{tab:2}
\end{table}

Finally, we investigate the performance behavior of the proposed GSKM and PASKM methods with state-of-the-art methods such as Interior point methods (IPMs) and Active set methods (ASMs) for several Netlib LP instances. The total CPU time is calculated in seconds (s). For a fair comparison, we run the algorithms 10 times and report the averaged performance throughout the experiments. Moreover, all the algorithms start from the same initial point that is far away from the feasible region.

\subsection{Experiments on Randomly Generated Instances}
\label{subsec:numRndm}
We considered the linear feasibility $Ax \leq b$, where the entries of matrices $A \in \R^{m \times n}$ and $b\in \R^{m}$ are chosen randomly from a certain distribution. To maintain the system consistency (i.e., $b \in \mathcal{R}(\mathbf{A})$), we first generated vectors $x_1, x_2 \in \R^n$ at random from the corresponding distributions, then set $b$ as the convex combination of vectors $Ax_1$ and $Ax_2$ (i.e., $b = \sigma Ax_1 + (1-\sigma) Ax_2,\ 0 \leq \sigma \leq 1$). Two types of random data sets are considered: highly correlated, and Gaussian. For the correlated systems, data matrices $A$ and $x_1, x_2$ are chosen uniformly at random between $[0.9,1.0]$ (i.e., $a_{ij}, x_j \in [0.9,1.0], \ i = 1,2,...,m, \ j =1,2,...,n$). For the Gaussian system data matrices, $A$ and $x_1, x_2$ are chosen uniformly at random from standard normal distribution (i.e., $a_{ij}, x_j \in \mathcal{N}(0,1), \forall i,j$). Moreover, the vector $b \in \R^m$ is generated by following the above-mentioned procedure.

\begin{figure}[h!]
    \includegraphics[scale = 0.525]{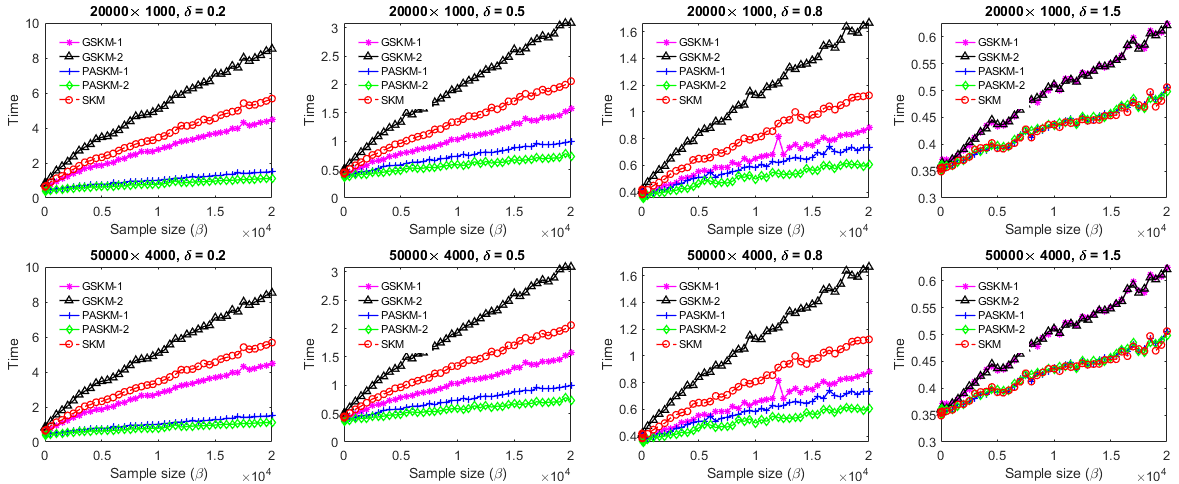}
    \caption{Sample size $\beta$ VS average CPU time comparison among SKM, GSKM, PASKM variants for $\delta = 0.2, 0.5, 0.8, 1.5$ on correlated systems. Problem size: $20000 \times 1000$ (Top panel), $50000 \times 4000$ (Bottom panel).}
    \label{fig:1}
\end{figure}

\paragraph{\textbf{CPU time vs Sample size $\beta$}} We first compared the total CPU time of the proposed algorithms (GSKM-1, GSKM-2, PASKM-2, PASKM-2) with the original SKM algorithm. The comparison is carried out by varying the sample size $\beta$ from $1$ to the total row size $m$. The positive residual error tolerance is chosen as $10^{-05}$ (i.e., $\|\left(Ax-b\right)^+\|_2 \leq 10^{-05}$). The comparison is carried out for $\delta = 0.2, 0.5, 0.8$ and $1.5$. In Figure \ref{fig:1}, we compared the above-mentioned algorithms for two randomly generated highly correlated linear feasibility problems of size $20000 \times 1000$ and $50000 \times 4000$. From Figure \ref{fig:1}, we see that the proposed GSKM-1, PASKM-1, PASKM-2 algorithms outperform the SKM algorithm in terms of average CPU time when $\delta = 0.2, 0.5, 0.8$. For $\delta = 1.5$, the performance of SKM, PASKM-1 and PASKM-2 are fairly similar whereas, the performance of GSKM-1 and GSKM-2 are worse compared to all other algorithms.

\begin{figure}[h!]
    \includegraphics[scale = 0.5]{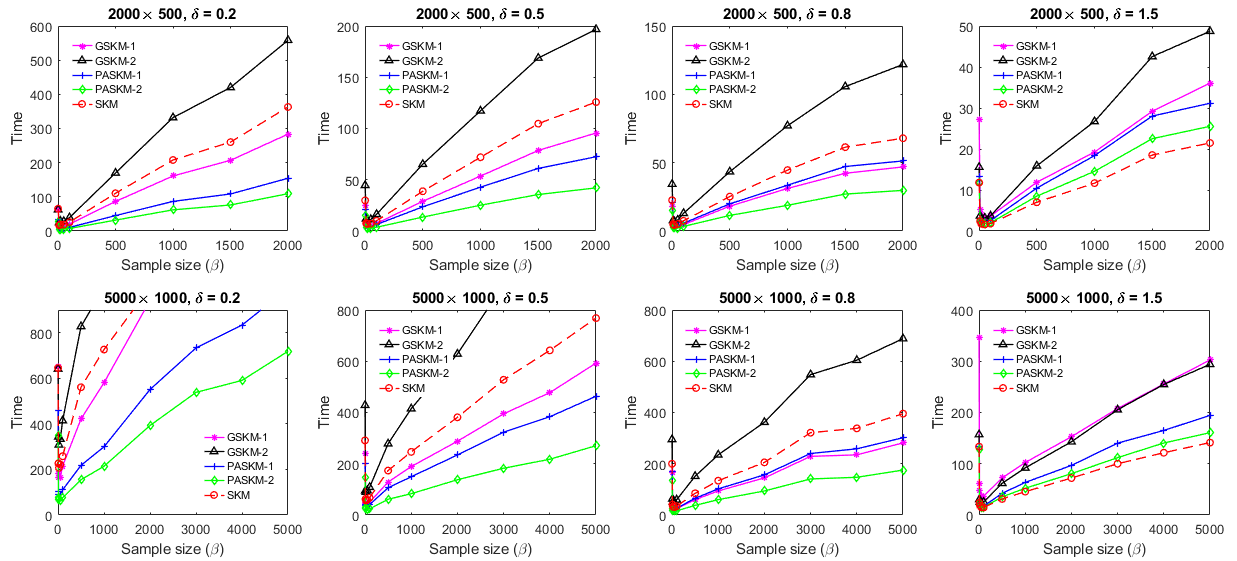}
    \caption{Sample size $\beta$ VS average CPU time comparison among SKM, GSKM, PASKM variants for $\delta = 0.2, 0.5, 0.8, 1.5$ on Gaussian systems. Problem size: $2000 \times 500$ (Top panel), $5000 \times 1000$ (Bottom panel).}
    \label{fig:2}
\end{figure}
\begin{figure}[h!]
    \includegraphics[scale = 0.54]{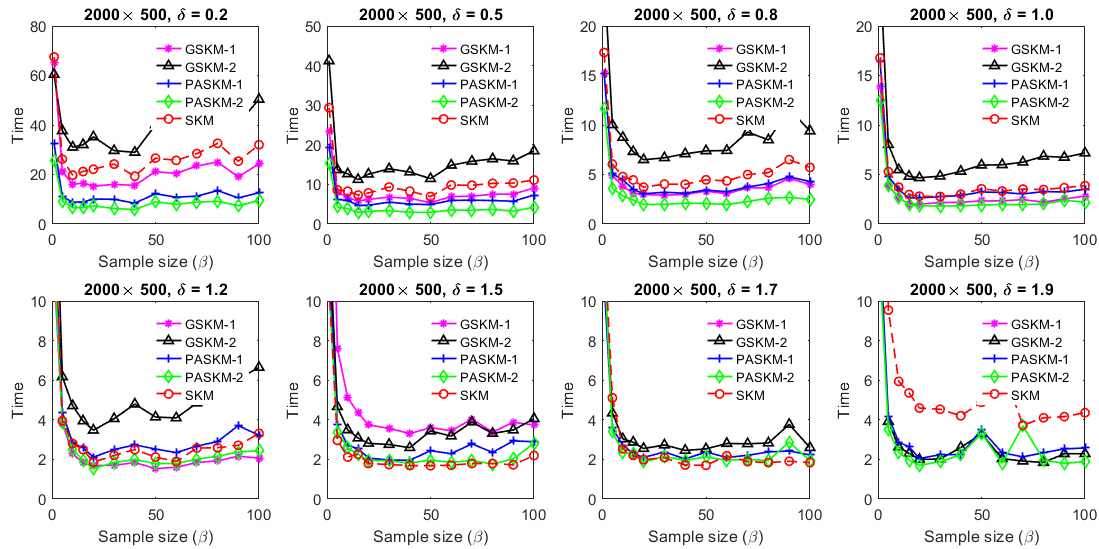}
    \caption{Sample size $\beta$ VS average CPU time comparison among SKM, GSKM, PASKM variants for $\delta = 0.2, 0.5, 0.8,1, 1.2, 1.5, 1.7, 1.9$ and samller sample size (i.e., $1 \leq \beta \leq 100$) on a $2000 \times 500$ Gaussian system.}
    \label{fig:small}
\end{figure}

We present the time versus sample size plot for two randomly generated Gaussian system of size $2000 \times 500$ and $5000 \times 1000$ in Figure \ref{fig:2}. All the algorithms show similar performance patterns as shown in the correlated systems (Figure \ref{fig:1}) for the choice of $0 < \delta < 1$. However, for the case of $\delta = 1.5$, SKM and PASKM-2 perform marginally better than the other algorithms. Since all of the considered methods perform significantly well whenever $\beta$ is small (i.e., $1 < \beta \leq 100$). For a better understanding, we compare the proposed algorithms for $1 < \beta \leq 100$. In Figure \ref{fig:small}, we plot the time vs $\beta$ graph for a $2000 \times 500$ Gaussian problem for smaller $\beta$.

In a nutshell, we can conclude that for the choice of $0 < \delta < 1$, PASKM-1, PASKM-2 and GSKM-1 outperform the original SKM method. And in that region, PASKM-2 is the best performing algorithm. Moreover, for $1.5 \leq \delta \leq 1.7$, all of the proposed algorithms perform similarly as the SKM method. However, for the case of $\delta = 1.9$, the proposed algorithms significantly outperform the SKM method. Furthermore, we believe with correct parameter choice one can find better-performing variants of GSKM and PASKM compared to the SKM algorithm for the case of $1.5 \leq \delta \leq 1.7$. Finally from Figure \ref{fig:small}, we can deduce that the best sample size choice for all of the considered methods occurs at $1 < \beta \ll m$. This amplifies the importance of the special sampling distribution selection.

\begin{figure}[h!]
    \includegraphics[scale = 0.7]{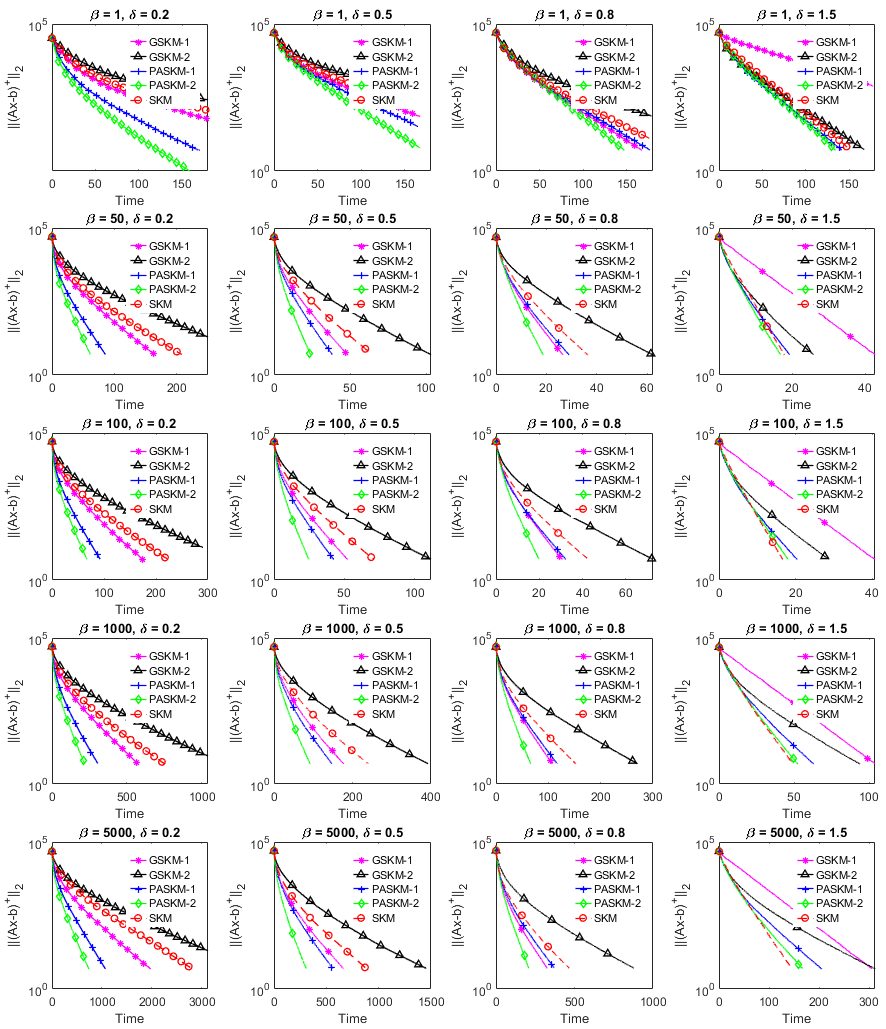}
    \caption{Positive residual error $\|\left(Ax-b\right)^+\|_2$ VS CPU time comparison among SKM, GSKM, PASKM variants for $\delta = 0.2, 0.5, 0.8, 1.5$ and $\beta = 1,50,100,1000,5000$ on $5000 \times 1000$ Gaussian system.}
    \label{fig:4}
\end{figure}

\paragraph{\textbf{Positive residual error $\|\left(Ax-b\right)^+\|_2$ VS No. of iterations and Time}} Now, we compare the respective convergence trend for the considered algorithms with respect to the number of iterations and CPU time. 
\begin{figure}[h!]
\centering
    \includegraphics[scale = 0.69]{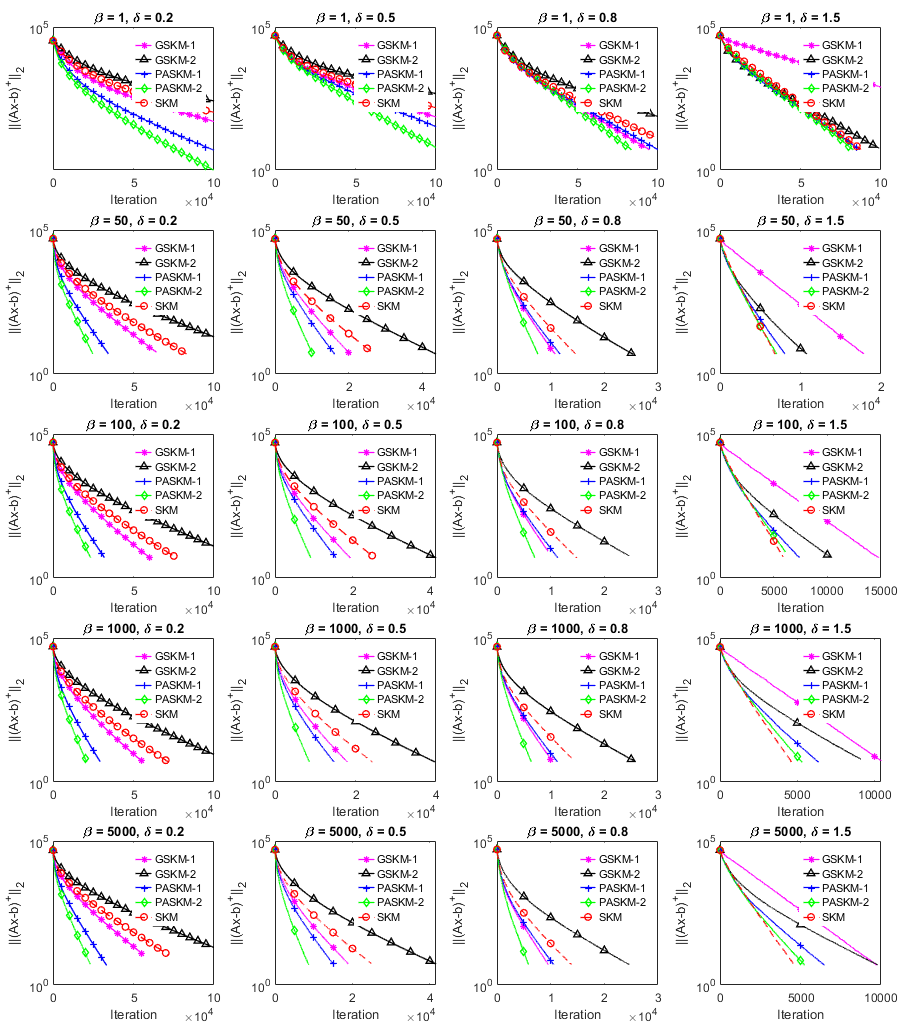}
    \caption{Positive residual error $\|\left(Ax-b\right)^+\|_2$ VS No. of iteration comparison among SKM, GSKM, PASKM variants for $\delta = 0.2, 0.5, 0.8, 1.5$ and $\beta = 1,50,100,1000,5000$ on $5000 \times 1000$ Gaussian system.}
    \label{fig:3}
\end{figure}
We choose positive residual error $\|(Ax-b)^+\|_2$ as the convergence measure and considered $5000 \times 1000$ Gaussian system. We carried out the analysis for several choices of sample sizes, $\beta = 1, 100, 1000, m$ and the choice of $\delta$ values remains the same as before. In Figures \ref{fig:3} and \ref{fig:4}, we provide the respective positive residual decay results for different sample sizes and different projection parameters. We plot positive residual error VS iteration and positive residual error VS time in Figures \ref{fig:3} and \ref{fig:4}, respectively. From Figures \ref{fig:3} and \ref{fig:4}, we see that irrespective of sample size, $\|(Ax_k-b)^+\|_2$ converges to zero much faster for the proposed PASKM-1 and PASKM-2, GSKM-1 compared to SKM whenever $\delta < 1$. For the case of $\delta = 1.5$, SKM and PASKM-2 has a similar kind of performance whereas the GSKM-1 performs poorly compared to SKM and PASKM method. As expected, the choice $\beta =1$ produces the slowest rate and the choice $\beta =100$ produces the best convergence graph.

\begin{figure}[h!]
    \includegraphics[scale = 0.7]{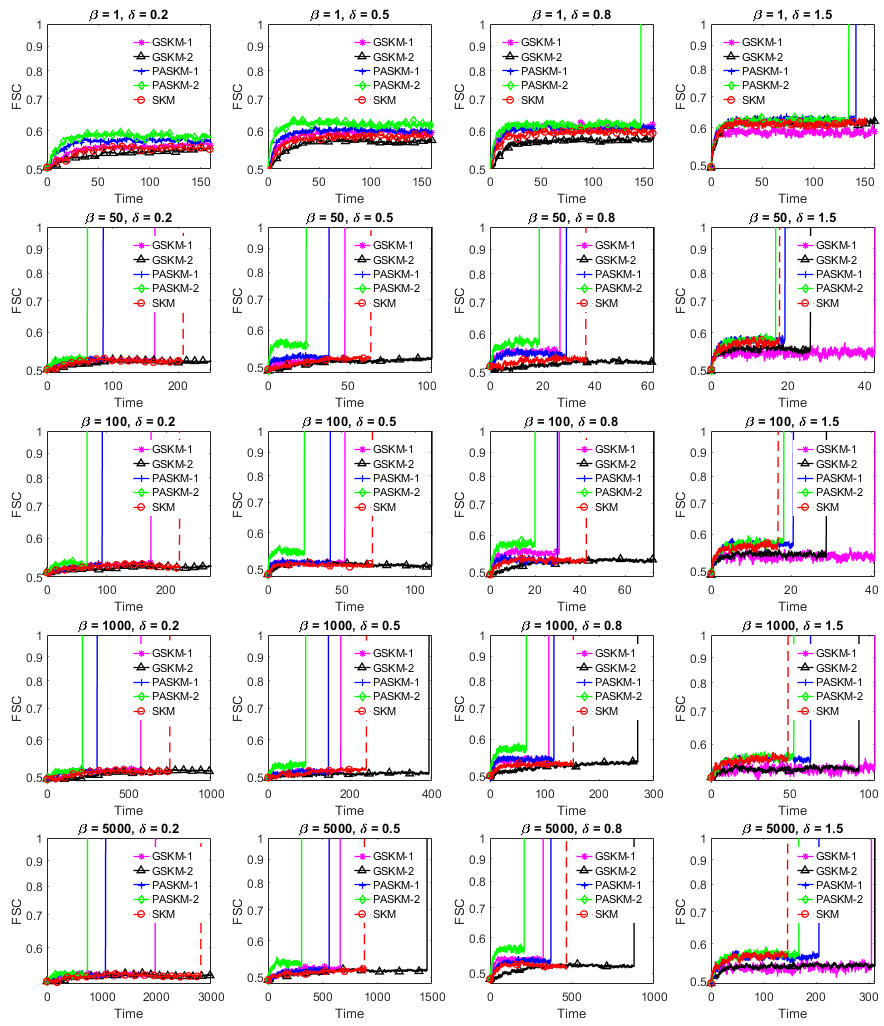}
    \caption{CPU time vs fraction of satisfied constraints (FSC) comparison among SKM, GSKM, PASKM variants for $\delta = 0.2, 0.5, 0.8, 1.5$ and $\beta = 1,50,100,1000,5000$ on $5000 \times 1000$ Gaussian system.}
    \label{fig:6}
\end{figure}

\paragraph{\textbf{Fraction of satisfied constraints (FSC) VS No. of iterations and Time}} To investigate the generated solution quality of the above-mentioned algorithms of Table \ref{tab:2}, we measure the number of satisfied constraints at each iteration, for that we define,
\begin{align*}
    \text{Fraction of Satisfied Constraints (FSC)} \ = \frac{\text{Number of satisfied constraints}}{\text{Total number of constraints ($m$)}}
\end{align*}
Note that, at any particular iteration we have, $0 \leq \text{FSC} \leq 1$.
\begin{figure}[h!]
    \includegraphics[scale = 0.69]{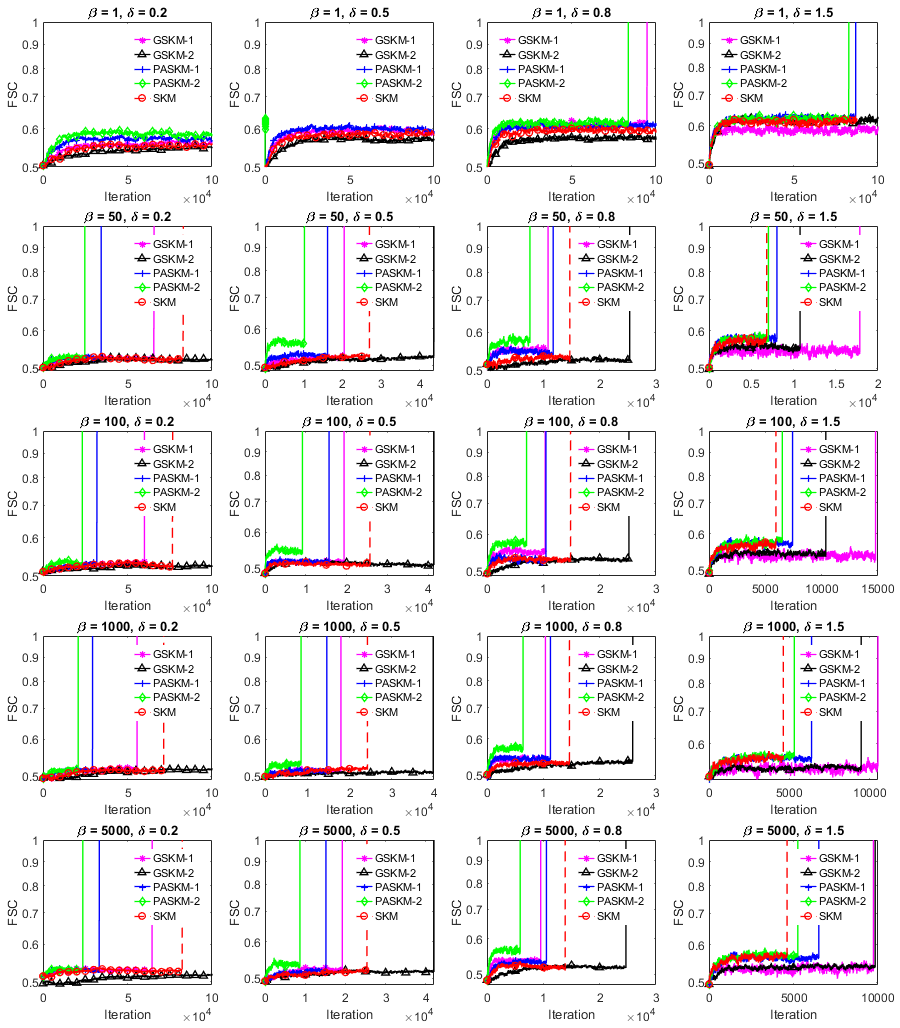}
    \caption{No. of iteration vs fraction of satisfied constraints (FSC) comparison among SKM, GSKM, PASKM variants for $\delta = 0.2, 0.5, 0.8, 1.5$ and $\beta = 1,50,100,1000,5000$ on $5000 \times 1000$ Gaussian system.}
    \label{fig:5}
\end{figure}
In Figures \ref{fig:5} and \ref{fig:6}, we plot the value of FSC with respect to No. of iterations and CPU time of each algorithm respectively. From Figures \ref{fig:5} and \ref{fig:6}, we can see that the choice of $\beta = 1$ is the worst choice for all algorithms as the improvement of FSC is much slower compared to other choices of $\beta$. And for the choice $\beta = 100$, we get the best solution quality for each algorithm. Our proposed GSKM-1, PASKM-1 and PASKM-2 algorithms outperform the other methods significantly for $0 < \delta <1$ but, for $\delta = 1.5$ only PASKM-2 performs similar to SKM. 

\subsection{Experiments on real-world Instances}
\label{subsec:reallife}

In this subsection, we consider some nonrandom, real-world test instances. For the sake of unbiased performance analysis, we consider the following two types of real-world data-sets: standard Machine Learning (ML) data-sets for Support Vector Machine (SVM) classifier \cite{yeh:2009,lichman:2013,haddock:2017}, and sparse linear feasibility problems extracted from benchmark Netlib LP problems \cite{netlib}.

\paragraph{\textbf{SVM classifier instances}} We first consider two linear feasibility problems arising from binary classification with SVM. We compare the proposed algorithms with SKM to the linear classification problem using the SVM model for the following two data sets: 1) Wisconsin (diagnostic) breast cancer data set and 2) Credit card default data set. The Wisconsin breast cancer data set consists of data points whose features are calculated from images. There are two types of data points: 1) malignant and 2) benign cancer cells. As shown by the researchers \cite{calafiore:2014,haddock:2017}, the SVM classifier problem can be re-written as an equivalent homogeneous system of linear inequalities ($Ax \leq 0$), which represents the separating hyper-plane between malignant and benign data points. The constraint matrix $A$ has $569$ rows (data points) and $30$ columns (features). Since the data set is not perfectly separable, we allow tolerance for the positive residual $\|(Ax)^+\|$. For our experiments, we fixed the tolerance as $10^{-3}$ (i.e., we ran the algorithm until $\|(Ax_k)^+\| \leq 10^{-3}$ is satisfied).

Similarly, we consider the credit card default data set described in \cite{yeh:2009,haddock:2017}. This data set consists of features denoting the payment profile of a user and binary variables describing payment conditions in a certain billing cycle: 1 for payment made on time and 0 for late payment. The SVM classification problem for the data set can be transformed into an equivalent homogeneous system of inequalities ($Ax \leq 0$) like before. The solution $x^{*}$ denotes the coefficients of the separating hyper-plane between on-time and default data points. The transformed data matrix $A$ has $30000$ rows ($30000$ user profiles) and $23$ columns ($22$ profile features). As the data set is not separable, like the previous problem we allow a tolerance error. In this case, we ran the algorithms until the condition: $\|(Ax_k)^+\|/\|(Ax_0)^+\| \leq 10^{-3}$ is satisfied.

\begin{figure}[h!]
    \includegraphics[scale = 0.455]{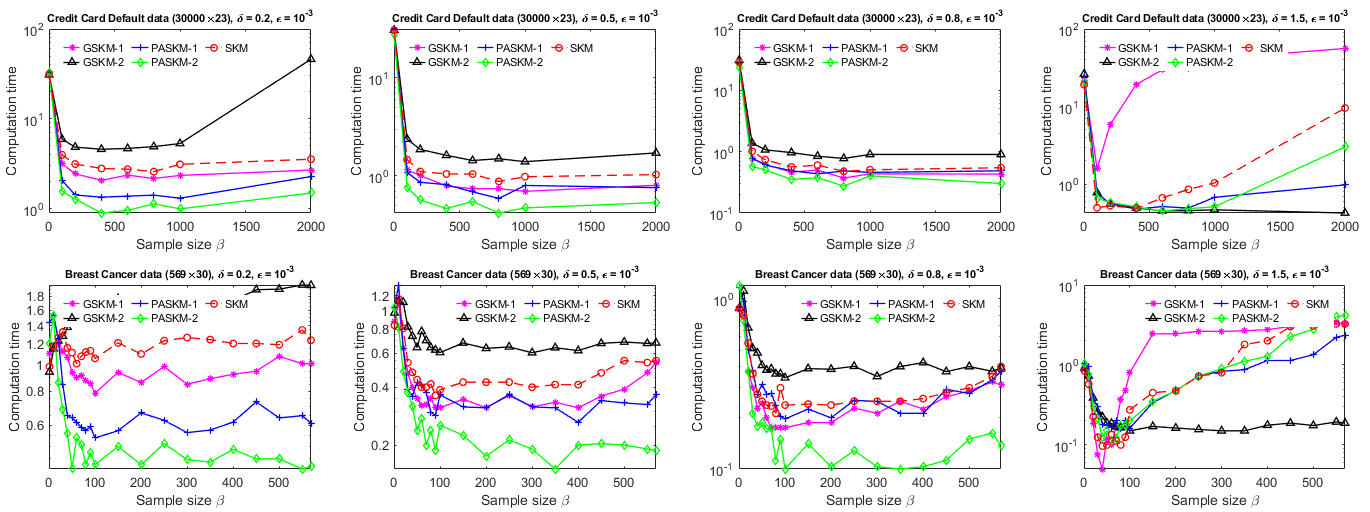}
    \caption{Average CPU time VS Sample size $\beta$ comparison among SKM, GSKM, PASKM variants for $\delta = 0.2, 0.5, 0.8, 1.5$ on Support Vector Machine problems; Top panel: Credit card data set, Bottom panel: Wisconsin breast cancer data set.}
    \label{fig:7}
\end{figure}

\paragraph{\textbf{CPU time vs Sample size $\beta$}}
We plot the CPU time VS sample size $\beta$ graphs for SVM problems in Figure \ref{fig:7}. To be consistent with our previous experiments, we choose $\delta = 0.2, 0.5, 0.8, 1.5$. From Figure \ref{fig:7}, we see that the proposed GSKM-1, PASKM-1 and PASKM-2 algorithms outperform the other algorithms including SKM for $\delta = 0.2, 0.5, 0.8$. However, for $\delta = 1.5$, GSKM-2 performs significantly well compared to the other methods. On the other hand, SKM, PASKM-1 and PASKM-2 follow a similar trend across different sample sizes. PASKM-1 and PASKM-2 marginally outperform SKM for this regime. Another interesting point can be noted that the comparison graphs for the credit card data set are not as smooth as the breast cancer data set graphs \footnote{The credit card data matrix has $30,000$ rows. From our earlier experiments, we observe that the choice of $1 < \beta \leq 100$, the proposed algorithms produce the best performance. For that reason, we plot the credit card graph up-to $\beta = 2000$. The irregularity of the credit card graph occurs when $\beta > 2000$.}, which can be attributed to the irregularity of the constraint matrix $A$. 

\paragraph{\textbf{Netlib LP instances}}
We also investigate the comparative performance of the proposed algorithms with SKM on real-world sparse data sets. For this experiment, we consider some Netlib LP \cite{netlib} test instances. Each of these problems is formulated as a standard linear programming problem ( $\min c^Tx$ subject to $Ax = b, \ l \leq x \leq u$). To conduct the above-mentioned experiments, we transform each of these problems into an equivalent linear feasibility problem.  

\begin{figure}[h!]
    \includegraphics[scale = 0.69]{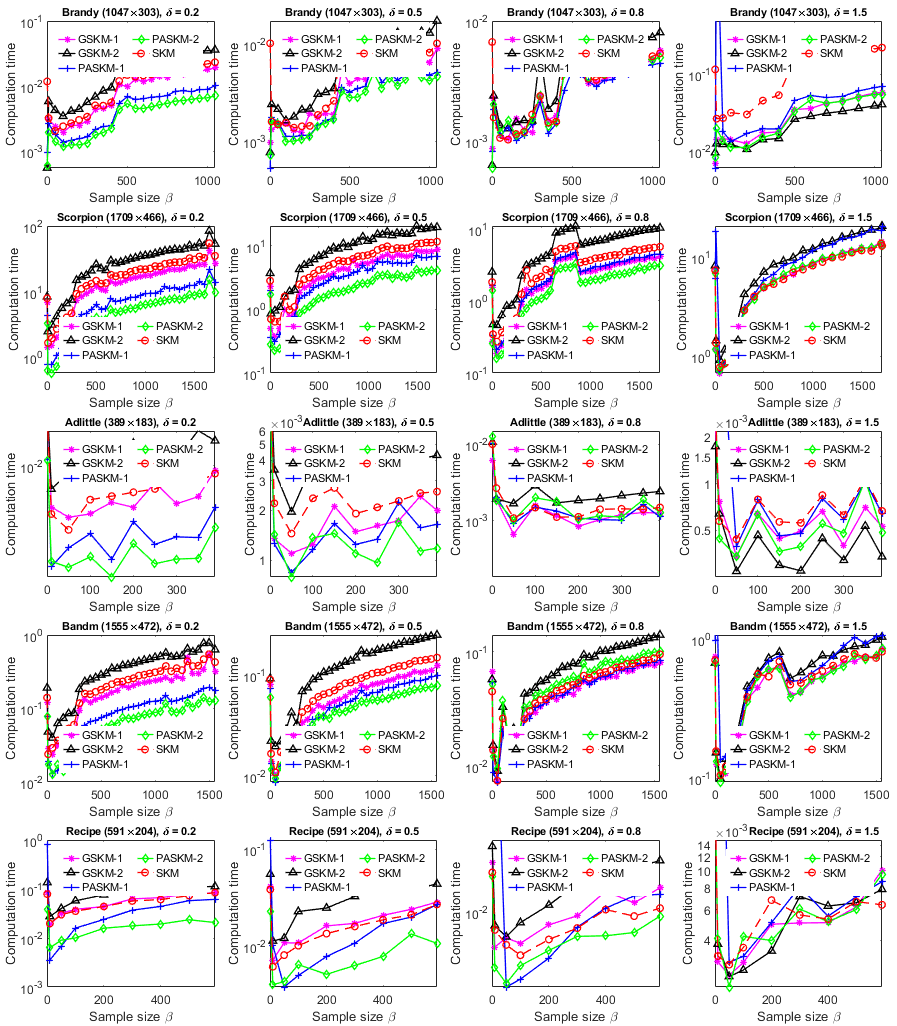}
    \caption{Average CPU time VS Sample size $\beta$ comparison among SKM, GSKM, PASKM variants for $\delta = 0.2, 0.5, 0.8, 1.5$ on Netlib LP instances.}
    \label{fig:8}
\end{figure}

\paragraph{\textbf{CPU time vs Sample size $\beta$}}
Now, we plot the CPU time VS sample size $\beta$ graphs for five Netlib LP instances in Figure \ref{fig:8}. Later in subsection \ref{subsec:3}, we consider a total of ten Netlib LP instances including the five considered here. In Figure \ref{fig:8}, we provide comparison graphs for the following Netlib LP test instances: lp-brandy, lp-addlittle, lp-scorpion, lp-bandm, lp-recipe. Furthermore, we consider different error tolerances for these problems (see Table \ref{tab:3} for details). From Figure \ref{fig:8}, we see that the proposed GSKM-1, PASKM-1, PASKM-2 algorithms outperform the SKM algorithm for $\delta = 0.2, 0.5, 0.8$. In the case of $\delta = 1.5$, the performance of SKM, GSKM-1 and PASKM-2 are fairly similar for the problems lp-scorpion, lp-bandm and lp-recipe. For lp-brandy and lp-adlittle, all of the proposed variants of GSKM and PASKM outperform the original SKM.

\subsection{Comparison with IPM and ASM for Netlib LP instances}
\label{subsec:3}
In this subsection, we compare the performance of GSKM and PASKM variants with SKM and benchmark commercial solvers for solving Netlib LP test instances. 
We follow the standard framework used by De Loera \textit{et. al} \cite{haddock:2017} and Morshed \textit{et. al} \cite{Morshed2019} in their work for linear feasibility problems. The problem instances are transformed from standard LP problems (i.e., $\min c^Tx$ subject to $Ax = b, \ l \leq x \leq u$ with optimum value $p^*$) to an equivalent linear feasibility formulation (i.e.,$\mathbf{A}x \leq \mathbf{b}$, where $\mathbf{A} = [A^T \ -A^T \ I \ -I \ c]^T$ and $\mathbf{b} = [b^T \ -b^T \ u^T \ -l^T \ p^*]^T$). For all of the experiments, we compared the proposed algorithms for $0 < \delta < 1$, since from our experiments in subsection \ref{subsec:numRndm} and \ref{subsec:reallife}, this is the domain where the proposed GSKM and PASKM variants significantly outperform the SKM method.

In Table \ref{tab:3}, we list the total CPU time in seconds for each of the above-mentioned algorithms in Table \ref{tab:2}. In addition to that, we provide the CPU time for Interior point method (IPM) and Active set method (ASM) algorithms for solving the selected Netlib LP problems. For a better and fair comparison, the pseudo-code of the proposed methods and SKM is written in MATLAB and Optimization Toolbox function \texttt{fmincon} is used to implement IPM and ASM methods. We first solve the linear feasibility problem ($\mathbf{A} x \leq \mathbf{b}$) with SKM, GSKM and PASKM variants and record the CPU time in Table \ref{tab:3}. Note that, we can't use \texttt{fmincon}'s IPM and ASM algorithms directly to solve the linear feasibility problem ($\min 0, \ s.t \ \mathbf{A}x \leq \mathbf{b} $) since both methods fail to solve the linear feasibility problems. The reason for that is, in IPM the \textit{Karush Kuhn Tucker} (KKT) system at each iteration becomes singular, and ASM stops in the first step of finding a feasible solution. 
\begin{table}[h!]
\centering
\caption{CPU time comparisons among the state-of-the-art methods (using MATLB's \texttt{fmincon} function) solving LP, and SKM, GSKM and PASK solving LF. $^*$ implies that the solver was unable to solve the problem with predetermined accuracy within 100,000 function evaluations. CPU time of the best performing algorithm for a problem is represented in bold letters.}
\label{tab:3}
\adjustbox{max width=\textwidth}{
\begin{tabular}{@{}|c|c|c|c|c|c|c|c|c|@{}}
\hline
Instance      & Dimensions & \begin{tabular}[c]{@{}c@{}} GSKM \\  $\times10^{-2}$  \end{tabular}  & \begin{tabular}[c]{@{}c@{}} PASKM \\  $\times10^{-2}$  \end{tabular}     & \begin{tabular}[c]{@{}c@{}} SKM \\  $\times10^{-2}$  \end{tabular}  & \begin{tabular}[c]{@{}c@{}} Interior \\   Point\end{tabular} & \begin{tabular}[c]{@{}c@{}}Active\\    set\end{tabular}                  & $\beta$ &  \begin{tabular}[c]{@{}c@{}} $\epsilon $ \\   $\times10^{-2}$ \end{tabular}\\ \hline
adlittle & $389 \times 138$    & 0.027 & \textbf{0.173}  & 0.032                                                        & 2.16                                                        & 4.96 & 150  & 0.1 \\ \hline
agg      & $2207 \times 615$   & 0.22  & \textbf{0.196}  & 0.23                                   & $66.54^*$         & $315.91^*$ & 50   & 1 \\ \hline
bandm    & $1555 \times 472$   & 9.82   & \textbf{4.057}   & 9.2                                & 14.57        & $529.43^* $ & 50  & 1 \\ \hline
blend    & $337 \times 114$    & 1.48   & \textbf{0.581}   & 1.28                                                        & 2.28                                                    & 4.62 & 50  & 0.1 \\ \hline
brandy   & $1047 \times 303$ & 0.53  & \textbf{0.491} & 14.06                                                       & 16.97                                                   & 63.11 & 1 & 1  \\ \hline
degen2   & $2403 \times 757$   & 26.26   & \textbf{10.139}  & 20.73 & 7.13  & 21038 & 100  & 1 \\ \hline
finnis   & $3123 \times 1064$  & 0.53   & 0.532   & \textbf{0.527}                                   & $66.16^*$          & $237750^*$                & 10  & 0.1 \\ \hline
recipe   & $591 \times 204$    & 0.60   & \textbf{0.164}   & 0.52                                                        & 0.89                                                        & 63.24                 & 50  & 0.1 \\ \hline
scorpion & $1709 \times 466$   & 156.9   & \textbf{42.712}  & 125                                                       & 17.68                                                       & 8.02                 & 50 & 1 \\ \hline
stocfor1 & $565 \times 165$    & 1.05   & \textbf{0.553}   & 0.95                                                       & 2.13                                                        & 2.52                 & 50  & 0.1 \\ \hline
\end{tabular}}
\end{table}
For a fair comparison, in Table \ref{tab:3}, we list the total CPU consumption time as follows: for the SKM method we solve the feasibility problem ($\mathbf{A} x \leq \mathbf{b}$) for a certain $\beta $ and $\delta$ \footnote{we note the best possible time from our previous experiments}, for GSKM and PASKM variants we solve the same feasibility problem and report the best performing method from each of the two, and finally for \texttt{fmincon} algorithms, we use the original LPs ($\min c^Tx \ s.t \ Ax \leq b, \ l \leq x\leq u $). Note that, this is not an ideal or obvious comparison, for a better suitable comparison we follow the framework used in \cite{haddock:2017,Morshed2019}. We set the stopping criterion for SKM, GSKM and PASKM variants as  $\frac{\max(\mathbf{A}x_k-\mathbf{b})}{\max(\mathbf{A}x_0-\mathbf{b})} \leq \epsilon$ and the halting  criterion for the \texttt{fmincon}'s algorithms (IPM, ASM) are set as $\frac{\max(Ax_k-b, l-x_k, x_k-u)}{\max(Ax_0-b, l-x_0, x_0-u)} \leq \epsilon$ and $\frac{c^Tx_k}{c^Tx_0} \leq \epsilon$, where $\epsilon$ is the tolerance gap listed in Table \ref{tab:3}. To avoid any biased conclusion, for each problem we set the initial update as far as possible from the feasible region.

From the comparison in Table \ref{tab:3}, we can see that the proposed algorithms work much faster than IPM and ASM but work marginally better than the existing SKM method. Notice that the improvement of PASKM and GSKM algorithms over the SKM method for most problems are marginal as the proposed algorithms are designed explicitly for dense matrices. One can develop special algorithmic variants of the proposed PASKM and GSKM methods for sparse problems by following some standard aggregation techniques. A possible technique is to combine multiple steps by using the sparsity of the test instances. For instance, after $k^{th}$ iteration when we have $x_k, y_k$ and $v_k$, instead of moving forward with the sequences $x_{k+1}, y_{k+1}$ and $v_{k+1}$, for any $T \gg 1$ we can skip $T$ iterations and update $x_{k+T}, y_{k+T}$ and $v_{k+T}$ using a generalized recurrence relation that can enhance the computational efficiency.

\section{Conclusion}
\label{sec:colc}

In this work, we propose a general algorithmic framework (GSKM) for solving linear feasibility problems that unify various SKM type algorithms with the addition of a relaxation parameter $\xi$. From our convergence analysis of the GSKM method, one can recover convergence Theorems of several well-known algorithms such as Randomized Kaczmarz, Motzkin Method and Sampling Kaczmarz Motzkin method. In addition to the general framework, we propose a Nesterov type acceleration scheme in the SKM method called as PASKM. Our proposed PASKM method provides a bridge between Nesterov type acceleration of Machine Learning to sampling Kaczmarz methods for solving linear feasibility problems. To show the effectiveness of the proposed algorithms, we performed a wide range of numerical experiments on various types of random and standard benchmark data sets. For a better understanding of the behavior of the proposed algorithms, we numerically analyze two variants for both GSKM and PASKM algorithms in comparison with the original SKM method. Furthermore, we compare our proposed methods to commercially available methods such as IPM and ASM. In the majority of the test instances, the proposed algorithms significantly outperform the state-of-the-art methods. Furthermore, as shown in our numerical experiments, the correct choice of parameters can lead to much faster and accelerated methods for different types of test instances.

\paragraph{\textbf{Future Research}} In the future, the proposed algorithms and the technical analysis can be adopted effectively to various types of extensions such as sparse variants, optimally tuned PASKM, and GSKM, PASKM variants with greedy sampling strategies. First, we plan to extend our work to design efficient sparse variations of the proposed methods that can handle large-scale real-world problems with greater sparsity in the data matrix $A$. Second, we intend to design a test instance dependent scheme for identifying optimal parameter selection (i.e., $\beta$, $\delta$, $\xi$, $\lambda$, $\tau$) for both GSKM and PASKM. For the GSKM algorithm, adaptive parameter selection (i.e., $\beta_k$, $\delta_k$, $\xi_k$) policy can be a great area of future research. One can also derive connecting ideas between the proposed GSKM and induced projection plane generation of Chubanov \cite{Chubanov:2012,Chubanov:2015} which can produce faster algorithms. Finally, we aspire to develop adaptive sampling strategies and integrate the greedy Kaczmarz \cite{greedbai:2018} type method into the GSKM framework to further speed up the convergence.





\section{Acknowledgements}

The authors are truly grateful to the anonymous referees and the editors for their valuable comments and suggestions in the earlier version of the paper. The comments helped immensely in the revision process and greatly improved the quality of this paper.


\section*{Appendix 1}

\paragraph{Proof of Lemma \ref{lem2}}
Using the expectation expression given in the first section we have,
\begin{align*}
     \E_{\mathbb{S}} \left[a_{i^*}a_{i^*}^T\right] & = \frac{1}{\binom{m}{\beta}} \sum\limits_{j = 0}^{m-\beta} \binom{\beta-1+j}{\beta-1} (A^TA)_{\underline{\mathbf{i_j}}} \\
     &  \preceq \ \frac{\binom{m-1}{\beta-1}}{\binom{m}{\beta}} \sum\limits_{j = 0}^{m-\beta}  (A^TA)_{\underline{\mathbf{i_j}}} \ \preceq \ \frac{\beta}{m} \sum\limits_{i = 1}^{m}  a_{i} a_{i}^T  \ = \ \frac{\beta}{m} A^TA.
\end{align*}
Here, the notation $(A^TA)_{\underline{\mathbf{i_j}}}$ denotes the matrix $a_l a_l^T$ where the index $l$ belongs to the list \eqref{eq:sampling}. Furthermore, the index $l$ corresponds to the $(\beta+ j)^{th}$ entry on the list \eqref{eq:sampling}. This proves the Lemma.



\paragraph{Proof of Lemma \ref{lem3}}
Using the definition of the expectation from \eqref{def:exp}, we have,
\allowdisplaybreaks{
\begin{align*}
 \E_{\mathbb{S}} \left[\big | (a_{i^*}^Tx-b_{i^*})^+\big |^2\right] & \overset{\eqref{def:exp}}{=} \frac{1}{\binom{m}{\beta}} \sum\limits_{j = 0}^{m-\beta} \binom{\beta-1+j}{\beta-1}  | (Ax-b)^{+}_{\underline{\mathbf{i_j}}} |^2 \\
    & \overset{\text{Lemma} \ \ref{lem:skmseq}}{\geq} \frac{1}{\binom{m}{\beta}} \sum\limits_{j = 0}^{m-\beta}  \frac{\sum\limits_{l = 0}^{m-\beta} \binom{\beta-1+l}{\beta-1}}{m-\beta +1} | (Ax-b)^{+}_{\underline{\mathbf{i_j}}} |^2 \\
    & \geq \frac{1}{m-\beta+1} \sum\limits_{j = 0}^{m-\beta}  \big |(Ax - b)^+_{\underline{\mathbf{i_j}}} \big |^2  \\
    & \geq \frac{1}{m-\beta+1} \min\{\frac{m-\beta+1}{m-s}, 1\} \ \|(Ax - b)^+\|^2 \\
    &  \overset{\text{Lemma} \ \ref{lem0}}{\geq} \ \frac{1}{m L^2}  \ d(x,P)^2. 
\end{align*}}
Here, $s$ is the number of zero entries in the residual $(Ax - b)^+$, which also corresponds to the number of satisfied constraints for $x$. Since $A \mathcal{P} (x) \leq b$, we have the following:
\begin{align*}
    \E_{\mathbb{S}} \left[\big | (a_{i^*}^Tx-b_{i^*})^+\big |^2\right]   & \overset{\eqref{def:exp}}{=}  \frac{1}{\binom{m}{\beta}} \sum\limits_{j = 0}^{m-\beta} \binom{\beta-1+j}{\beta-1} | (Ax-b)^{+}_{\underline{\mathbf{i_j}}}|^2 \\
    & \leq  \frac{1}{\binom{m}{\beta}} \sum\limits_{j = 0}^{m-\beta} \binom{\beta-1+j}{\beta-1} \ \big | (Ax-A\mathcal{P} (x))_{\underline{\mathbf{i_j}}}   \big|^2 \\
    & =  \frac{1}{\binom{m}{\beta}} (x-\mathcal{P} (x))^T\sum\limits_{j = 0}^{m-\beta} \binom{\beta-1+j}{\beta-1} \ (A^TA)_{\underline{\mathbf{i_j}}} (x-\mathcal{P} (x) ) \\
    & =  (x-\mathcal{P} (x) )^T \E_{\mathbb{S}} \left[a_{i^*}a_{i^*}^T\right] (x-\mathcal{P} (x) ) \\
    & \overset{\text{Lemma} \ \ref{lem1} \ \& \ \ref{lem2} }{\leq} \  \min \left\{1, \frac{\beta}{m} \lambda_{\max}\right\} \big \| x- \mathcal{P} (x) \big \|^2 \\
    & = \  \min \left\{1, \frac{\beta}{m} \lambda_{\max}\right\} d(x,P)^2.
\end{align*}
Combining the above identities and using the expression for $f(x)$ from \eqref{def:function} we get,
\begin{align*}
    \frac{\mu_1}{ 2} \ d(x,P)^2 \ \leq \ f(x) \ \leq \  \frac{\mu_2}{2} \ d(x,P)^2,
\end{align*}
which proves the Lemma. 

\paragraph{Proof of Lemma \ref{lem:grad}}

From the definition of $f(x)$, it can be easily checked that $f(x)$ is a convex function. Now, by the convexity property of $f(x)$, for any $x, y \in \R^n$, we have the following:
\begin{align}
\label{eq:grad}
    \langle x-y, \nabla f(y) \rangle \leq f(x) - f(y).
\end{align}
Therefore, we have
\begin{align*}
\big \langle  x-y, \E_{\mathbb{S}}  \left[(a_{i^*}^Ty-b_{i^*})^{+} a_{i^*}\right] \big \rangle & =  \langle x-y, \nabla f(y) \rangle  \leq  \ f(x) - f(y) \\
&  \overset{\text{Lemma} \  \ref{lem3} }{\leq}  \    \frac{\mu_2}{2} \ d(x,P)^2 -  \frac{\mu_1}{2} \ d(y,P)^2.
\end{align*}
This completes the proof.

\paragraph{Proof of Lemma \ref{lem:grad1}}

Since, $\bar{y} \in P$, from the definition we have,
\begin{align*}
   \langle  \bar{y}-y, \E_{\mathbb{S}}  \left[a_{i^*}(a_{i^*}^Ty-b_{i^*})^{+}\right]  \rangle \ & = \E_{\mathbb{S}} \left[(a_{i^*}^Ty-b_{i^*})^{+} \left(a_{i^*}^T\bar{y}- a_{i^*}^T y\right)\right] \\
   & \leq \E_{\mathbb{S}} \left[(a_{i^*}^Ty-b_{i^*})^{+} \left(b_{i^*}- a_{i^*}^T y\right)\right] \\
   & = - \E_{\mathbb{S}} \left[\big | (a_{i^*}^Ty-b_{i^*})^{+} \big |^2\right] \\
   &= - 2f(y) \overset{\text{Lemma} \ \ref{lem3} }{\leq} \ \ -\mu_1 \ d(y,P)^2.
\end{align*}
Here, we used the identity $xx^+ = |x^+|^2$. This proves the Lemma.

\paragraph{Proof of Lemma \ref{lem4}}

Since, $\mathcal{P}(x) \in P$, we have
{\allowdisplaybreaks
\begin{align*}
    \E_{\mathbb{S}}  \left[ d(z,P)^2\right] &  \overset{\text{Lemma} \ \ref{lem:distance}}{\leq} \  \E_{\mathbb{S}}  \left[ \|z- \mathcal{P}(x)\|^2\right] = \E_{\mathbb{S}} \left[ \big \| x- \mathcal{P} (x) - \delta \left(a_{i^*}^Tx-b_{i^*}\right)^+ a_{i^*} \big \|^2 \right] \\
    & \overset{\eqref{def:function}}{=} \ \|x - \mathcal{P} (x) \|^2 + 2\delta^2 f(x)  + 2 \delta \ \big  \langle  \mathcal{P}(x)-x, \nabla f(x) \big \rangle \\
    & \overset{\text{Lemma} \ \ref{lem3} }{\leq} \ \|x - \mathcal{P} (x) \|^2 -2(2\delta-\delta^2) f(x) \\
    & \leq \ \|x-  \mathcal{P} (x) \|^2 - (2\delta-\delta^2) \ \mu_1  \|x-  \mathcal{P} (x) \|^2 = h(\delta) \ d(x,P)^2.
\end{align*}}
Here, we used the lower bound of the expected value from Lemma \ref{lem3}.

\paragraph{Proof of Theorem \ref{th:0}}

Since, $\phi_1, \phi_2 \geq 0$, the largest root $\phi$ of equation $\phi^2 + \phi_1 \phi - \phi_2 = 0$ can written as
\begin{align*}
    \phi = \frac{-\phi_1+ \sqrt{\phi_1^2+4\phi_2}}{2} \geq  \frac{-\phi_1+ \phi_1}{2} = 0.
\end{align*}
Then using the given recurrence we have,
\begin{align*}
    G_{k+1} + \phi G_k & \leq (\phi+\phi_1) G_k + \phi_2 G_{k-1} \\
    & = (\phi+\phi_1) \left(G_{k} + \phi G_{k-1}\right) \\
    & \vdots \\
    & \leq (\phi+\phi_1)^k \left(G_{1} + \phi G_{0}\right) \\
    & = (\phi+\phi_1)^k (1+\phi) G_0.
\end{align*}
This proves the first part. Also note that since $\phi_1 + \phi_2 < 1$, we have,
\begin{align*}
    \phi + \phi_1 & = \frac{\phi_1+ \sqrt{\phi_1^2+4\phi_2}}{2}  < \frac{\phi_1+ \sqrt{\phi_1^2+4(1-\phi_1)}}{2}  = \frac{\phi_1+2-\phi_1}{2} = 1.
\end{align*}

For the second part, notice that from the recurrence inequality, we can deduce the following matrix inequality:
{\allowdisplaybreaks
\begin{align}
\label{eq:th00}
 \begin{bmatrix}
G_{k+1} \\
G_{k}
\end{bmatrix}  \leq  \begin{bmatrix}
\phi_1^2+\phi_2 & \phi_1 \phi_2 \\
\phi_1 & \phi_2
\end{bmatrix}  \begin{bmatrix}
G_{k-1} \\
G_{k-2}
\end{bmatrix}.
\end{align}}
The Jordan decomposition of the matrix in the above expression is given by,
\begin{align}
\label{jordan}
\begin{bmatrix}
\phi_1^2+\phi_2 & \phi_1 \phi_2 \\
\phi_1 & \phi_2
\end{bmatrix}  =  \begin{bmatrix}
-\phi &  \phi + \phi_1 \\
1 & 1 
\end{bmatrix} \begin{bmatrix}
\phi^2 & 0 \\
0 & \rho^2 
\end{bmatrix}   \begin{bmatrix}
\frac{-1}{\phi_1+2 \phi} & \frac{1}{2}+ \frac{\phi_1}{2(\phi_1+2 \phi)} \\
\frac{1}{\phi_1+2 \phi} & \frac{1}{2}-\frac{\phi_1}{2(\phi_1+2 \phi)}
\end{bmatrix}.
\end{align}

Next, we discuss two possible cases of values of $k$. Also, we substituted $\phi_2 = \phi(\phi+\phi_1)$ in the Jordan decomposition of equation \eqref{jordan}.
\paragraph{Case 1: $k$ even}

{\allowdisplaybreaks
\begin{align}
\label{eq:th01}
 \begin{bmatrix}
G_{k+1} \\
G_{k} 
\end{bmatrix} & \overset{\eqref{eq:th00}}{\leq}   \begin{bmatrix}
\phi_1^2+\phi^2 + \phi \phi_1 & \phi^2 \phi_1 + \phi \phi_1^2 \\
\phi_1 & \phi^2 + \phi \phi_1 \nonumber 
\end{bmatrix}  \begin{bmatrix}
G_{k-1} \\
G_{k-2}
\end{bmatrix} \nonumber \\
&  \quad \vdots \nonumber \\
& \leq \  \begin{bmatrix}
\phi_1^2+\phi^2 + \phi \phi_1 & \phi^2 \phi_1 + \phi \phi_1^2 \\
\phi_1 & \phi^2 + \phi \phi_1 \nonumber 
\end{bmatrix}^{\frac{k}{2}}  \begin{bmatrix}
G_{1} \\
G_{0}
\end{bmatrix} \\
& \overset{\eqref{jordan} }{=}   \begin{bmatrix}
-\phi &  \phi + \phi_1 \\
1 & 1 
\end{bmatrix} \begin{bmatrix}
\phi^k & 0 \\
0 & \rho^k 
\end{bmatrix}   \begin{bmatrix}
\frac{-1}{\phi_1+2 \phi} & \frac{1}{2}+ \frac{\phi_1}{2(\phi_1+2 \phi)} \\
\frac{1}{\phi_1+2 \phi} & \frac{1}{2}-\frac{\phi_1}{2(\phi_1+2 \phi)}
\end{bmatrix} \begin{bmatrix}
G_{0} \\
G_{0} 
\end{bmatrix} \nonumber \\
& = \begin{bmatrix}
(1+\phi) \rho^{k+1} + (1-\phi-\phi_1) \phi^{k+1} \\
(1+\phi) \rho^{k} - (1-\phi-\phi_1) \phi^{k}
\end{bmatrix} \begin{bmatrix}
\frac{G_{0}}{\phi_1+ 2 \phi}
\end{bmatrix} \nonumber \\
& \overset{\eqref{def:0} }{=}   \begin{bmatrix}
R_1 \rho^{k+1}+ R_2 \phi^{k+1}  \\[6pt]
R_1 \rho^{k}- R_2 \phi^{k} 
\end{bmatrix} \ G_0.
\end{align}}
Here, we used $G_0 = G_1$.

\paragraph{Case 2: $k$ odd}

{\allowdisplaybreaks
\begin{align}
\label{eq:th02}
 \begin{bmatrix}
G_{k+1} \\
G_{k} 
\end{bmatrix} & \overset{\eqref{eq:th00}}{\leq}   \begin{bmatrix}
\phi_1^2+\phi^2 + \phi \phi_1 & \phi^2 \phi_1 + \phi \phi_1^2 \\
\phi_1 & \phi^2 + \phi \phi_1 \nonumber 
\end{bmatrix}  \begin{bmatrix}
G_{k-1} \\
G_{k-2}
\end{bmatrix} \nonumber \\
&  \quad \vdots \nonumber \\
& \leq   \  \begin{bmatrix}
\phi_1^2+\phi^2 + \phi \phi_1 & \phi^2 \phi_1 + \phi \phi_1^2 \\
\phi_1 & \phi^2 + \phi \phi_1 \nonumber 
\end{bmatrix}^{\frac{k-1}{2}}  \begin{bmatrix}
G_{2} \\
G_{1}
\end{bmatrix} \\
& \overset{\eqref{jordan}}{=}   \begin{bmatrix}
-\phi &  \phi + \phi_1 \\
1 & 1 
\end{bmatrix} \begin{bmatrix}
\phi^{k-1} & 0 \\
0 & \rho^{k-1} 
\end{bmatrix}   \begin{bmatrix}
\frac{-1}{\phi_1+2 \phi} & \frac{1}{2}+ \frac{\phi_1}{2(\phi_1+2 \phi)} \\
\frac{1}{\phi_1+2 \phi} & \frac{1}{2}-\frac{\phi_1}{2(\phi_1+2 \phi)}
\end{bmatrix} \begin{bmatrix}
(\phi_1+\phi_2)G_{0} \\
G_{0} 
\end{bmatrix} \nonumber \\
& = \begin{bmatrix}
(\phi+\phi_1+\phi_2)\rho^{k} - (\phi-\phi_2) \phi^{k} \\
(\phi+\phi_1+\phi_2)\rho^{k-1} + (\phi-\phi_2) \phi^{k-1}
\end{bmatrix} \begin{bmatrix}
\frac{G_{0}}{\phi_1+ 2 \phi}
\end{bmatrix} \nonumber \\
& \overset{\eqref{def:0}}{=}   \begin{bmatrix}
R_3 \rho^{k}- R_4 \phi^{k}  \\[6pt]
R_3 \rho^{k-1} + R_4 \phi^{k-1} 
\end{bmatrix} \ G_0.
\end{align}}
Here, we used the inequality $G_2 \leq \phi_1 G_1 + \phi_2 G_0$. Now combining the relations from equation \eqref{eq:th01} and \eqref{eq:th02}, we can prove the second part of Theorem \ref{th:0}.

\paragraph{Proof of Theorem \ref{th:seq2}}

From the given recurrence relation, we have
{\allowdisplaybreaks
\begin{align}
\label{t2}
 \begin{bmatrix}
H_{k+1} \\
F_{k+1}
\end{bmatrix}  & \leq  \begin{bmatrix}
\Pi_1 & \Pi_2 \\
\Pi_3  & \ \Pi_4
\end{bmatrix} \begin{bmatrix}
H_{k} \\
F_{k}
\end{bmatrix}   \leq  \begin{bmatrix}
\Pi_1 & \Pi_2 \\
\Pi_3  & \ \Pi_4
\end{bmatrix}^k  \begin{bmatrix}
H_{1} \\
F_{1}
\end{bmatrix}. 
\end{align}}
Using the definitions of \eqref{t1}, we can write the Jordan decomposition of the above matrix as follows
\begin{align}
\label{t3}
\begin{bmatrix}
\Pi_1 & \Pi_2 \\
\Pi_3  & \ \Pi_4
\end{bmatrix}  =  \begin{bmatrix}
\Gamma_2 & \Gamma_1 \\
1  & \ 1
\end{bmatrix} \begin{bmatrix}
\rho_1 & 0 \\
0 & \rho_2
\end{bmatrix}    \begin{bmatrix}
-\Gamma_3 & \Gamma_1 \Gamma_3 \\
\Gamma_3  & \ \Gamma_2 \Gamma_3
\end{bmatrix}.
\end{align}
Now, substituting the matrix decomposition into equation \eqref{t2} and simplifying we have
{\allowdisplaybreaks
\begin{align}
\label{t4}
 \begin{bmatrix}
H_{k+1} \\
F_{k+1} 
\end{bmatrix} & \leq \   \begin{bmatrix}
\Pi_1 & \Pi_2 \\
\Pi_3  & \ \Pi_4
\end{bmatrix}^k \begin{bmatrix}
H_{1} \\
F_{1}
\end{bmatrix}  \overset{\eqref{t3} }{=}   \begin{bmatrix}
\Gamma_2 & \Gamma_1 \\
1  & \ 1
\end{bmatrix} \begin{bmatrix}
\rho_1^k & 0 \\
0 & \rho_2^k
\end{bmatrix}    \begin{bmatrix}
-\Gamma_3 & \Gamma_1 \Gamma_3 \\
\Gamma_3  & \ \Gamma_2 \Gamma_3
\end{bmatrix} \begin{bmatrix}
H_{1} \\
F_1 
\end{bmatrix} \nonumber \\
& = \begin{bmatrix}
\Gamma_2 \Gamma_3 (\Gamma_1-1) \ \rho_1^{k}+ \Gamma_1 \Gamma_3 (\Gamma_2+1)\ \rho_2^{k} \\[6pt]
\Gamma_3 (\Gamma_1-1) \ \rho_1^{k}+ \Gamma_3 (\Gamma_2+1)\ \rho_2^{k}
\end{bmatrix} \ \begin{bmatrix}
H_{1} \\
F_1 
\end{bmatrix}.
\end{align}}
Since $\Pi_1, \Pi_2, \Pi_3, \Pi_4 \geq 0$, one can easily verify that $\Gamma_1,  \Gamma_3 \geq 0$. Now, it remains to show that $ 0 \leq \rho_1 \leq \rho_2 < 1$. To show that, first note that
\begin{align}
    \label{t5}
    (\Pi_1 - \Pi_4)^2 + 4\Pi_2 \Pi_3 & \overset{\eqref{t0} }{<}  (\Pi_1 - \Pi_4)^2 + 4 - 4 \Pi_1 \Pi_4 - 4(\Pi_1 + \Pi_4) \nonumber \\ 
    & = \left(2-\Pi_1 - \Pi_4\right)^2.
\end{align}
Now, we have,
\begin{align*}
     \rho_1  & = \frac{1}{2}  \left[\Pi_1+\Pi_4 - \sqrt{(\Pi_1-\Pi_4)^2+4\Pi_2\Pi_3}\right] \\
     & \geq \frac{1}{2}  \left[\Pi_1+\Pi_4 - \sqrt{(\Pi_1-\Pi_4)^2+4\Pi_1\Pi_4}\right] =  \frac{1}{2} \left[\Pi_1+\Pi_4 - (\Pi_1+\Pi_4)\right] =0.
\end{align*}
Moreover, since $2-\Pi_1 - \Pi_4 \geq 0 $, we have
\begin{align*}
      \rho_1 \leq   \ \rho_2  & = \frac{1}{2} \left[\Pi_1+\Pi_4 + \sqrt{(\Pi_1-\Pi_4)^2+4\Pi_2\Pi_3}\right] \\
     &  \overset{\eqref{t5} }{<}  \frac{1}{2} \left[\Pi_1+\Pi_4 + \sqrt{ \left(2-\Pi_1 - \Pi_4\right)^2 }\right] \\
     &   \overset{\eqref{t0} }{=} \frac{1}{2} \left[\Pi_1+\Pi_4 + 2-\Pi_1 - \Pi_4 \right] = 1.
\end{align*}
As $0 \leq  \rho_1 \leq \rho_2 < 1$, considering \eqref{t4} we can deduce that the sequence $\{H_k\}$ and $ \{F_k\}$ converges.

\section*{Appendix 2}

\paragraph{Proof of Theorem \ref{th:1}}

From the update formula of Algorithm \ref{alg:gskm}, we have $z_{k} = x_k - (A_{\tau_k}x_k - b_{\tau_k})^+_{i^*} a_{i^*}$ where,
\begin{align}
\label{eq:p1}
i^* = \argmax_{i \in \tau_k} \{a_i^Tx_k-b_i, 0\} \ = \ \argmax_{i \in \tau_k} (A_{\tau_k}x_k-b_{\tau_k})^+_i.
\end{align}
Similarly, the previous update formula can be written as, $z_{k-1} = x_{k-1} - (A_{\tau_{k-1}}x_{k-1} - b_{\tau_{k-1}})^+_{j^*} a_{j^*}$; where,
\begin{align}
\label{eq:p2}
j^* = \argmax_{j \in \tau_{k-1}} \{a_j^Tx_{k-1}-b_j, 0\} \ = \ \argmax_{j \in \tau_{k-1}} (A_{\tau_{k-1}}x_{k-1}-b_{\tau_{k-1}})^+_j.
\end{align}
Note that, the notation is consistent with the definition of \eqref{def:i1}. Since for any $\xi \in Q_1$, $ (1-\xi) \mathcal{P}(x_{k}) + \xi \mathcal{P}(x_{k-1}) \in P$ we have,
\begin{align}
    d(x_{k+1},P)^2 & = \  \ \big \| x_{k+1}- \mathcal{P}(x_{k+1}) \big \|^2  \nonumber \\ 
     & \overset{\text{Lemma} \ \ref{lem:distance}}{ \leq} \  \big \| x_{k+1}- (1-\xi) \mathcal{P}(x_{k}) - \xi \mathcal{P}(x_{k-1})\big \|^2 \nonumber \\
    & \overset{\eqref{eq:3a}}{=} \big \| (1-\xi) z_k + \xi z_{k-1}-(1-\xi) \mathcal{P}(x_{k}) - \xi \mathcal{P}(x_{k-1}) \big \|^2 \nonumber \\
    & \overset{\eqref{eq:3b}}{=} \big \| (1-\xi) \left\{x_k- \mathcal{P}(x_{k})- \delta \left(a_{i^*}^Tx_{k}-b_{i^*}\right)^+ a_{i^*}\right\}  \nonumber \\
    & \qquad \qquad   + \xi \left\{x_{k-1}- \mathcal{P}(x_{k-1})- \delta \left(a_{j^*}^Tx_{k-1}-b_{j^*}\right)^+ a_{j^*}\right\} \big \|^2 \nonumber \\
    & \leq (1-\xi) \ \big \| x_k- \mathcal{P}(x_{k})- \delta \left(a_{i^*}^Tx_{k}-b_{i^*}\right)^+ a_{i^*} \big \|^2 \nonumber \\
    & \qquad   \qquad  + \xi \ \big \| x_{k-1}- \mathcal{P}(x_{k-1})- \delta \left(a_{j^*}^Tx_{k-1}-b_{j^*}\right)^+ a_{j^*} \big \|^2. \label{eq:p3}
\end{align}
We used the fact that the function $\|\cdot \|^2$ is convex and $0 \leq \xi \leq 1$. Now, taking expectation in both sides of the equation \eqref{eq:p3} and using Lemma \ref{lem4}, we get the following:
\begin{align}
\label{eq:p4}
\E [d(x_{k+1},P)^2 \ | & \ \mathbb{S}_{k}, \mathbb{S}_{k-1}]  \leq (1-\xi) \ \E_{\mathbb{S}_k} \left[\big \| x_k- \mathcal{P}(x_{k})- \delta \left(a_{i^*}^Tx_{k}-b_{i^*}\right)^+ a_{i^*} \big \|^2 \right] \nonumber \\
&     + \xi \ \E_{\mathbb{S}_{k-1}} \left[\big \| x_{k-1}- \mathcal{P}(x_{k-1})- \delta \left(a_{j^*}^Tx_{k-1}-b_{j^*}\right)^+ a_{j^*} \big \|^2\right] \nonumber \\
& \overset{\text{Lemma} \ \ref{lem4}}{\leq}  (1-\xi) \ h(\delta) \ d(x_k,P)^2 + \xi \ h(\delta) \ d(x_{k-1}, P)^2,
\end{align}
where, $h(\delta)$ is defined in Lemma \ref{lem4}. Taking expectation again in equation \eqref{eq:p4} and letting $G_{k+1} = \E \left[d(x_{k+1},P)^2\right]$, we get the following:
\begin{align}
    \label{eq:p5}
    G_{k+1} \leq \phi_1 G_{k} + \phi_2 G_{k-1}.
\end{align}
Since, $\phi_1 , \phi_2 \geq 0$, $\phi_1+ \phi_2 < 1$ and $z_0 = z_1$, using first part of Theorem \ref{th:0}, we have the following:
\begin{align}
\label{eq:p500}
   \E \left[d(x_{k+1},P)^2\right] \leq (1+\phi) (\phi+ \phi_1)^k G_0 = (1+\phi) \rho^k \E \left[d(x_{0},P)^2\right].
\end{align}
Moreover, considering \eqref{eq:p500} with Lemma \ref{lem3} we get the bound of $\E[f(x_k)]$ which proves the first part of Theorem \ref{th:1}. Furthermore, using the second part of Theorem \ref{th:0} and equation \eqref{eq:p5}, we get the second part of Theorem \ref{th:1}. Now, to prove the third part first note that $\frac{1}{k} \sum \limits_{l=1}^{k} \mathcal{P}(x_l) \in P$, using Lemma \ref{lem:distance} we have
\begin{align}
\label{eq:b3}
  \E[d(\Tilde{x}_k,P)^2]  & = \E[\|\Tilde{x}_k-  \mathcal{P}(\Tilde{x}_k)\|^2]  \overset{\text{Lemma} \ \ref{lem:distance}}{ \leq}  \E \left[\Big \| \frac{1}{k} \sum \limits_{l=1}^{k} \left(x_l-\mathcal{P}(x_l)\right)\Big \|^2\right] \nonumber \\
  & \leq  \E \left[\frac{1}{k} \sum \limits_{l=1}^{k} \big \| x_l-\mathcal{P}(x_l)\big \|^2\right] = \frac{1}{k} \sum \limits_{l=1}^{k} \E[d(x_l,P)^2] \nonumber \\
  & \leq \frac{(1+\phi) d(x_0,P)^2}{k} \sum \limits_{l=1}^{k} \rho^{l-1} \leq \frac{(1+\phi) d(x_0,P)^2}{ k(1-\rho)}.
\end{align}
Furthermore, using a more simplifies version of \eqref{eq:p3} we have the following:
\begin{align*}
  G_{l+1} - G_l \leq \xi (G_l-G_{l-1})+ 2 \xi \delta (2-\delta) [f(x_l)-f(x_{l-1})] - 2 \delta (2-\delta) f(x_{l}), 
\end{align*}
for any $l \geq 1$. Summing up the above identity for $l =1,2,...,k$, we have the following:
\begin{align}
    \label{eq:b100}
    2 \delta (2-\delta) \sum \limits_{l=1}^{k} f(x_l) & \leq \xi G_0 + G_1- \xi G_k- G_{k+1}+ 2 \xi \delta (2-\delta) [f(x_k)-f(x_{0})] \nonumber \\
    & \leq (1+\xi) G_0 + \xi [2 \eta f(x_k)-G_k] \nonumber \\
    & \leq (1+\xi) G_0 + \xi [\eta \mu_2 -1] G_k \leq (1+\xi) d(x_0,P)^2,
\end{align}
where, $\eta = 2\delta - \delta^2$. We used the non-negativity of the sequences $G_k$ and $ f(x_k)$. We also used the upper bound from Lemma \ref{lem3}. Then, we get
\begin{align*}
  \E[f(\Tilde{x}_k)]   & \leq  \E \left[ \frac{1}{k} \sum \limits_{l=1}^{k} f(x_l)\right] = \frac{1}{k} \sum \limits_{l=1}^{k} \E[f(x_l)] \leq  \frac{(1+\xi) \ d(x_0,P)^2}{2 \delta k (2-\delta)}.
\end{align*}
This proves the second part of Theorem \ref{th:1}.


\paragraph{Proof of Theorem \ref{th:cesaro}}

For any natural number $l \geq 1$ define, $\vartheta_l = \frac{\xi}{1+\xi}[x_{l-1}-x_l-\delta (a_{j^*}^Tx_{l-1}-b_{j^*})^+a_{j^*}]$,  $ \ \Delta_l = x_l + \vartheta_l$ and $\chi_l = \|x_l+\vartheta_l-\mathcal{P}(\Delta_l)\|^2$, then using the update formulas \eqref{eq:3a} and \eqref{eq:3b}, we have
\begin{align*}
    x_{l+1} + \vartheta_{l+1} \overset{\eqref{eq:3a} \ \& \ \eqref{eq:3b}}{=} x_l + \vartheta_l - \frac{\delta}{1+\xi} \left(a_{i^*}^Tx_l-b_{i^*}\right)^+ a_{i^*},
\end{align*}
here, the index $i^{*}$ and $j^*$ are defined based on \eqref{def:i1} respectively for the sequences $x_l$ and $x_{l-1}$.  Using the above relation, we can write
\begin{align}
\label{ces:1}
    \chi_{l+1} & = \|x_{l+1}+\vartheta_{l+1}-\mathcal{P}(\Delta_{l+1})\|^2  \overset{\text{Lemma} \ \ref{lem:distance}}{ \leq}  \|x_{l+1}+\vartheta_{l+1}-\mathcal{P}(\Delta_l)\|^2 \nonumber \\
    & = \big \| x_l + \vartheta_l - \frac{\delta}{1+\xi} \left(a_{i^*}^Tx_l-b_{i^*}\right)^+ a_{i^*}-\mathcal{P}(\Delta_l) \big \|^2 \nonumber \\
    & =   \underbrace{\|x_l+\vartheta_l-\mathcal{P}(\Delta_l)\|^2}_{= \chi_l} + \frac{\delta^2}{(1+\xi)^2} \underbrace{\|(a_{i^*}^Tx_l- b_{i^*})^+a_{i^*}\|^2}_{J_1} \nonumber \\
    & - \frac{2 \delta}{1+\xi}   \underbrace{\big \langle x_l+\vartheta_l-\mathcal{P}(\Delta_l) \ ,\  a_{i^*} (a_{i^*}^Tx_l-b_{i^{*}})^{+} \big \rangle }_{J_2} \nonumber \\
    & = \chi_l + \frac{\delta^2}{(1+\xi)^2} J_1- \frac{2 \delta}{1+\xi} J_2.
\end{align}
Taking expectation with respect to $\mathbb{S}_l$ we have,
\begin{align}
  \label{ces:2}
  \frac{\delta^2}{(1+\xi)^2} \E_{\mathbb{S}_l} [J_1] \overset{\eqref{def:function} }{=}  \frac{2\delta^2}{(1+\xi)^2} f(x_l).
\end{align}
Similarly, we can simplify the third term of \eqref{ces:1} as
\begin{align}
  \label{ces:3}
  & - \frac{2\delta}{1+\xi}  \E_{\mathbb{S}_l} [J_2] \nonumber \\
  & \overset{\eqref{def:function} }{=}  -\frac{2\delta}{1+\xi} \big \langle x_l-\mathcal{P}(\Delta_l)  , \nabla f(x_l) \big \rangle - \frac{2\delta \xi}{(1+\xi)^2} \big \langle x_{l-1}-x_l - \delta \nabla f(x_{l-1}) ,  \nabla f(x_l) \big \rangle \nonumber \\
  & =  -\frac{2\delta}{1+\xi} \big \langle x_l-\mathcal{P}(\Delta_l)  , \nabla f(x_l) \big \rangle - \frac{2\delta \xi}{(1+\xi)^2} \big \langle x_{l-1}-x_l ,  \nabla f(x_l) \big \rangle \nonumber \\
  &  \quad \quad \quad \quad \quad  +\frac{\delta^2 \xi}{(1+\xi)^2} \left[ \|\nabla f(x_l)+\nabla f(x_{l-1})\|^2-\|\nabla f(x_l)\|^2-\|\nabla f(x_{l-1})\|^2 \right] \nonumber \\
  & \overset{\text{Lemma} \ \ref{lem:grad} \ \& \ \ref{lem:grad1}}{\leq}  - \frac{4\delta}{1+\xi} f(x_l) -  \frac{2\delta \xi}{(1+\xi)^2} \left[f(x_{l-1}) -  f(x_{l}) \right] - \frac{2\delta^2 \xi}{(1+\xi)^2} \left[f(x_{l-1}) +  f(x_{l}) \right] \nonumber \\ 
  & = -\frac{2 \delta \xi (1+\delta)}{(1+\xi)^2} f(x_{l-1}) +\frac{2 \delta \xi (1+\delta)}{(1+\xi)^2} f(x_{l}) - \frac{4 \delta  (1+\xi + \delta \xi)}{(1+\xi)^2} f(x_{l}).
\end{align} 
Using the expressions of equation \eqref{ces:2} and \eqref{ces:3} in \eqref{ces:1} and simplifying further, we have
\begin{align}
    \label{ces:4}
    \E[\chi_{l+1}] - \frac{2\delta \xi (1+\delta)}{(1+\xi)^2} f(x_l) + \varpi f(x_l) \ \leq \ \E[\chi_l] - \frac{2\delta \xi (1+\delta)}{(1+\xi)^2} f(x_{l-1}),  
\end{align}
here,
\begin{align}
 \label{ces:10}
   \varpi =   \frac{4 \delta  (1+\xi + \delta \xi)}{(1+\xi)^2}  -  \frac{2 \delta^2  }{(1+\xi)^2} =  \frac{2 \delta  (2+2 \xi + 2\delta \xi -\delta)}{(1+\xi)^2} \ > \ 0.
\end{align}
Now, taking expectation again in \eqref{ces:4} and using the tower property, we get,
\begin{align}
\label{ces:5}
    q_{l+1} + \varpi \E[f(x_l)] \leq q_l, \quad l = 1,2,3...,
\end{align}
where, $q_l = \E[\chi_l] - \frac{2\delta \xi (1+\delta)}{(1+\xi)^2} \E[f(x_{l-1})]$. Summing up \eqref{ces:5} for $l=1,2,...,k$ we get
\begin{align}
    \label{ces:6}
    \sum \limits_{l=1}^{k} \E [f(x_l)] \ \leq \ \frac{q_1-q_{k+1}}{\varpi} \ \leq \ \frac{q_1}{\varpi}.
\end{align}
Now, using Jensen's inequality, we have
\begin{align*}
    \E \left[f(\bar{x_k})\right] = \E \left[f\left(\sum \limits_{l=1}^{k} \frac{x_l}{k}\right)\right] \ \leq \ \E \left[\frac{1}{k} \sum \limits_{l=1}^{k}f(x_l)\right] \ = \ \frac{1}{k}  \sum \limits_{l=1}^{k} \E [f(x_l)] \ \overset{\eqref{ces:6}}{\leq}  \frac{q_1}{\varpi k}.
\end{align*}
Since, $x_0=x_1$, we have $\vartheta_1 = \frac{-\delta \xi}{1+\xi} (a_{i^*}^Tx_0-b_{i^*})^+ a_{i^*} $. Furthermore, 
\begin{align}
    \label{ces:7}
    \E[\chi_1] & = \E \left[\|x_1+\vartheta_1 -\mathcal{P}(\Delta_1)\|^2 \right]  \overset{\text{Lemma} \ \ref{lem:distance}}{ \leq} \E \left[\|x_1+\vartheta_1 -\mathcal{P}(x_0)\|^2 \right] \nonumber \\
    & = \E \left[\|x_0 -\mathcal{P}(x_0)-\frac{\delta \xi}{1+\xi} (a_{i^*}^Tx_0-b_{i^*})^+ a_{i^*} \|^2 \right] \nonumber \\
    & = \|x_0 -\mathcal{P}(x_0)\|^2 + \frac{\delta^2\xi^2}{(1+\xi)^2} \E [|(a_{i^*}^Tx_0-b_{i^*})^+|^2] \nonumber \\
    & \quad  \quad \quad \quad  \quad \quad - \frac{2\delta \xi}{1+\xi} \langle x_0 -\mathcal{P}(x_0) , \E [(a_{i^*}^Tx_0-b_{i^*})^+ a_{i^*}] \rangle \nonumber \\
    & \overset{\text{Lemma} \ \ref{lem3} }{\leq} \|x_0 -\mathcal{P}(x_0)\|^2 + \frac{2\delta^2\xi^2}{(1+\xi)^2} f(x_0)-\frac{2\delta \xi \mu_2}{1+\xi} \|x_0 -\mathcal{P}(x_0)\|^2.
\end{align}
Now, from our construction we get
\begin{align*}
   q_1 = \E[\chi_1] - \frac{2\delta \xi (1+\delta)}{(1+\xi)^2} \E[f(x_{0})] \leq  (1-\frac{2\delta \xi \mu_2}{1+\xi}) \ d(x_0,P)^2+ \frac{2\delta\xi(\delta \xi -1-\delta )}{(1+\xi)^2} f(x_0).
\end{align*}
Substituting the values of $\varpi$ and $q_1$ in the expression of $\E \left[f(\bar{x_k})\right] $, we have the following:
\begin{align*}
    \E \left[f(\bar{x}_k)\right] \leq \frac{ (1+\xi) (1+\xi-2\delta \xi \mu_2) \ d(x_0,P)^2+ 2 \xi \delta (\delta \xi -\delta  -1) f(x_0)}{2 \delta k \left(2+2 \xi + 2 \delta \xi -\delta\right)}.
\end{align*}

\paragraph{Proof of Theorem \ref{th:2}}

Since, the term $\| x_{k+1}-  \mathcal{P}(x_{k+1})  \|$ is constant under  From the update formula of the GSKM algorithm, we get,
\begin{align}
\label{1}
   \E [ \| x_{k+1}-  \mathcal{P}(x_{k+1})  & \| \ | \ \mathbb{S}_{k+1}, \mathbb{S}_k]  =  \E [ \| x_{k+1}-  \mathcal{P}(x_{k+1})  \| \ | \  \mathbb{S}_k] \nonumber \\ 
   & \overset{\text{Lemma} \ \ref{lem:distance}}{ \leq}\  \  \E_{\mathbb{S}_k} [\| x_{k+1}- \mathcal{P}(x_{k})  \| ] \nonumber \\
    & = \E_{\mathbb{S}_k} [\|z_k-\mathcal{P}(x_k)- \xi (z_k-z_{k-1}) \|] \nonumber \\
    & \leq  \E_{\mathbb{S}_k}[\|z_k-\mathcal{P}(x_k)\|] + |\xi| \E_{\mathbb{S}_k} [\|z_k-z_{k-1}\|] \nonumber \\
    & \leq \left\{\E_{\mathbb{S}_k}[ \|z_k-\mathcal{P}(x_k)\|^2]\right\}^{\frac{1}{2}} +  |\xi| \|z_k-z_{k-1}\| \nonumber \\
     & \overset{\text{Lemma} \ \ref{lem4} }{\leq}  \sqrt{h(\delta)} \  \|x_k-\mathcal{P}(x_k)\| + |\xi| \|z_k-z_{k-1}\|.
\end{align}
We performed the two expectations in order, from the innermost to the outermost. Now, taking expectation in \eqref{1} and using the tower property of expectation we have,
\begin{align}
\label{3}
\E [\|x_{k+1}-\mathcal{P}(x_{k+1})\|]  & \leq  \sqrt{h(\delta)} \  \E [\|x_k-\mathcal{P}(x_k)\|] + |\xi| \ \E [\|z_k-z_{k-1}\|].
\end{align}
Similarly, using the update formula for $z_{k+1}$, we have
\begin{align}
    \label{2}
   \E [\| z_{k+1}& -z_k  \|  \ | \  \mathbb{S}_{k+1}, \mathbb{S}_k]  = \E[\E_{\mathbb{S}_{k+1}}[\| x_{k+1} - \delta \left(a_{i^*}^Tx_{k+1}-b_{i^*}\right)^+ a_{i^*}-z_k \|] \ | \ \mathbb{S}_k] \nonumber \\
    & = \E [\E_{\mathbb{S}_{k+1}} [\|-\xi (z_k-z_{k-1}) - \delta \left(a_{i^*}^Tx_{k+1}-b_{i^*}\right)^+ a_{i^*}\|]  \ | \ \mathbb{S}_k] \nonumber \\
    & \leq  |\xi| \ \|z_k-z_{k-1}\| + \delta  \E [\E_{\mathbb{S}_{k+1}}[|(a_{i^*}^Tx_{k+1}-b_{i^*})^+ |]  \ | \ \mathbb{S}_k] \nonumber \\
    & \leq  |\xi| \ \|z_k-z_{k-1}\| + \delta  \E [\left\{\E_{\mathbb{S}_{k+1}}[|(a_{i^*}^Tx_{k+1}-b_{i^*})^+ |^2]\right\}^{\frac{1}{2}}  \ | \ \mathbb{S}_k] \nonumber \\
    &  \overset{\text{Lemma} \ \ref{lem3} }{\leq}  |\xi| \ \|z_k-z_{k-1}\| + \delta \sqrt{\mu_2} \ \E[\|x_{k+1}-\mathcal{P}(x_{k+1})\| \ | \ \mathbb{S}_k].
\end{align}
Taking expectation in \eqref{2} and using \eqref{3} along with the tower property, we have,
\begin{align}
    \label{4}
    \E [\|& z_{k+1}-  z_{k}\|]  \   \leq \ |\xi| \ \E [\|z_k-z_{k-1}\|] + \delta \sqrt{\mu_2} \ \E [\|x_{k+1}-\mathcal{P}(x_{k+1})\|] \nonumber \\
    &  \overset{\eqref{3}}{\leq}  |\xi| \left(1+ \delta \sqrt{\mu_2}\right)  \E [\|z_k-z_{k-1}\|] + \delta \sqrt{\mu_2 h(\delta)}  \E [\|x_k-\mathcal{P}(x_k)\|].
\end{align}
Combining both \eqref{3} and \eqref{4},
we can deduce the following matrix inequality:
{\allowdisplaybreaks
\begin{align}
\label{5}
\E \begin{bmatrix}
\|x_{k+1}-\mathcal{P}(x_{k+1})\|  \\[6pt]
\|z_{k+1}-z_k\| 
\end{bmatrix}  & \leq  \begin{bmatrix}
\sqrt{h(\delta)} & |\xi| \\
\delta \sqrt{\mu_2 h(\delta)}  &  \ |\xi| \left(1+ \delta \sqrt{\mu_2}\right)
\end{bmatrix} \E \begin{bmatrix}
 \|x_k-\mathcal{P}(x_k)\| \\
 \|z_k-z_{k-1}\|
\end{bmatrix}.
\end{align}}
Now, from the definition, it can be easily checked that $\Pi_1, \Pi_2, \Pi_3, \Pi_4 \geq 0$. Since, $\xi \in Q_2$, we have
\begin{align}
    \label{7}
    \Pi_2 \Pi_3 - \Pi_1 \Pi_4 & = |\xi| \delta \sqrt{\mu_2 h(\delta)} - |\xi| \sqrt{h(\delta)} - |\xi| \delta \sqrt{\mu_2 h(\delta)} =  - |\xi| \sqrt{h(\delta)} \leq 0.
\end{align}
Also, we have
\begin{align}
\label{8}
\Pi_1  + \Pi_4-  \Pi_1 \Pi_4+ & \Pi_2 \Pi_3  = \sqrt{h(\delta)} + |\xi| \left(1+ \delta \sqrt{\mu_2}\right)  - |\xi| \sqrt{h(\delta)} < 1.
\end{align}
Here, in the last inequality we used the given condition. Considering \eqref{8}, we can check that $\Pi_1  + \Pi_4 < 1+|\xi| \sqrt{h(\delta)} = 1+ \min\{1, |\xi| \sqrt{h(\delta)}\} = 1+ \min\{1,\Pi_1 \Pi_4-\Pi_2 \Pi_3\}$. Also from \eqref{7}, we have $\Pi_2 \Pi_3 - \Pi_1 \Pi_4 \leq 0$, which is precisely the condition provided in \eqref{t0}. Let's define the sequences $F_k = \E [\|z_k-z_{k-1}\|]$ and $H_k = \E [\|x_k-\mathcal{P}(x_k)\|]$. Now, using Theorem \ref{th:seq2}, we have
\allowdisplaybreaks{\begin{align} \label{9}
\begin{bmatrix}
H_{k+1}  \\[6pt]
F_{k+1} 
\end{bmatrix} & \leq   \begin{bmatrix}
\Gamma_2 \Gamma_3 (\Gamma_1-1) \ \rho_1^{k}+ \Gamma_1 \Gamma_3 (\Gamma_2+1)\ \rho_2^{k} \\[6pt]
\Gamma_3 (\Gamma_1-1) \ \rho_1^{k}+ \Gamma_3 (\Gamma_2+1)\ \rho_2^{k}
\end{bmatrix} \ \begin{bmatrix}
H_{1} \\
F_1 
\end{bmatrix}.
\end{align}}
where, $\Gamma_1 , \Gamma_2, \Gamma_3, \rho_1, \rho_2$ can be derived from \eqref{t1} using the parameter choice of \eqref{6}. Note that, from the GSKM algorithm we have, $x_1 = x_0$ and $z_1 = z_0$.Therefore we can easily check that, $F_1 = \E [\|z_1-z_{0}\|] = 0$ and $H_1 = \E [\|x_1-\mathcal{P}(x_1)\|] = \E [\|x_0-\mathcal{P}(x_0)\|] = \|x_0-\mathcal{P}(x_0)\| = H_0 $. Now, substituting the values of $H_1$ and $F_1$ in \eqref{9}, we have
{\allowdisplaybreaks
\begin{align}
\label{10}
 \begin{bmatrix}
H_{k+1} \\
F_{k+1} 
\end{bmatrix}  =  \E \begin{bmatrix}
d(x_{k+1}, P)  \\[6pt]
\|z_{k+1}-z_k\| 
\end{bmatrix} \leq \begin{bmatrix}
-\Gamma_2 \Gamma_3 \ \rho_1^{k}+ \Gamma_1 \Gamma_3 \ \rho_2^{k} \\[6pt]
- \Gamma_3 \ \rho_1^{k}+ \Gamma_3 \ \rho_2^{k} 
\end{bmatrix} \ d(x_0,P).
\end{align}}
Also from Theorem \ref{th:seq2} we have, $\Gamma_1,  \Gamma_3 \geq 0$ and $ 0 \leq \rho_1 \leq \rho_2 < 1$. Which proves the Theorem.

\paragraph{Proof of Theorem \ref{th:3}}

Note that, since $Ax \leq b$ is feasible, then from Lemma \ref{lem:skm4}, we know that there is a feasible solution $x^*$ with $|x^{*}_j| \leq \frac{2^{\sigma}}{2n}$ for $j = 1, ..., n$. Thus, we have,
\begin{align}
\label{eq:th40}
   d(x_0,P) = \|x_0-\mathcal{P}(x_0)\| \ \leq \ \|x^*\| \ \leq \ \frac{2^{\sigma -1}}{\sqrt{n}},
\end{align}
as $x_0 = 0$. Then if the system $Ax \leq b$ is infeasible, by using Lemma \ref{lem:skm1}, we have,
\begin{align*}
 \theta (x) \ \geq \ 2^{1-\sigma}.
\end{align*}
This implies when GSKM runs on the system $Ax \leq b$, the system is feasible when $\theta (x) < 2^{1-\sigma}$. Furthermore, since every point of the feasible region $P$ is inside the half-space defined by $\Tilde{H}_i = \{x \ | \ a_i^T x \leq b_i\}$ for all $i = 1,2,...,m$, we have the following:
\begin{align}
\label{eq:th41}
  \theta(x) \ = \   \left[\max_{i}\{a_i^Tx-b_i\}\right]^{+} \ \leq \ \|a_i^T(x-\mathcal{P}(x))\|   \ \leq \ d(x,P).
\end{align}
Then, for $\xi \in Q_1$ whenever the system $Ax \leq b$ is feasible, we have,
\begin{align}
\label{eq:th420}
  \E \left[\theta(x_k)\right]  \overset{\eqref{eq:th41}}{\leq}  \E \left[d(x_{k+1},P)\right]  & \leq  \sqrt{\E \left[d(x_{k+1},P)^2\right]}   \overset{\text{Theorem} \ \ref{th:1}}{\leq} \sqrt{1+\phi}  \rho^{\frac{k}{2}} \ d(x_0,P).
\end{align}
Similarly for $\xi \in Q_2$ whenever the system $Ax \leq b$ is feasible, we have,
\begin{align}
\label{eq:th421}
  \E \left[\theta(x_k)\right]  \overset{\eqref{eq:th41}}{\leq} \ \E \left[d(x_{k+1},P)\right]  &  \overset{\text{Theorem} \ \ref{th:2}}{\leq} \sqrt{1+\phi} \ \rho_2^{k} \ d(x_0,P).
\end{align}
Take, $\bar{\rho} = \max\{\rho, \rho_2^2\}$ \footnote{Note that, since $\Gamma_1 \Gamma_2 \leq 0$ and $\Gamma_1 \Gamma_3 \leq 1 \leq \sqrt{(1+\phi)}$, from Theorem \ref{th:2} we have $\E[d(x_{k+1},P)] \leq \sqrt{(1+\phi)} \ \rho_2^k \ d(x_0,P)$ for any $\xi \in Q_2$.}. Now combining \eqref{eq:th420} and \eqref{eq:th421}, for any $\xi \in Q = Q_1 \cup Q_2$, whenever the system $Ax \leq b$ is feasible, we have,
\begin{align}
\label{eq:th42}
  \E \left[\theta(x_k)\right] \overset{\eqref{eq:th420} \ \& \ \eqref{eq:th421}}{\leq} \sqrt{1+\phi} \ \bar{\rho}^{\frac{k}{2}} \ d(x_0,P) \overset{\eqref{eq:th40}}{\leq} \ \sqrt{1+\phi} \ \bar{\rho}^{\frac{k}{2}} \ \frac{2^{\sigma -1}}{\sqrt{n}}.
\end{align}
Here, we used Theorems \ref{th:1} \& \ref{th:2} and the identities from equations \eqref{eq:th40} \& \eqref{eq:th41}. Now, for detecting feasibility we need to have, $\E [\theta(x_k)] < 2^{1-\sigma}$. That gives us,
\begin{align*}
    \sqrt{1+\phi} \ \bar{\rho}^{\frac{k}{2}} \ \frac{2^{\sigma -1}}{\sqrt{n}} < 2^{1-\sigma}.
\end{align*}
Simplifying the above identity further we get the following lower bound for $k$:
\begin{align*}
   k \ > \ \frac{4 \sigma - 4 -\log n + \log (1+\phi)}{\log \left(\frac{1}{\bar{\rho}}\right)}. 
\end{align*}
Moreover, if the system $Ax \leq b$ is feasible, then the probability of not having a certificate of feasibility is bounded as follows,
\begin{align*}
  p = \mathbb{P} \left(\theta(x_k) \geq 2^{1-\sigma}\right) \ \leq \ \frac{\E \left[\theta(x_k)\right]}{2^{1-\sigma}} \ < \ \sqrt{\frac{1+\phi}{n}} \ 2^{2\sigma -2} \ \bar{\rho}^{\frac{k}{2}}.
\end{align*}
Here, we used the Markov's inequality $\mathbb{P} (x \geq t) \leq \frac{\E[x]}{t}$. This completes the proof of Theorem \ref{th:3}.

\section*{Appendix 3}

\paragraph{Proof of Theorem \ref{th:4}}

From the update formula of the PASKM algorithm, we get,
\begin{align}
\label{20}
   & \E_{\mathbb{S}_{k}} [\|  v_{k+1}-  \mathcal{P}(v_{k+1})  \|^2]  \nonumber \\
   & \overset{\text{Lemma} \ \ref{lem:distance}}{ \leq}\  \  \E_{\mathbb{S}_{k}} [\| v_{k+1}- \omega \mathcal{P}(v_{k})-(1-\omega) \mathcal{P}(y_{k})  \|^2]  \nonumber \\
    & =  \E_{\mathbb{S}_{k}} [\|\omega (v_{k}-  \mathcal{P}(v_{k})) + (1-\omega) (y_{k}-  \mathcal{P}(y_{k})) - \gamma \left(a_{i^*}^Ty_{k}-b_{i^*}\right)^+ a_{i^*}\|^2] \nonumber \\
    & = \E_{\mathbb{S}_{k}} [\|\omega (v_{k}-  \mathcal{P}(v_{k})) + (1-\omega) (y_{k}-  \mathcal{P}(y_{k}))\|^2] + \gamma^2 \E_{\mathbb{S}_{k}} [|(a_{i^*}^Ty_{k}-b_{i^*})^+|^2]  \nonumber \\ 
    & - 2 \gamma (1-\omega) \big \langle y_k-\mathcal{P}(y_{k}), \E_{\mathbb{S}_{k}} [(a_{i^*}^Ty_{k}-b_{i^*})^+ a_{i^*}] \big \rangle   \nonumber \\ 
    & - 2 \gamma \omega \big \langle v_k-\mathcal{P}(v_{k}), \E_{\mathbb{S}_{k}} [(a_{i^*}^Ty_{k}-b_{i^*})^+ a_{i^*}] \big \rangle  \nonumber \\
    & \leq \omega  \| v_{k}-  \mathcal{P}(v_{k})  \|^2 + (1-\omega) \| y_{k}-  \mathcal{P}(y_{k}) \|^2 + \gamma^2 \E_{\mathbb{S}_{k}} [|(a_{i^*}^Ty_{k}-b_{i^*})^+|^2] \nonumber \\
    & - 2 \gamma (1-\omega) \E_{\mathbb{S}_{k}} [|(a_{i^*}^Ty_{k}-b_{i^*})^+|^2] + \omega \gamma \E_{\mathbb{S}_{k}} [|(a_{i^*}^Ty_{k}-b_{i^*})^+|^2] + \omega \gamma \| v_{k}-  \mathcal{P}(v_{k})  \|^2 \nonumber \\
    & = \omega (1+\gamma)  \| v_{k}-  \mathcal{P}(v_{k})  \|^2 + (1-\omega) \| y_{k}-  \mathcal{P}(y_{k}) \|^2 + 2\gamma (\gamma + 3 \omega -2) f(y_k) \nonumber \\
    & \leq  \omega (1+\gamma)  \| v_{k}-  \mathcal{P}(v_{k})  \|^2 + \left\{1-\omega+\gamma \mu_1 (\gamma + 3 \omega -2)\right\} \| y_{k}-  \mathcal{P}(y_{k}) \|^2.
\end{align}
Here, we used the condition $\gamma + 3 \omega -2 \leq 0$. Similarly, using the update formula for $y_{k+1}$, we have
\begin{align}
    \label{21}
    \E_{\mathbb{S}_{k}}[\| & y_{k+1}- \mathcal{P}(y_{k+1}) \|^2] \nonumber \\  
    & \overset{\text{Lemma} \ \ref{lem:distance}}{ \leq}\  \E_{\mathbb{S}_{k}}[\| \alpha (v_{k+1}- \mathcal{P}(v_{k+1})) + (1-\alpha) (x_{k+1} - \mathcal{P}(y_{k}))  \|^2] \nonumber \\
    & \leq \alpha \E_{\mathbb{S}_{k}}[\|v_{k+1}- \mathcal{P}(v_{k+1})\|^2] + (1-\alpha) \E_{\mathbb{S}_{k}}[\|x_{k+1} - \mathcal{P}(y_{k})\|^2] \nonumber \\
    & \overset{\text{Lemma} \ \ref{lem4}}{ \leq}\ \alpha \E_{\mathbb{S}_{k}}[\|v_{k+1}- \mathcal{P}(v_{k+1})\|^2] + (1-\alpha) h(\delta) \|y_{k} - \mathcal{P}(y_{k})\|^2.
\end{align}
Following Theorem \ref{th:seq2}, let us define the sequences $H_k = E [\|v_k-\mathcal{P}(v_k)\|^2]$ and $F_k = \E [\|y_k-\mathcal{P}(y_k)\|^2]$. The goal is to prove that $H_k$ and $F_k$ satisfy the condition \eqref{t0}. Now, taking expectation in \eqref{20} and using the tower property of expectation we have,
\begin{align}
\label{22}
H_{k+1} \ & \leq \omega (1+\gamma)  H_k + \left\{1-\omega+\gamma \mu_1 (\gamma + 3 \omega -2)\right\} F_k. 
\end{align}
Similarly, taking expectation in \eqref{21} and using \eqref{22} along with the tower property of expectation we have,
\begin{align}
    \label{24}
    & F_{k+1} \  \leq  \alpha H_{k+1} + (1-\alpha) h(\delta) F_k \nonumber \\
    & \leq \ \alpha \omega (1+\gamma)  H_k + \{ (1-\alpha) h(\delta) + \alpha (1-\omega) + \alpha \gamma \mu_1 (\gamma + 3 \omega -2)\} F_k.
\end{align}
Combining both \eqref{22} and \eqref{24},
we can deduce the following matrix inequality:
{\allowdisplaybreaks
\begin{align}
\label{25}
 \begin{bmatrix}
H_{k+1} \\
F_{k+1}
\end{bmatrix}  & \leq  \begin{bmatrix}
\Pi_1 & \Pi_2 \\
\Pi_3  & \ \Pi_4
\end{bmatrix} \begin{bmatrix}
H_{k} \\
F_{k}
\end{bmatrix}   \leq  \begin{bmatrix}
\Pi_1 & \Pi_2 \\
\Pi_3  & \ \Pi_4
\end{bmatrix}^{k+1}  \begin{bmatrix}
H_{0} \\
F_{0}
\end{bmatrix}. 
\end{align}}
Here, we use the fact that $\Pi_1, \Pi_2, \Pi_3, \Pi_4 \geq 0$. Now we will use Theorem \ref{th:seq2} to simplify the expression of \eqref{25}. Before we can use Theorem \ref{th:seq2}, we need to make sure the sequences $H_k$ and $F_k$ satisfy the condition of \eqref{t0}. From the definition, we have
\begin{align}
    \label{28}
    \Pi_2 & \Pi_3  - \Pi_1 \Pi_4  = \alpha \omega (1-\omega)(1+\gamma) + \alpha \omega \gamma \mu_1 (1+\gamma)(\gamma+3w-2) \nonumber \\
    & -  \omega  h(\delta) (1-\alpha) (1+\gamma) -\alpha \omega (1-\omega)(1+\gamma) - \alpha \omega \gamma \mu_1 (1+\gamma)(\gamma+3w-2) \nonumber \\ 
    & = -  \omega  h(\delta) (1-\alpha) (1+\gamma)  \leq 0. 
\end{align}
Also, we have
\begin{align}
\label{29}
\Pi_1   + \Pi_4-  \Pi_1 \Pi_4+  & \Pi_2 \Pi_3 =  \omega (1+\gamma) + h(\delta) (1-\alpha) + \alpha (1-\omega) \nonumber \\
&  + \alpha \gamma \mu_1 (\gamma+3w-2)-  \omega  h(\delta) (1-\alpha) (1+\gamma) < 1.
\end{align}
Here, in the last inequality, we used the given condition. Considering \eqref{29}, we can check that $\Pi_1  + \Pi_4 < 1+\omega  h(\delta) (1-\alpha) (1+\gamma) = 1+ \min\{1, \omega  h(\delta) (1-\alpha) (1+\gamma)\} = 1+ \min\{1,\Pi_1 \Pi_4-\Pi_2 \Pi_3\}$. Also from \eqref{28}, we have $\Pi_2 \Pi_3 - \Pi_1 \Pi_4 \leq 0$, which is precisely the condition provided in \eqref{t0}. Now, using Theorem \ref{th:seq2}, we have
\allowdisplaybreaks{\begin{align} \label{30}
\begin{bmatrix}
H_{k+1}  \\[6pt]
F_{k+1} 
\end{bmatrix} & \leq   \begin{bmatrix}
\Gamma_2 \Gamma_3 (\Gamma_1-1) \ \rho_1^{k+1}+ \Gamma_1 \Gamma_3 (\Gamma_2+1)\ \rho_2^{k+1} \\[6pt]
\Gamma_3 (\Gamma_1-1) \ \rho_1^{k+1}+ \Gamma_3 (\Gamma_2+1)\ \rho_2^{k+1}
\end{bmatrix} \ \begin{bmatrix}
H_{0} \\
F_0 
\end{bmatrix}.
\end{align}}
where, $\Gamma_1 , \Gamma_2, \Gamma_3, \rho_1, \rho_2$ can be derived from \eqref{t1} using the given parameter. Note that, from the PASKM algorithm we have, $x_0 = v_0 = y_0$. Therefore we can easily check that, $ H_0 =  \|v_0-\mathcal{P}(v_0)\|^2 = \|y_0-\mathcal{P}(y_0)\|^2 = F_0 $. Now, substituting the values of $H_0$ and $F_0$ in \eqref{30}, we have
{\allowdisplaybreaks
\begin{align}
\label{eq:10000}
\E \begin{bmatrix}
d(v_{k+1}, P)^2  \\[6pt]
d(y_{k+1}, P)^2  
\end{bmatrix} \leq  \begin{bmatrix}
\Gamma_2 \Gamma_3 (\Gamma_1-1) \ \rho_1^{k+1}+ \Gamma_1 \Gamma_3 (\Gamma_2+1)\ \rho_2^{k+1} \\[6pt]
\Gamma_3 (\Gamma_1-1) \ \rho_1^{k+1}+ \Gamma_3 (\Gamma_2+1)\ \rho_2^{k+1}
\end{bmatrix} d(y_0,P)^2.
\end{align}}
Also from Theorem \ref{th:seq2} we have, $\Gamma_1,  \Gamma_3 \geq 0$ and $ 0 \leq \rho_1 \leq \rho_2 < 1$. Which proves the the first part of the Theorem. Now, considering Lemma \ref{lem3}, we get 
\begin{align}
    \label{5001}
   \E[ f(x_{k+1})] \  \leq  \ \frac{\mu_2}{2} \ \E[\|y_{k+1}-\mathcal{P}(y_{k+1})\|^2] = \ \frac{\mu_2}{2} \ \E[d(y_{k+1},P)^2].
\end{align}
Now, substituting the result of \eqref{eq:10000} in \eqref{5001}, we get the second part of the Theorem.

\paragraph{Proof of Theorem \ref{th:10}}

Let us define, $\mathcal{V} = \omega \mathcal{P}(v_{k}) +(1-\omega) \mathcal{P}(y_{k}) $. Since $\mathcal{V} \in P$, using the update formula of $v_{k+1}$ from equation \eqref{eq:askm2}, we have,
\allowdisplaybreaks{
\begin{align}
\label{eq:askm4}
    d(v_{k+1},P)^2 & = \|v_{k+1} - \mathcal{P}(v_{k+1})\|^2 \overset{\text{Lemma} \ \ref{lem:distance}}{\leq}  \|v_{k+1} - \mathcal{V}\|^2  \nonumber \\
    & \overset{\eqref{eq:askm2}}{=}  \big \|\omega v_k + (1-\omega) y_k- \mathcal{V} -\gamma (a_{i^*}^Ty_k- b_{i^*})^+ a_{i^*} \big \|^2 \nonumber \\
    & =   \underbrace{ \|\omega v_k + (1-\omega) y_k- \mathcal{V} \|^2}_{I_1} + \gamma^2 \underbrace{\|(a_{i^*}^Ty_k- b_{i^*})^+a_{i^*}\|^2}_{I_2} \nonumber \\
    & - 2 \gamma  \underbrace{\big \langle \omega v_k + (1-\omega) y_k- \mathcal{V} \ ,\  a_{i^*} (a_{i^*}^Ty_k-b_{i^{*}})^{+} \big \rangle }_{I_3} \nonumber \\
    & = I_1 + \gamma^2 I_2 -2 \gamma I_3.
\end{align}}
Since $\|\cdot\|^{2}$ is a convex function and $0 < \omega < 1$, we can bound the expected first term as follows,
\begin{align}
\label{eq:askm5}
    \E_{\mathbb{S}_k} [I_1]  &= \E_{\mathbb{S}_k} \left[\|\omega v_k + (1-\omega) y_k-\mathcal{V}\|^2\right] \nonumber \\
    & = \E_{\mathbb{S}_k} \left[\|\omega v_k + (1-\omega) y_k- \omega \mathcal{P}(v_{k}) - (1-\omega) \mathcal{P}(y_{k})\|^2\right] \nonumber \\
    & \leq \omega \| v_k-\mathcal{P}(v_{k})\|^2 + (1-\omega) \|y_k-\mathcal{P}(y_{k})\|^2\nonumber \\
    & =  \omega \ d(v_k,P)^2 + (1-\omega) \ d(y_k,P)^2.
\end{align}
Taking expectation with respect to the sampling distribution in the second term of equation \eqref{eq:askm4} and using Lemma \ref{lem4} with the choice $z = x_{k+1}, \ x = y_k$ and $\eta = 2\delta-\delta^2$, we get,
\begin{align}
\label{eq:askm6}
   \gamma^2 \E_{\mathbb{S}_k} \left[\|(a_{i^*}^Ty_k- b_{i^*})^+ a_{i^*}\|^2\right] & \ = \  \gamma^2 \E_{\mathbb{S}_k} \left[|(a_{i^*}^Ty_k- b_{i^*})^+ |^2\right] \nonumber \\
    \overset{\text{Lemma} \ \ref{lem4}}{\leq} & \ \frac{\gamma^2 }{\eta} \left[d(y_k,P)^2 - \E \left[ d(x_{k+1},P)^2 \right]\right].
\end{align}

Now, taking expectation in the third term of \eqref{eq:askm4} we get,
\begin{align}
\label{eq:askm9}
    - 2  \gamma  \E_{\mathbb{S}_k} & [I_3 ]  = - 2 \gamma \big \langle \omega v_k + (1-\omega) y_k-\mathcal{V}, \E_{\mathbb{S}_k} \left[ a_{i^*} (a_{i^*}^Ty_k-b_{i^*})^{+}\right]  \big \rangle \nonumber \\
    & \overset{\eqref{eq:yk} \ \& \ \eqref{def:function} }{=} - 2 \gamma \big \langle \frac{\omega}{\alpha} \left[y_k - (1-\alpha) x_k\right] + (1-\omega) y_k-\mathcal{V},  \nabla f(y_k)  \big \rangle \nonumber \\
    & = - 2 \gamma \big \langle \frac{\omega(1-\alpha)}{\alpha}(y_k-x_k)+ y_k-\mathcal{V}, \nabla f(y_k)  \big \rangle.
\end{align}
Using Lemma \ref{lem:grad} and Lemma \ref{lem:grad1} we can simplify equation \eqref{eq:askm9} as follows,
\begin{align}
\label{eq:askm10}
  & \ -  2 \gamma  \E_{\mathbb{S}_k}  [I_3 ]  =   - 2 \gamma \big \langle \frac{\omega(1-\alpha)}{\alpha}(y_k-x_k)+ y_k-\mathcal{V}, \nabla f(y_k)  \big \rangle \nonumber \\
    & \quad =    2 \gamma \frac{\omega(1-\alpha)}{\alpha} \big \langle x_k-y_k, \nabla f(y_k)  \big \rangle + 2 \gamma \big \langle \mathcal{V}-y_k, \nabla f(y_k)  \big \rangle \nonumber \\
  & \overset{\text{Lemma} \ \ref{lem:grad} \ \& \ \ref{lem:grad1}}{\leq} \frac{\gamma \omega (1-\alpha)}{\alpha} \  d(x_k,P)^2 - \frac{ \mu_1 \gamma  \omega (1-\alpha)}{\alpha} \ d(y_k,P)^2 - 2 \mu_1 \gamma \ d(y_k,P)^2.
\end{align}
Now, substituting the values of equation \eqref{eq:askm5}, \eqref{eq:askm6} \& \eqref{eq:askm10} in equation \eqref{eq:askm4} we get the following:
\allowdisplaybreaks{
\begin{align*}
 \E [ & d(v_{k+1} ,P)^2  ]  =  I_1 + \gamma^2 \E_{\mathbb{S}_k}[I_2 ]- 2 \gamma \E_{\mathbb{S}_k}[I_3] \nonumber \\
 & = \omega \ d(v_{k},P)^2 + (1-\omega) \ d(y_k,P)^2 + \frac{\gamma^2 }{\eta} \left\{d(y_k,P)^2 - \E \left[ d(x_{k+1},P)^2 \right]\right\} \nonumber \\
 & \quad + \frac{\gamma \omega (1-\alpha)}{\alpha}  d(x_k,P)^2 - \gamma \mu_1 \left(2+\frac{ \omega (1-\alpha)}{\alpha}\right)  d(y_k,P)^2.
\end{align*}}
With further simplification, the above identity can be written as follows:
\allowdisplaybreaks{
\begin{align}
\label{eq:askm11}
 & \E \left[ d(v_{k+1},P)^2 + \frac{\gamma^2}{\eta} d(x_{k+1},P)^2 \right]  = \omega \ \left[d(v_{k},P)^2 + \frac{\gamma (1-\alpha)}{\alpha} d(x_{k},P)^2 \right] \nonumber \\
 & \quad \quad \quad \quad \quad \quad \quad + d(y_k,P)^2 \ \left\{1-\omega+ \frac{\gamma^2 }{\eta}  - 2\gamma \mu_1 - \frac{\gamma \omega \mu_1 (1-\alpha)}{\alpha}  \right\}.
\end{align}}
Now, let's choose the parameters as in equation \eqref{eq:askm12} along with $0 < \zeta  < \frac{4\eta \mu_1}{(1-\mu_1)^2} $. We can easily see that $\frac{\gamma^2 }{\eta} = \frac{ \gamma(1-\alpha)}{\alpha}$ and $ \alpha \in (0,1)$. Also note that,
\allowdisplaybreaks{\begin{align}
 2 \mu_1 \gamma \ = \ 2 \mu_1 \sqrt{\eta \zeta \mu_1} \ > \ 2 \mu_1 \sqrt{\frac{\zeta (1-\mu_1)^2 \zeta}{4}} = \mu_1 \zeta(1-\mu_1),
\end{align}}
which implies $\omega < 1$. Similarly, whenever $\mu_1 < 1$ we have
\begin{align}
    2\gamma - \zeta - \frac{1}{\mu_1} < 2\sqrt{\eta \mu_1 \zeta} - \zeta - 1 \leq  2 \sqrt{\zeta} - \zeta - 1  = -(\sqrt{\zeta}-1)^2 \leq 0,
\end{align}
which implies $\omega > 0$. Also, using the parameter choice of \eqref{eq:askm12}, we have, 
\begin{align}
\label{eq:askm13}
  1-\omega+ \frac{\gamma^2 }{\eta}  - 2\gamma \mu_1 - & \frac{\gamma \omega \mu_1 (1-\alpha)}{\alpha}   = 1-\omega + \zeta \mu_1 - 2\gamma \mu_1 - \omega \zeta \mu_1^2 \nonumber \\
  & = 1- 2\gamma \mu_1 + \zeta \mu_1 - \omega (1+\zeta \mu_1^2) = 0.
\end{align}
Now, using all of the above relations (equation \eqref{eq:askm12}, \eqref{eq:askm13}) in equation \eqref{eq:askm11}, we get the following:
\begin{align}
\label{eq:askm14}
  \E \left[ d(v_{k+1},P)^2 + \frac{\gamma^2}{\eta} d(x_{k+1},P)^2 \right] &  \ \leq \ \omega \Big[d(v_{k},P)^2 + \underbrace{\frac{\gamma (1-\alpha)}{\alpha}}_{= \frac{\gamma^2 }{\eta}} \ d(x_{k},P)^2 \Big] \nonumber \\
  + & \Big [\underbrace{1-\omega+ \frac{\gamma^2 }{\eta}  - 2\gamma \mu_1 - \frac{\gamma \omega \mu_1 (1-\alpha)}{\alpha} }_{= \ 0} \Big ] \ d(y_k,P)^2 \nonumber \\
   = \ & \omega \ \left[ d(v_{k},P)^2+ \frac{\gamma^2}{\eta} d(x_{k},P)^2 \right].
\end{align}
Finally, taking expectation again with tower rule and substituting $\frac{\gamma^2 }{\eta}= \zeta \mu_1 $ we have,
\begin{align*}
\E \left[ d(v_{k+1},P)^2 + \zeta \mu_1 \  d(x_{k+1},P)^2 \right] \ & \leq \ \omega^{k+1} \ \E \left[ d(v_{0},P)^2+ \zeta \mu_1 \  d(x_{0},P)^2 \right] \\
& = (1+\zeta \mu_1) \ \omega^{k+1} \ d(x_{0},P)^2.
\end{align*}
This proves the Theorem. Furthermore, for faster convergence, we need to choose parameters such that, $\omega$ becomes as small as possible. In the proof, we assumed $\mu_1 <1$ holds which is the most probable scenario. Whenever $\mu_1 = 1$, we must have $\mu_1 = \mu_2 =1$ and Lemma 6 holds with both equality, i.e., $f(x) = d(x,P)^2$. Therefore if we choose $\alpha, \gamma, \omega$ as $\alpha = \frac{\eta}{\eta+ \gamma}, \ \gamma = \sqrt{\zeta \eta}, \ \omega = 1- \frac{2\gamma}{1+\zeta}$, we can check that condition (89) holds and $0 < \omega < 1$ holds for any $\zeta > 0$.

\paragraph{Proof of Theorem \ref{th:cesaro2}}

For any natural number $l \geq 1$, using the update formula of $v_{l+1}$, we have
\begin{align}
    \label{ces:20}
    v_{l+1} & = \omega v_l +(1-\omega) y_l - \gamma (a_{i^*}^Ty_{l}-b_{i^*})^+a_{i^*} \nonumber \\
    & \overset{\eqref{eq:yk}}{=} \ \left(1-\omega + \frac{\omega}{\alpha}\right) y_l - \frac{\omega (1-\alpha)}{\alpha} x_l - \gamma  (a_{i^*}^Ty_{l}-b_{i^*})^+a_{i^*}.
\end{align}
Let $\varphi = \omega (1-\alpha)$. It can be easily checked that $0 \leq \varphi < 1$. Now, considering equation \eqref{eq:yk}, we have
\begin{align}
    \label{ces:21}
    & y_{l+1}  = \alpha v_{l+1} + (1-\alpha) x_{l+1} \nonumber \\
    & \overset{\eqref{eq:askm1} \ \& \ \eqref{ces:20} }{=} \ (1+\omega -\alpha \omega) y_l - \omega (1-\alpha) y_{l-1} + \omega \delta (1-\alpha) (a_{j^*}^Ty_{l-1}-b_{j^*})^+a_{j^*}  \nonumber \\
    & \quad \quad \quad \quad \quad \quad - [\alpha \gamma+ (1-\alpha) \delta] \ (a_{i^*}^Ty_{l}-b_{i^*})^+a_{i^*} \nonumber \\
    & = (1+\omega -\alpha \omega) y_l - \omega (1-\alpha) y_{l-1} + \omega \delta (1-\alpha) (a_{j^*}^Ty_{l-1}-b_{j^*})^+a_{j^*}  \nonumber \\
    & \quad \quad \quad \quad \quad \quad - \delta (1+\omega - \alpha \omega ) \ (a_{i^*}^Ty_{l}-b_{i^*})^+a_{i^*} \nonumber \\
    & = (1+\varphi) y_l - \varphi y_{l-1} +  \delta \varphi  (a_{j^*}^Ty_{l-1}-b_{j^*})^+a_{j^*} - \delta (1+\varphi) \ (a_{i^*}^Ty_{l}-b_{i^*})^+a_{i^*}.
\end{align}
here, the index $i^{*}$ and $j^*$ are defined based on \eqref{def:i1} respectively for the sequences $y_l$ and $y_{l-1}$. Furthermore, with the choice of $x_0 = v_0$, the points $y_0$ and $y_1$ generated by the PASKM method (i.e, algorithm \ref{alg:paskm} with arbitrary parameter choice) can be calculated as 
\begin{align}
\label{ces:100}
    y_1  = \alpha v_1 + (1-\alpha) x_1 & = x_0 - (\alpha \gamma + \delta (1-\alpha)) (a_{i^*}^Ty_{0}-b_{i^*})^+a_{i^*} \nonumber \\
    & = y_0 - \delta (1+\varphi) (a_{i^*}^Ty_{0}-b_{i^*})^+a_{i^*},
\end{align}
since $y_0 = x_0 = v_0$. Now, let's define, $\bar{\vartheta}_l = \frac{\varphi }{1-\varphi}[y_{l}-y_{l-1}+\delta (a_{j^*}^Ty_{l-1}-b_{j^*})^+a_{j^*}]$, $ \bar{\Delta}_l = y_l + \bar{\vartheta}_l$ and $ \bar{\chi}_l = \|y_l+\bar{\vartheta}_l-\mathcal{P}(\bar{\Delta}_l)\|^2$, then using the update formula \eqref{ces:21}, we have
\begin{align*}
    y_{l+1} + \bar{\vartheta}_{l+1} & = y_{l+1} +\frac{\varphi }{1-\varphi}[y_{l+1}-y_{l}+\delta (a_{i^*}^Ty_{l}-b_{i^*})^+a_{i^*}]  \\
    & = \frac{1}{1-\varphi} y_{l+1} - \frac{\varphi}{1-\varphi} y_{l}+ \frac{\delta \varphi }{1-\varphi} (a_{i^*}^Ty_{l}-b_{i^*})^+a_{i^*} \\
    & \overset{\eqref{ces:21}}{=} y_l + \bar{\vartheta}_l - \frac{\delta}{1-\varphi} \left(a_{i^*}^Ty_l-b_{i^*}\right)^+ a_{i^*}.
\end{align*}
Using the above relation, we can write
\begin{align}
\label{ces:22}
    \bar{\chi}_{l+1} & = \|y_{l+1}+\bar{\vartheta}_{l+1}-\mathcal{P}(\bar{\Delta}_{l+1})\|^2 \overset{\text{Lemma} \ \ref{lem:distance}}{ \leq}  \|y_{l+1}+\bar{\vartheta}_{l+1}-\mathcal{P}(\bar{\Delta}_l)\|^2 \nonumber \\
    & = \big \| y_l + \bar{\vartheta}_l - \frac{\delta}{1-\varphi} \left(a_{i^*}^Ty_l-b_{i^*}\right)^+ a_{i^*}-\mathcal{P}(\bar{\Delta}_l) \big \|^2 \nonumber \\
    & =   \underbrace{\|y_l+\bar{\vartheta}_l-\mathcal{P}(\bar{\Delta}_l)\|^2}_{= \bar{\chi}_l} + \frac{\delta^2}{(1-\varphi)^2} \underbrace{\|(a_{i^*}^Ty_l- b_{i^*})^+a_{i^*}\|^2}_{I_1} \nonumber \\
    & - \frac{2 \delta}{1-\varphi}   \underbrace{\big \langle y_l+\bar{\vartheta}_l-\mathcal{P}(\bar{\Delta}_l) \ ,\  a_{i^*} (a_{i^*}^Ty_l-b_{i^{*}})^{+} \big \rangle }_{	I_2} \nonumber \\
    & = \bar{\chi}_l + \frac{\delta^2}{(1-\varphi)^2} 	I_1- \frac{2 \delta}{1-\varphi} 	I_2.
\end{align}
Taking expectation with respect to $\mathbb{S}_l$ we have,
\begin{align}
  \label{ces:23}
  \frac{\delta^2}{(1-\varphi)^2} \E_{\mathbb{S}_l} [I_1] \overset{\eqref{def:function} }{=}  \frac{2\delta^2}{(1-\varphi)^2} f(y_l).
\end{align}
Similarly, we can simplify the third term of \eqref{ces:22} as
\begin{align}
  \label{ces:24}
  & - \frac{2\delta}{1-\varphi}  \E_{\mathbb{S}_l} [I_2] \nonumber \\
  & \overset{\eqref{def:function} }{=}  -\frac{2\delta}{1-\varphi} \big \langle y_l-\mathcal{P}(\bar{\Delta}_l)  , \nabla f(y_l) \big \rangle + \frac{2\delta \varphi}{(1-\varphi)^2} \big \langle y_{l-1}-y_l - \delta \nabla f(y_{l-1}) ,  \nabla f(y_l) \big \rangle \nonumber \\
  & =  -\frac{2\delta}{1-\varphi} \big \langle y_l-\mathcal{P}(\bar{\Delta}_l)  , \nabla f(y_l) \big \rangle + \frac{2\delta \varphi}{(1-\varphi)^2} \big \langle y_{l-1}-y_l ,  \nabla f(y_l) \big \rangle \nonumber \\
  &  \quad \quad \quad \quad \quad  -\frac{\delta^2 \varphi}{(1-\varphi)^2} \left[ \|\nabla f(y_l)+\nabla f(y_{l-1})\|^2-\|\nabla f(y_l)\|^2-\|\nabla f(y_{l-1})\|^2 \right] \nonumber \\
  & \overset{\text{Lemma} \ \ref{lem:grad} \ \& \ \ref{lem:grad1}}{\leq}  - \frac{4\delta}{1-\varphi} f(y_l) + \frac{2\delta \varphi}{(1-\varphi)^2} \left[f(y_{l-1}) -  f(y_{l}) \right] + \frac{2\delta^2 \varphi}{(1-\varphi)^2} \left[f(y_{l-1}) +  f(y_{l}) \right] \nonumber \\ 
  & = \frac{2 \delta \varphi (1+\delta)}{(1-\varphi)^2} f(y_{l-1}) -\frac{2 \delta \varphi (1+\delta)}{(1-\varphi)^2} f(y_{l}) + \frac{4 \delta  (\varphi + \delta \varphi-1)}{(1-\varphi)^2} f(y_{l}).
\end{align} 
Using the expressions of equation \eqref{ces:23} and \eqref{ces:24} in \eqref{ces:22} and simplifying further, we have
\begin{align}
    \label{ces:25}
    \E[\bar{\chi}_{l+1}] + \frac{2\delta \varphi (1+\delta)}{(1-\varphi)^2} f(y_l) + \varsigma f(y_l) \ \leq \ \E[\bar{\chi}_l] + \frac{2\delta \varphi (1+\delta)}{(1-\varphi)^2} f(y_{l-1}),  
\end{align}
here,
\begin{align}
 \label{ces:26}
   \varsigma =   \frac{4 \delta  (1-\varphi - \delta \varphi)}{(1-\varphi)^2}  -  \frac{2 \delta^2  }{(1-\varphi)^2} =  \frac{2 \delta  (2-2 \varphi - 2\delta \varphi -\delta)}{(1-\varphi)^2} \ > \ 0.
\end{align}
Now, taking expectation again in \eqref{ces:25} and using the tower property, we get,
\begin{align}
\label{ces:27}
    \bar{q}_{l+1} + \varsigma \E[f(y_l)] \leq \bar{q}_l, \quad l = 1,2,3...,
\end{align}
where, $\bar{q}_l = \E[\bar{\chi}_l] + \frac{2\delta \varphi (1+\delta)}{(1-\varphi)^2} \E[f(y_{l-1})]$. Summing up \eqref{ces:27} for $l=1,2,...,k$ we get
\begin{align}
    \label{ces:28}
    \sum \limits_{l=1}^{k} \E [f(y_l)] \ \leq \ \frac{\bar{q}_1-\bar{q}_{k+1}}{\varsigma} \ \leq \ \frac{\bar{q}_1}{\varsigma}.
\end{align}
Now, using Jensen's inequality, we have
\begin{align*}
    \E \left[f(\bar{y_k})\right] = \E \left[f\left(\sum \limits_{l=1}^{k} \frac{y_k}{k}\right)\right] \ \leq \ \E \left[\frac{1}{k} \sum \limits_{l=1}^{k}f(y_l)\right] \ = \ \frac{1}{k}  \sum \limits_{l=1}^{k} \E [f(y_l)] \ \overset{\eqref{ces:28}}{\leq}  \frac{\bar{q}_1}{\varsigma k}.
\end{align*}
From \eqref{ces:100}, $y_1 = y_0 - \delta (1+\varphi) (a_{i^*}^Ty_0-b_{i^*})^+ a_{i^*}$ and $\bar{\vartheta}_1 = \frac{-\varphi^2 \delta }{1-\varphi} (a_{i^*}^Ty_0-b_{i^*})^+ a_{i^*} $. Then, 
\begin{align}
    \label{ces:29}
    \E[\bar{\chi}_1] & = \E \left[\|y_1+\bar{\vartheta}_1 -\mathcal{P}(\bar{\Delta}_1)\|^2 \right] \overset{\text{Lemma} \ \ref{lem:distance}}{ \leq} \E \left[\|y_1+\bar{\vartheta}_1 -\mathcal{P}(y_0)\|^2 \right] \nonumber \\
    & = \E \left[\|y_0 -\mathcal{P}(y_0)-\frac{\delta }{1-\varphi} (a_{i^*}^Ty_0-b_{i^*})^+ a_{i^*} \|^2 \right] \nonumber \\
    & = \|y_0 -\mathcal{P}(y_0)\|^2 + \frac{\delta^2}{(1-\varphi)^2} \E [|(a_{i^*}^Ty_0-b_{i^*})^+|^2] \nonumber \\
    & \quad  \quad \quad \quad  \quad \quad - \frac{2\delta }{1-\varphi} \langle y_0 -\mathcal{P}(y_0) , \E [(a_{i^*}^Ty_0-b_{i^*})^+ a_{i^*}] \rangle \nonumber \\
    & \overset{\text{Lemma} \ \ref{lem:grad1} }{\leq} \|y_0 -\mathcal{P}(y_0)\|^2 + \frac{2\delta^2}{(1-\varphi)^2} f(y_0)-\frac{4\delta }{1-\varphi} f(y_0).
\end{align}
Now, from our construction we get
\begin{align*}
   \bar{q}_1  = \E[\bar{\chi}_1] + \frac{2\delta \varphi (1+\delta)}{(1-\varphi)^2} \E[f(y_{0})] \leq   d(y_0,P)^2+ \frac{2\delta (\delta -2+ 3\varphi + \delta \varphi)}{(1-\varphi)^2} f(y_0).
\end{align*}
Substituting the values of $\varsigma$ and $q_1$ in the expression of $\E \left[f(\bar{y_k})\right] $, we have the following:
\begin{align*}
    \E \left[f(\bar{y}_k)\right] \leq \frac{ (1-\omega + \alpha \omega )^2 \ d(y_0,P)^2+ 2\delta (\delta -2+ 3 \omega - 3 \alpha \omega + \delta \omega - \delta \alpha \omega) f(y_0)}{2 \delta k \left(2-2 \omega + 2 \alpha \omega - 2\delta \omega + 2 \delta \alpha \omega -\delta\right)}.
\end{align*}


\bibliographystyle{unsrt}
\bibliography{template}

\begin{thebibliography}{10}

\bibitem{strohmer:2008}
Thomas Strohmer and Roman Vershynin.
\newblock A randomized kaczmarz algorithm with exponential convergence.
\newblock {\em Journal of Fourier Analysis and Applications}, 15(2):262, Apr
  2008.

\bibitem{lewis:2010}
Dennis Leventhal and Adrian~S. Lewis.
\newblock Randomized methods for linear constraints: Convergence rates and
  conditioning.
\newblock {\em Mathematics of Operations Research}, 35(3):641--654, 2010.

\bibitem{needell:2010}
Deanna Needell.
\newblock Randomized kaczmarz solver for noisy linear systems.
\newblock {\em BIT Numerical Mathematics}, 50(2):395--403, Jun 2010.

\bibitem{drineas:2011}
Petros Drineas, Michael~W. Mahoney, Shan Muthukrishnan, and Tam{\'a}s
  Sarl{\'o}s.
\newblock Faster least squares approximation.
\newblock {\em Numerische Mathematik}, 117(2):219--249, Feb 2011.

\bibitem{zouzias:2013}
Anastasios Zouzias and Nikolaos~M. Freris.
\newblock Randomized extended kaczmarz for solving least squares.
\newblock {\em SIAM Journal on Matrix Analysis and Applications},
  34(2):773--793, 2013.

\bibitem{lee:2013}
Yin~Tat Lee and Aaron Sidford.
\newblock Efficient accelerated coordinate descent methods and faster
  algorithms for solving linear systems.
\newblock In {\em Proceedings of the 2013 IEEE 54th Annual Symposium on
  Foundations of Computer Science}, FOCS '13, pages 147--156, Washington, DC,
  USA, 2013. IEEE Computer Society.

\bibitem{ma:2015}
Anna Ma, Deanna Needell, and Aaditya Ramdas.
\newblock Convergence properties of the randomized extended gauss seidel and
  kaczmarz methods.
\newblock {\em {SIAM} Journal on Matrix Analysis and Applications},
  36(4):1590--1604, Jan 2015.

\bibitem{gower:2015}
Robert~M. Gower and Peter Richtárik.
\newblock Randomized iterative methods for linear systems.
\newblock {\em SIAM Journal on Matrix Analysis and Applications},
  36(4):1660--1690, 2015.

\bibitem{qu:2016}
Zheng Qu, Peter Richtarik, Martin Takac, and Olivier Fercoq.
\newblock {SDNA: Stochastic Dual Newton Ascent for Empirical Risk
  Minimization}.
\newblock In {\em Proceedings of The 33rd International Conference on Machine
  Learning}, volume~48, pages 1823--1832, New York, USA, 20--22 Jun 2016. PMLR.

\bibitem{needell:2016}
Deanna Needell, Nathan Srebro, and Rachel Ward.
\newblock Stochastic gradient descent, weighted sampling, and the randomized
  kaczmarz algorithm.
\newblock {\em Mathematical Programming}, 155(1):549--573, Jan 2016.

\bibitem{haddock:2017}
Jes{\'u}s De~Loera, Jamie Haddock, and Deanna Needell.
\newblock A sampling kaczmarz--motzkin algorithm for linear feasibility.
\newblock {\em SIAM Journal on Scientific Computing}, 39(5):S66--S87, 2017.

\bibitem{razaviyayn:2019}
Meisam Razaviyayn, Mingyi Hong, Navid Reyhanian, and Zhi-Quan Luo.
\newblock A linearly convergent doubly stochastic gauss--seidel algorithm for
  solving linear equations and a certain class of over-parameterized
  optimization problems.
\newblock {\em Mathematical Programming}, 176(1):465--496, Jul 2019.

\bibitem{kaczmarz:1937}
Stefan Kaczmarz.
\newblock Angenaherte auflsung von systemen linearer gleichungen.
\newblock {\em Bulletin International de l'Acadmie Polonaise des Sciences et
  des Letters}, 35:355--357, 1937.

\bibitem{gordon:1970}
Richard Gordon, Robert Bender, and Gabor~T. Herman.
\newblock Algebraic reconstruction techniques (art) for three-dimensional
  electron microscopy and x-ray photography.
\newblock {\em Journal of Theoretical Biology}, 29(3):471 -- 481, 1970.

\bibitem{Censor:1988}
Yair Censor.
\newblock Parallel application of block-iterative methods in medical imaging
  and radiation therapy.
\newblock {\em Mathematical Programming}, 42(1):307--325, Apr 1988.

\bibitem{herman:2009}
Gabor~T. Herman.
\newblock {\em Fundamentals of Computerized Tomography: Image Reconstruction
  from Projections}.
\newblock Springer Publishing Company, Incorporated, 2nd edition, 2009.

\bibitem{lorenz:2014}
D.~A. {Lorenz}, S.~{Wenger}, F.~{Schöpfer}, and M.~{Magnor}.
\newblock A sparse kaczmarz solver and a linearized bregman method for online
  compressed sensing.
\newblock In {\em 2014 IEEE International Conference on Image Processing
  (ICIP)}, pages 1347--1351, Oct 2014.

\bibitem{elble:2010}
Joseph~M. Elble, Nikolaos~V. Sahinidis, and Panagiotis Vouzis.
\newblock Gpu computing with kaczmarz’s and other iterative algorithms for
  linear systems.
\newblock {\em Parallel Computing}, 36(5):215 -- 231, 2010.
\newblock Parallel Matrix Algorithms and Applications.

\bibitem{fabio:2012}
Fabio Pasqualetti, Ruggero Carli, and Francesco Bullo.
\newblock Distributed estimation via iterative projections with application to
  power network monitoring.
\newblock {\em Automatica}, 48(5):747 -- 758, 2012.

\bibitem{censor:1981}
Yair Censor.
\newblock Row-action methods for huge and sparse systems and their
  applications.
\newblock {\em SIAM Review}, 23(4):444--466, 1981.

\bibitem{agamon}
Shmuel Agamon.
\newblock The relaxation method for linear inequalities.
\newblock {\em Canadian J. Math}, pages 382--392, 1954.

\bibitem{motzkin}
Theodore~S. Motzkin and Issac~J. Schoenberg.
\newblock The relaxation method for linear inequalities.
\newblock {\em Canadian J. Math}, pages 393--404, 1954.

\bibitem{rosenblatt}
Frank Rosenblatt.
\newblock The perceptron: A probabilistic model for information storage and
  organization in the brain.
\newblock {\em Psychological Review}, pages 65--386, 1958.

\bibitem{ramdas:2014}
Aaditya Ramdas and Javier Peña.
\newblock Margins, kernels and non-linear smoothed perceptrons.
\newblock In {\em Proceedings of the 31st International Conference on Machine
  Learning}, volume~32, pages 244--252, Bejing, China, 22--24 Jun 2014. PMLR.

\bibitem{ramdas:2016}
Aaditya Ramdas and Javier Pe\~{n}a.
\newblock Towards a deeper geometric, analytic and algorithmic understanding of
  margins.
\newblock {\em Optimization Methods and Software}, 31(2):377--391, 2016.

\bibitem{nutini:2016}
Julie Nutini, Behrooz Sepehry, Issam Laradji, Mark Schmidt, Hoyt Koepke, and
  Alim Virani.
\newblock Convergence rates for greedy kaczmarz algorithms, and faster
  randomized kaczmarz rules using the orthogonality graph.
\newblock In {\em Proceedings of the Thirty-Second Conference on Uncertainty in
  Artificial Intelligence}, UAI'16, pages 547--556, Arlington, Virginia, United
  States, 2016. AUAI Press.

\bibitem{petra:2016}
Stefania Petra and Constantin Popa.
\newblock Single projection kaczmarz extended algorithms.
\newblock {\em Numerical Algorithms}, 73(3):791--806, Nov 2016.

\bibitem{telgen:1982}
Jan Telgen.
\newblock On relaxation methods for systems of linear inequalities.
\newblock {\em European Journal of Operational Research}, 9(2):184 -- 189,
  1982.

\bibitem{Maurras1981}
J.~F. Maurras, K.~Truemper, and M.~Akg{\"u}l.
\newblock Polynomial algorithms for a class of linear programs.
\newblock {\em Mathematical Programming}, 21(1):121--136, Dec 1981.

\bibitem{Chubanov:2012}
Sergei Chubanov.
\newblock A strongly polynomial algorithm for linear systems having a binary
  solution.
\newblock {\em Mathematical Programming}, 134(2):533--570, Sep 2012.

\bibitem{Chubanov:2015}
Sergei Chubanov.
\newblock A polynomial projection algorithm for linear feasibility problems.
\newblock {\em Mathematical Programming}, 153(2):687--713, Nov 2015.

\bibitem{wright:2016}
Ji~Liu and Stephen~J. Wright.
\newblock An accelerated randomized kaczmarz algorithm.
\newblock {\em Math. Comput.}, 85(297):153--178, 2016.

\bibitem{NEEDELL:2015}
Deanna Needell, Ran Zhao, and Anastasios Zouzias.
\newblock Randomized block kaczmarz method with projection for solving least
  squares.
\newblock {\em Linear Algebra and its Applications}, 484:322 -- 343, 2015.

\bibitem{anna:2015}
Anna. Ma, Deanna. Needell, and Aaditya. Ramdas.
\newblock Convergence properties of the randomized extended gauss--seidel and
  kaczmarz methods.
\newblock {\em SIAM Journal on Matrix Analysis and Applications},
  36(4):1590--1604, 2015.

\bibitem{hefny:2017}
Ahmed. Hefny, Deanna. Needell, and Aaditya. Ramdas.
\newblock Rows versus columns: Randomized kaczmarz or gauss--seidel for ridge
  regression.
\newblock {\em SIAM Journal on Scientific Computing}, 39(5):S528--S542, 2017.

\bibitem{gower2016linearly}
Robert~M. Gower and Peter Richtárik.
\newblock Linearly convergent randomized iterative methods for computing the
  pseudoinverse, 2016.

\bibitem{gower:2017}
Robert~M. Gower and Peter. Richtárik.
\newblock Randomized quasi-newton updates are linearly convergent matrix
  inversion algorithms.
\newblock {\em SIAM Journal on Matrix Analysis and Applications},
  38(4):1380--1409, 2017.

\bibitem{richtrik2017stochastic}
Peter Richtárik and Martin Takáč.
\newblock Stochastic reformulations of linear systems: Algorithms and
  convergence theory, 2017.

\bibitem{NIPS:2018}
Robert Gower, Filip Hanzely, Peter Richtarik, and Sebastian~U Stich.
\newblock Accelerated stochastic matrix inversion: General theory and speeding
  up bfgs rules for faster second-order optimization.
\newblock In S.~Bengio, H.~Wallach, H.~Larochelle, K.~Grauman, N.~Cesa-Bianchi,
  and R.~Garnett, editors, {\em Advances in Neural Information Processing
  Systems 31}, pages 1619--1629. Curran Associates, Inc., 2018.

\bibitem{NEEDELL:2014}
Deanna Needell and Joel~A. Tropp.
\newblock Paved with good intentions: Analysis of a randomized block kaczmarz
  method.
\newblock {\em Linear Algebra and its Applications}, 441:199 -- 221, 2014.
\newblock Special Issue on Sparse Approximate Solution of Linear Systems.

\bibitem{blockneddel:2015}
Jonathan Briskman and Deanna Needell.
\newblock Block kaczmarz method with inequalities.
\newblock {\em J. Math. Imaging Vis.}, 52(3):385–396, July 2015.

\bibitem{needell2019block}
Deanna Needell and Elizaveta Rebrova.
\newblock On block gaussian sketching for the kaczmarz method, 2019.

\bibitem{basu:2014}
Amitabh Basu, Jesús A.~De Loera, and Mark Junod.
\newblock On chubanov's method for linear programming.
\newblock {\em INFORMS Journal on Computing}, 26(2):336--350, 2014.

\bibitem{VEGH:2014}
László~A. Végh and Giacomo Zambelli.
\newblock A polynomial projection-type algorithm for linear programming.
\newblock {\em Operations Research Letters}, 42(1):91 -- 96, 2014.

\bibitem{eldar:2011}
Yonina~C. Eldar and Deanna Needell.
\newblock Acceleration of randomized kaczmarz method via the
  johnson--lindenstrauss lemma.
\newblock {\em Numerical Algorithms}, 58(2):163--177, Oct 2011.

\bibitem{agaskar:2014}
A.~{Agaskar}, C.~{Wang}, and Y.~M. {Lu}.
\newblock Randomized kaczmarz algorithms: Exact mse analysis and optimal
  sampling probabilities.
\newblock In {\em 2014 IEEE Global Conference on Signal and Information
  Processing (GlobalSIP)}, pages 389--393, Dec 2014.

\bibitem{bai:2018}
Zhong-Zhi Bai and Wen-Ting Wu.
\newblock On relaxed greedy randomized kaczmarz methods for solving large
  sparse linear systems.
\newblock {\em Applied Mathematics Letters}, 83:21 -- 26, 2018.

\bibitem{greedbai:2018}
Zhong-Zhi. Bai and Wen-Ting. Wu.
\newblock On greedy randomized kaczmarz method for solving large sparse linear
  systems.
\newblock {\em SIAM Journal on Scientific Computing}, 40(1):A592--A606, 2018.

\bibitem{kovachki2019analysis}
Nikola~B Kovachki and Andrew~M Stuart.
\newblock Analysis of momentum methods.
\newblock {\em arXiv preprint arXiv:1906.04285}, 2019.

\bibitem{ruder2016overview}
Sebastian Ruder.
\newblock An overview of gradient descent optimization algorithms.
\newblock {\em arXiv preprint arXiv:1609.04747}, 2016.

\bibitem{polyak1964some}
Boris~T Polyak.
\newblock Some methods of speeding up the convergence of iteration methods.
\newblock {\em USSR Computational Mathematics and Mathematical Physics},
  4(5):1--17, 1964.

\bibitem{nesterov:1983}
Yuri Nesterov.
\newblock A method for solving the convex programming problem with convergence
  rate $o(1/k\sp{2})$.
\newblock {\em Soviet Mathematics Doklady}, Vol. 27:p(372--376), 1983.

\bibitem{nesterov:2005}
Yuri Nesterov.
\newblock Smooth minimization of non-smooth functions.
\newblock {\em Mathematical Programming}, 103(1):127--152, May 2005.

\bibitem{nesterov:2013}
Yuri Nesterov.
\newblock Gradient methods for minimizing composite functions.
\newblock {\em Mathematical Programming}, 140(1):125--161, Aug 2013.

\bibitem{nesterov:2014}
Yuri Nesterov.
\newblock {\em Introductory Lectures on Convex Optimization: A Basic Course}.
\newblock Springer Publishing Company, Incorporated, 1 edition, 2014.

\bibitem{nesterov:2012}
Yuri Nesterov.
\newblock Efficiency of coordinate descent methods on huge-scale optimization
  problems.
\newblock {\em SIAM Journal on Optimization}, 22(2):341--362, 2012.

\bibitem{loizou:2017}
Nicolas Loizou and Peter Richtárik.
\newblock Momentum and stochastic momentum for stochastic gradient, newton,
  proximal point and subspace descent methods, 2017.

\bibitem{morshed:2018}
Md~Sarowar Morshed and Md. Noor-E-Alam.
\newblock Generalized affine scaling algorithms for linear programming
  problems.
\newblock {\em Computers \& Operations Research}, 114:104807, 2020.

\bibitem{peter:2019}
Michael~Rabbat Nicolas~Loizou and Peter Richtárik.
\newblock Provably accelerated randomized gossip algorithms.
\newblock {\em Arxiv}, 2018.

\bibitem{Morshed2019}
Md~Sarowar Morshed, Md~Saiful Islam, and Md. Noor-E-Alam.
\newblock Accelerated sampling kaczmarz motzkin algorithm for the linear
  feasibility problem.
\newblock {\em Journal of Global Optimization}, Oct 2019.

\bibitem{haddock:2019}
Jamie Haddock and Anna Ma.
\newblock Greed works: An improved analysis of sampling kaczmarz-motkzin, 2019.

\bibitem{morshed:momentum}
Md~Sarowar Morshed and Md. Noor-E-Alam.
\newblock Heavy ball momentum induced sampling kaczmarz motzkin methods for
  linear feasibility problems.
\newblock {\em arXiv preprint arXiv:200908251}, 2020.

\bibitem{hoffman}
Alan~J Hoffman.
\newblock On approximate solutions of systems of linear inequalities.
\newblock In {\em Selected Papers Of Alan J Hoffman: With Commentary}, pages
  174--176. World Scientific, 2003.

\bibitem{karimi:2016}
Hamed Karimi, Julie Nutini, and Mark Schmidt.
\newblock Linear convergence of gradient and proximal-gradient methods under
  the polyak-{\l}ojasiewicz condition.
\newblock In Paolo Frasconi, Niels Landwehr, Giuseppe Manco, and Jilles
  Vreeken, editors, {\em Machine Learning and Knowledge Discovery in
  Databases}, pages 795--811, Cham, 2016. Springer International Publishing.

\bibitem{KHACHIYAN:1980}
L.G. Khachiyan.
\newblock Polynomial algorithms in linear programming.
\newblock {\em USSR Computational Mathematics and Mathematical Physics},
  20(1):53 -- 72, 1980.

\bibitem{yeh:2009}
I-Cheng Yeh and Che-hui Lien.
\newblock The comparisons of data mining techniques for the predictive accuracy
  of probability of default of credit card clients.
\newblock {\em Expert Syst. Appl.}, 36(2):2473--2480, Mar 2009.

\bibitem{lichman:2013}
Moshe Lichman.
\newblock {UCI} machine learning repository, 2013.

\bibitem{netlib}
{Netlib}.
\newblock The netlib linear programming library.

\bibitem{calafiore:2014}
Giuseppe Calafiore and Laurent {El Ghaoui}.
\newblock {\em Optimization Models}.
\newblock Control systems and optimization series. Cambridge University Press,
  October 2014.

\end{thebibliography}


\end{document}